\documentclass[12pt,english]{article}
\usepackage[T1]{fontenc}
\usepackage[latin9]{inputenc}
\usepackage{amsmath}
\usepackage{amssymb}
\usepackage{graphicx}
\usepackage{setspace}
\usepackage{esint}
\makeatletter
\providecommand{\tabularnewline}{\\}

\usepackage{ifpdf}

\newif\ifpdf
\ifx\pdfoutput\undefined
  \pdffalse
\else
  \pdfoutput=1
  \pdftrue
\fi

\RequirePackage{xspace} %
\RequirePackage{subfigure} %

\ifpdf
  \RequirePackage[ pdftex, plainpages = false, pdfpagelabels,
                 pdfpagelayout = useoutlines,
                 bookmarks,
                 breaklinks = true,
                 linktocpage,
                 colorlinks = true,
                 linkcolor = blue,
                 urlcolor  = blue,
                 citecolor = blue,
                 anchorcolor = blue,
                 hyperindex = true,
                 hyperfigures
                 ]{hyperref}

\else
  \RequirePackage{color}
  \RequirePackage{colortbl}
   \RequirePackage{array}
  \RequirePackage[dvips]{graphicx}
  \RequirePackage{hyperref}
  \usepackage{rotating}
\fi

\usepackage{makeidx} 
\usepackage{setspace} 
\usepackage{rotating} 
\usepackage{ecltree}
\usepackage{epic}
\usepackage{supertabular}  
\usepackage{color}
\usepackage{exscale}
\usepackage{fontenc}
\usepackage{ifthen}
\usepackage{latexsym}
\usepackage{makeidx}
\usepackage{syntonly}
\usepackage{inputenc}
\usepackage{graphicx}
\usepackage{setspace}
\usepackage{caption2}
\usepackage[english]{babel}
\usepackage[square, comma,numbers,sort&compress]{natbib}
\usepackage{hypernat}
\usepackage{boxedminipage}
\usepackage{framed}
\usepackage{longtable}
\usepackage[all]{hypcap}

\newcommand{\note}[1]{\marginpar[left]{\singlespace \tiny #1}}

\renewcommand{\sectionmark}[1]%
      {\markright{\thesection\ #1}}

\renewcommand{\note}[1]{}
\renewcommand\@biblabel[1]{$[{#1}]$}

\usepackage[bottom]{footmisc}

\makeatother

\usepackage{babel}
\begin{document}

\title{Principles of Differential Geometry\\
\vspace{6cm}}

\author{Taha Sochi\thanks{Department of Physics \& Astronomy, University College London, Gower
Street, London, WC1E 6BT. Email: t.sochi@ucl.ac.uk.}\vspace{8cm}}

\maketitle
\pagebreak{}

\phantomsection \addcontentsline{toc}{section}{Preface}

\section*{Preface\label{secPreface}}

The present text is a collection of notes about differential geometry
prepared to some extent as part of tutorials about topics and applications
related to tensor calculus. They can be regarded as continuation to
the previous notes on tensor calculus \cite{SochiTC1, SochiTC2} as
they are based on the materials and conventions given in those documents.
They can be used as a reference for a first course on the subject
or as part of a course on tensor calculus.

We generally follow the same notations and conventions employed in
\cite{SochiTC1,SochiTC2}, however the following points should be
observed:

$\bullet$ Following the convention of several authors, when discussing
issues related to 2D and 3D manifolds the Greek indices range over
1,2 while the Latin indices range over 1,2,3. Therefore, the use of
Greek and Latin indices indicates the type of the intended manifold
unless it is stated otherwise.

$\bullet$ The indexed $u$ with Greek letters are used for the surface
curvilinear coordinates while the indexed $x$ with Latin letters
are used largely for the space Cartesian coordinates but sometimes
they are used for the space curvilinear coordinates. Comments are
added when necessary to clarify the situation.

$\bullet$ For curvilinear coordinates of surfaces we use both $u,v$
and $u^{1},u^{2}$ as each has advantages in certain contexts; in
particular the latter is necessary for expressing the equations of
differential geometry in tensor forms.

$\bullet$ Unless stated otherwise, ``surface'' and ``space''
in the present notes mean 2D and 3D manifolds respectively.

$\bullet$ The indexed $\mathbf{E}$ are largely used for the surface,
rather than space, basis vectors where they are labeled with Greek
indices, e.g. $\mathbf{E}_{\alpha}$ and $\mathbf{E}^{\beta}$. However,
in a few cases indexed $\mathbf{E}$ are also used for space basis
vectors in which case they are distinguished by using Latin indices,
e.g. $\mathbf{E}_{i}$ and $\mathbf{E}^{j}$. When the basis vectors
are numbered rather than indexed, the distinction should be obvious
from the context if it is not stated explicitly.

$\bullet$ The Christoffel symbols may be based on the space metric
or the surface metric, hence when a number of Christoffel symbols
in a certain context or equation are based on more than one metric,
the type of indices (Greek or Latin) can be used as an indicator to
the underlying metric where the Greek indices represent surface (e.g.
$[\alpha\beta,\gamma]$ and $\Gamma_{\alpha\beta}^{\gamma}$) while
the Latin indices represent space (e.g. $[ij,k]$ and $\Gamma_{ij}^{k}$).
Nevertheless, comments are generally added to remove any ambiguity.
In particular, when the Christoffel symbols are numbered (e.g. $\Gamma_{22}^{1}$)
comments will be added to clarify the situation.

$\bullet$ The present notes are largely based on curves and surfaces
embedded in a 3D flat space coordinated by a rectangular Cartesian
system ($x,y,z$).

$\bullet$ Some of the definitions provided in the present text which
are related to concepts from other subjects of mathematics such as
calculus and topology are elementary because of the limits on the
text size, moreover the notes are not prepared for these subjects.
The purpose of these definitions is to provide a basic understanding
of the related ideas in general. The readers are advised to refer
to textbooks on those subjects for more technical and detailed definitions.

$\bullet$ For brevity, convenience and clean notation in certain
contexts, we use overdot (e.g. $\mathbf{\dot{r}}$) to indicate derivative
with respect to a general parameter $t$ while we use prime (e.g.
$\mathbf{r}'$) to indicate derivative with respect to a natural parameter
$s$ representing arc length.

$\bullet$ The materials of differential geometry are strongly interlinked
and hence any elementary text about the subject, like the present
one, will face the problem of arranging the materials in a natural
order to ensure gradual development of concepts. In this text we largely
followed such a scheme; however this is not always possible and hence
in some cases references are provided for materials in later parts
of the text for concepts needed in earlier parts.

$\bullet$ To facilitate linking related concepts and theorems, and
hence ensuring a better understanding of the provided materials, we
use hyperlinks (which are colored blue) extensively in the text.
The reader, therefore, is advised to use these links.

$\bullet$ Twisted curves can reside in a 2D manifold (surface) or
in a higher dimensionality manifold (usually 3D space). Hence we use
``surface curves'' and ``space curves'' to refer to the type of
the manifold of residence. However, in most cases a single curve can
be viewed as a resident of more than one manifold and hence it is
a surface and space curve at the same time. For example, a curve embedded
in a surface which in its turn is embedded in a 3D space is a surface
curve and a space curve at the same time. Consequently, in the present
text these terms should be interpreted flexibly. Many statements formulated
in terms of a particular type of manifolds can be correctly and easily
extended to another type with minimal adjustments of dimensionality
and symbolism. Moreover, ``space'' in many statements should be
understood in its general meaning as a manifold embracing the curve
not as opposite to ``surface'' and hence it can include a 2D space,
i.e. surface.

$\bullet$ We deliberately use a variety of notations for the same
concepts (e.g. the above-mentioned $u^{1},u^{2}$ and $u,v$) for
convenience and to familiarize the reader with different notations
all of which are in common use in the literature of differential geometry
and tensor calculus. Having proficiency in these subjects requires
familiarity with these various, and sometimes conflicting, notations.

$\bullet$ Of particular importance is an issue related to the previous
point that is the use of different symbols for the coefficients of
the first and second fundamental forms $E,F,G,e,f,g$ and the coefficients
of the surface covariant metric and curvature tensors $a_{11},a_{12},a_{22},b_{11},b_{12},b_{22}$
despite the equivalence of these coefficients, i.e. $\left(E,F,G,e,f,g\right)=\left(a_{11},a_{12},a_{22},b_{11},b_{12},b_{22}\right)$
and hence all these different formulations can be replaced by just
one. However, we keep both notations as they are both in common use
in the literature of differential geometry and tensor calculus; moreover
in many situations the use of one of these notations or the other
is advantageous depending on the context.

\pagebreak{}

\phantomsection \addcontentsline{toc}{section}{Contents}

\tableofcontents{}

\pagebreak{}

\section{Preliminaries\label{subPreliminaries}}

\subsection{Differential Geometry}

$\bullet$ Differential geometry is a branch of mathematics that largely
employs methods and techniques of other branches of mathematics such
as differential and integral calculus, topology and tensor analysis
to investigate geometric issues related to abstract objects, mainly
space curves and surfaces, and their properties where these investigations
are mostly focused on these properties at small scales. The investigations
also include characterizing categories of these objects. There is
also a close link between differential geometry and the disciplines
of differential topology and differential equations.

$\bullet$ Differential geometry may be contrasted with ``Algebraic
geometry'' which is another branch of geometry that uses algebraic
tools to investigate geometric issues mainly of global nature.

$\bullet$ The investigation of the properties of curves and surfaces
in differential geometry are closely linked. For instance, investigating
the characteristics of space curves is largely exploited in the investigation
of surfaces since common properties of surfaces are defined and quantified
in terms of the properties of curves embedded in the surface. For
example, several aspects of the surface curvature at a point are defined
and quantified in terms of the parameters of the surface curves passing
through that point.

\subsection{Categories of Curve and Surface Properties}

$\bullet$ The properties of curves and surfaces may be categorized
into two main groups: local and global where these properties describe
the geometry of the curves and surfaces \textit{in the small} and
\textit{in the large} respectively. The local properties correspond
to the characteristics of the object in the immediate neighborhood
of a point on the object such as the curvature of a curve or surface
at that point, while the global properties correspond to the characteristics
of the object on a large scale and over extended parts of the object
such as the number of stationary points of a curve or a surface or
being a one-side surface like Mobius strip which is locally a double-sided
surface. As indicated earlier, differential geometry of space curves
and surfaces is mainly concerned with the local properties. The investigation
of global properties normally involve topological treatments.\footnote{There is a branch of differential geometry dedicated to the investigation
of global (or \textit{in the large}) properties. The focus of the
present text is largely differential geometry \textit{in the small}
although a number of global differential geometric issues are investigated
casually.}

$\bullet$ Another classification of the properties of curves and
surfaces, based on their relation to the embedding external space
in which they reside, may be made where the properties are divided
into intrinsic and extrinsic. The first category corresponds to those
properties which are independent in their existence and definition
from the ambient space which embraces the object such as the distance
along a given curve or the Gaussian curvature of a surface at a given
point (see $\S$ \ref{subGaussianCurvature}), while the second category
is related to those properties which depend in their existence and
definition on the external embedding space such as having a normal
vector at a point on the curve or surface. The idea of intrinsic and
extrinsic properties may be illustrated by an inhabitant of a surface
with a 2D perception (hereafter this creature will be called ``2D
inhabitant'') where he can detect and measure intrinsic properties
but not extrinsic properties as the former do not require appealing
to an external embedding 3D space in which the surface is immersed
while the latter do. Hence, in simple terms all the properties that
can be detected and measured by a 2D inhabitant are intrinsic to the
surface while the other properties are extrinsic. A 1D inhabitant
of a curve may also be used, to a lesser extent, analogously to distinguish
between intrinsic and extrinsic properties of space curves (refer
for example to $\S$ \ref{1Dinhabitant}).

$\bullet$ More technically, the intrinsic properties are defined
and expressed in terms of the metric tensor (formulated in differential
geometry as the first fundamental form; see $\S$ \ref{subFirstFundamentalForm})
while the extrinsic properties are expressed in terms of the surface
curvature tensor (formulated in differential geometry as the second
fundamental form; see $\S$ \ref{subSecondFundamentalForm}).

$\bullet$ The ``intrinsic geometry'' of the surface comprises the
collection of all the intrinsic properties of the surface.

$\bullet$ When two surfaces can have a coordinates system on each
such that the first fundamental forms of the two surfaces are identical
at each pair of corresponding points on the two surfaces then the
two surfaces have identical intrinsic geometry. Such surfaces are
isometric and can be mapped on each other by a transformation that
preserves the line lengths, the angles and the surface areas.

\subsection{Functions}

$\bullet$ The domain of a functional mapping: $F:\mathbb{R}^{m}\rightarrow\mathbb{R}^{n}$
is the largest set of $\mathbb{R}^{m}$ on which the mapping is defined.

$\bullet$ A bicontinuous function or mapping is a continuous function
with a continuous inverse.

$\bullet$ A scalar function is of class $C^{n}$ if the function
and all of its first $n$ (but not $n+1$) partial derivatives do
exist and are continuous. A vector function (e.g. a position vector
representing a space curve or a surface) is of class $C^{n}$ if one
of its components is of this class while all the other components
are of this class or higher. A curve or a surface is of class $C^{n}$
if it is mathematically represented by a function of this class.

$\bullet$ In gross terms, a smooth curve or surface means that the
functional relation that represents the object is sufficiently differentiable
for the intended objective, being of class $C^{n}$ at least where
$n$ is the minimum requirement for the differentiability index to
satisfy the required conditions.

$\bullet$ A deleted neighborhood of a point $P$ on a 1D interval
on the real line is defined as the set of all points $x\in\mathbb{R}$
in the interval such that $0<\left|x-x_{P}\right|<\epsilon$ where
$x_{P}$ is the coordinate of $P$ on the real line and $\epsilon$
is a positive real number. Hence, the deleted neighborhood includes
all the points in the open interval ($x_{P}-\epsilon,x_{P}+\epsilon$)
excluding $x_{P}$ itself. For a space curve (which is not straight
in general) represented by $\mathbf{r}=\mathbf{r}(t)$, where $\mathbf{r}$
is the spatial representation of the curve and $t$ is a general parameter
in the curve representation, the definition applies to the neighborhood
of $t_{P}$ where $t_{P}$ is the value of $t$ corresponding to the
point $P$ on the curve.

$\bullet$ A deleted neighborhood of a point $P$ on a 2D flat surface
is defined as the set of all points $(x,y)\in\mathbb{R}^{2}$ on the
surface such that $0<\sqrt{(x-x_{P})^{2}+(y-y_{P})^{2}}<\epsilon$
where $(x_{P},y_{P})$ are the coordinates of $P$ on the plane and
$\epsilon$ is a positive real number. Hence, the deleted neighborhood
includes all the points inside a circle of radius $\epsilon$ and
center $(x_{P},y_{P})$ excluding the center itself. For a space surface
(which is not flat in general) represented by $\mathbf{r}=\mathbf{r}(u,v)$,
where $\mathbf{r}$ is the spatial representation of the surface and
$u,v$ are the surface coordinates, the definition applies to the
neighborhood of ($u_{P},v_{P}$) where ($u_{P},v_{P}$) are the coordinates
on the 2D $uv$ plane corresponding to the point $P$ on the surface.

$\bullet$ A quadratic expression $Q(x,y)=a_{1}x^{2}+2a_{2}xy+a_{3}y^{2}$
of real coefficients $a_{1},a_{2},a_{3}$ and real variables $x,y$
is described as ``positive definite'' if it possesses positive values
($>0$) for all pairs $(x,y)\ne(0,0)$. The sufficient and necessary
condition for $Q$ to be positive definite is that $a_{1}>0$ and
$\left(a_{1}a_{3}-a_{2}a_{2}\right)>0$.\footnote{The conditions $a_{1}>0$ and $\left(a_{1}a_{3}-a_{2}a_{2}\right)>0$
necessitate $a_{3}>0$ since these coefficients are real.}

\subsection{Coordinate Transformations}

$\bullet$ An orthogonal coordinate transformation is a combination
of translation, rotation and reflection of axes. The Jacobian of orthogonal
transformations is unity, that is $J=\pm1$. The orthogonal transformation
is described as positive \textit{iff} $J=+1$ and negative \textit{iff}
$J=-1$. Positive orthogonal transformations consist solely of translation
and rotation (possibly trivial ones as in the case of an identity
transformation) while negative orthogonal transformations include
reflection, by applying an odd number of axes reversal, as well. Positive
transformations can be decomposed into an infinite number of continuously
varying infinitesimal positive transformations each one of which assimilates
an identity transformation. Such a decomposition is not possible in
the case of negative orthogonal transformations because the shift
from the identity transformation to reflection is impossible by a
continuous process.

$\bullet$ Coordinate curves, which are also called parametric curves
or parametric lines, on a surface are curves along which only one
coordinate variable ($u$ or $v$) varies while the other coordinate
variable ($v$ or $u$) remains constant.

$\bullet$ An invariant property of a curve or a surface is a property
which is independent of allowable coordinate transformations and parameterizations.

$\bullet$ A regular representation of class $C^{m}$ ($m>0$) of
a surface patch $S$ in a 3D Euclidean space is a functional mapping
of an open set $\Omega$ in the $uv$ plane onto $S$ that satisfies
the following conditions:

(A) The functional mapping relation is of class $C^{m}$ over the
entire $\Omega$.\footnote{It is noteworthy that the condition of being of class $C^{m}$ in
this context means that $m$ is the minimum requirement for differentiability
and hence the condition is satisfied by any function of class $C^{n}$
where $n\ge m$.}

(B) The Jacobian matrix\footnote{For a functional mapping of the form $\mathbf{S}(u,v)=\left(S_{1}(u,v),S_{2}(u,v),S_{3}(u,v)\right)$,
this Jacobian matrix is given by:
\begin{equation}
\left[\begin{array}{cc}
\partial_{u}S_{1} & \partial_{v}S_{1}\\
\partial_{u}S_{2} & \partial_{v}S_{2}\\
\partial_{u}S_{3} & \partial_{v}S_{3}
\end{array}\right]
\end{equation}
} for the transformation between the representations of the surface
in 3D and 2D spaces is of rank 2 for all the points in $\Omega$.

$\bullet$ Having a Jacobian matrix of rank 2 for the transformation
is equivalent to the condition that $\mathbf{E}_{1}\times\mathbf{E}_{2}\ne\mathbf{0}$
where $\mathbf{E}_{1}=\partial_{u}\mathbf{r}$ and $\mathbf{E}_{2}=\partial_{v}\mathbf{r}$
are the surface basis vectors, which are the tangents to the $u$
and $v$ coordinate curves respectively, and $\mathbf{r}=\mathbf{r}(u,v)$
is the 3D spatial representation of the curves.\footnote{``Rank'' here, and in similar contexts, refers to its meaning in
linear algebra and should not be confused with the rank of tensor.}

$\bullet$ Having a Jacobian matrix of rank 2 is also equivalent to
having a well-defined tangent plane to the surface at the related
point.

$\bullet$ A point on a surface which is not regular is called singular.
Singularity occurs either because of a geometric reason, which is
the case for instance for the apex of a cone, or because of the particular
parametric representation of the surface. While the first type of
singularity is inherent and hence it cannot be removed, the second
type can be removed by changing the representation.

$\bullet$ Corresponding points on two curves refer to two points,
one on each curve, with a common value of a common parameter of the
two curves. When the two curves have two different parameterizations
then a one-to-one correspondence between the two parameters should
be established and the corresponding points then refer to two points
with corresponding values of the two parameters. Corresponding points
on two surfaces are defined in a similar manner taking into account
that surfaces are parameterized by multiple values depending on the
dimensionality of the reference space (e.g. two when using surface
curvilinear coordinates $u^{1},u^{2}$ and three when using 3D space
coordinates $x^{1},x^{2},x^{3}$).

$\bullet$ In many cases of theoretical and practical situations,
a mixed tensor $A_{\alpha}^{i}$, which is contravariant with respect
to transformations in space coordinates $x^{i}$ and covariant with
respect to transformations in surface coordinates $u^{\alpha}$, is
defined. Following a coordinate transformation in which both the space
and surface coordinates change, the tensor $A_{\alpha}^{i}$ will
be given in the new (barred) system by:
\begin{equation}
\bar{A}_{\alpha}^{i}=A_{\beta}^{j}\frac{\partial\bar{x}^{i}}{\partial x^{j}}\frac{\partial u^{\beta}}{\partial\bar{u}^{\alpha}}
\end{equation}
More generally, tensors with space and surface contravariant indices
and space and surface covariant indices (e.g. $A_{j\beta}^{i\alpha}$)
can also be defined similarly. The extension of the above transformation
rule to include such tensors can be conducted trivially by following
the obvious pattern seen in the last equation.

\subsection{Intrinsic Distance}

$\bullet$ The intrinsic distance between two points on a surface
is the greatest lower bound (or infimum) of the lengths of all regular
arcs connecting the two points on the surface.

$\bullet$ The intrinsic distance is an intrinsic property of the
surface.

$\bullet$ The intrinsic distance $d$ between two points is invariant
under a local isometric mapping, that is $d(f(P_{1}),f(P_{2}))=d(P_{1},P_{2})$
where $f$ is an isometric mapping from a surface $S_{1}$ to a surface
$S_{2}$ (see $\S$ \ref{Llocalisometry}), $P_{1}$ and $P_{2}$
are the two points on $S_{1}$ and $f(P_{1})$ and $f(P_{2})$ are
their images on $S_{2}$.\footnote{In fact this may be taken as the definition of isometric mapping,
i.e. it is the mapping that preserves intrinsic distance.}

$\bullet$ The following conditions apply to the intrinsic distance
$d$ between points $P_{1}$, $P_{2}$ and $P_{3}$:

(A) Symmetry: $d(P_{1},P_{2})=d(P_{2},P_{1})$.

(B) Triangle inequality: $d(P_{1},P_{3})\le d(P_{1},P_{2})+d(P_{2},P_{3})$.

(C) Positive definiteness: $d(P_{1},P_{2})>0$ with $d(P_{1},P_{2})=0$
\textit{iff} $P_{1}$ and $P_{2}$ are the same point.

$\bullet$ An arc $C$ connecting two points, $P_{1}$ and $P_{2}$,
on a surface is described as an arc of minimum length\footnote{Such an arc may not be unique (see next point).}
between $P_{1}$ and $P_{2}$ if the length of $C$ is equal to the
intrinsic distance between $P_{1}$ and $P_{2}$.

$\bullet$ The existence and uniqueness of an arc of minimum length
between two specific points on a surface is not guaranteed, i.e. it
may not exit and if it does it may not be unique (refer to $\S$ \ref{subGeodesicCurves}).
Yes, for certain types of surface such an arc does exist and it is
unique. For example, on a plane in a Euclidean space there exists
an arc of minimum length between any two points on the plane and it
is unique; this arc is the straight line segment connecting the two
points.

\subsection{Basis Vectors\label{subBasisVectors}}

$\bullet$ The set of basis vectors in a given manifold plays a pivotal
role in the theoretical construction of the geometry of the manifold,
and this applies to the basis vectors in differential geometry where
these vectors are used in the definition and construction of essential
concepts and objects such as the surface metric tensor. They are also
employed to serve as moving coordinate frames for their underlying
constructions (see $\S$ \ref{subSpaceCurves} and $\S$ \ref{secSurfacesinSpace}).

$\bullet$ The differential geometry of curves and surfaces employs
two main sets of basis vectors:\footnote{Other sets of basis vectors are also defined and employed in differential
geometry, see e.g. $\S$ \ref{subCurvatureVector} and $\S$ \ref{subPrincipalCurvatures}.}

(A) One set is constructed on space curves and consists of the three
unit vectors: the tangent $\mathbf{T}$, the normal $\mathbf{N}$
and the binormal $\mathbf{B}$ to the curve.

(B) Another set is constructed on surfaces and consists of two linearly
independent vectors tangent to the coordinate curves of the surface,
$\mathbf{E}_{1}=\frac{\partial\mathbf{r}}{\partial u^{1}}$ and $\mathbf{E}_{2}=\frac{\partial\mathbf{r}}{\partial u^{2}}$,
plus the normal to the surface $\mathbf{n}$, where $\mathbf{r}(u^{1},u^{2})$
is the spatial representation of a surface coordinate curve, and $u^{1},u^{2}$
are the surface curvilinear coordinates\footnote{The surface curvilinear coordinates are also called the Gaussian coordinates.}
as will be explained in detail later on (refer to $\S$ \textbf{\ref{secSurfacesinSpace}}).

$\bullet$ Each one of the above basis sets is defined on each point
of the curve or the surface and hence in general the vectors in each
one of these basis sets vary from one point to another, i.e. they
are position dependent.\footnote{The non-unit vectors (i.e. $\mathbf{E}_{1}$ and $\mathbf{E}_{2}$)
vary in magnitude and direction while the unit vectors (which are
the rest) vary in direction.}

$\bullet$ In tensor notation, the surface basis vectors, $\mathbf{E}_{1}$
and $\mathbf{E}_{2}$, are given by $\frac{\partial x^{i}}{\partial u^{\alpha}}$
($i=1,2,3$ and $\alpha=1,2$) which is usually abbreviated as $x_{\alpha}^{i}$.
These vectors can be seen as contravariant space vectors or covariant
surface vectors (see $\S$ \ref{subSurfaceMetric} for further details).

\subsection{Flat and Curved Spaces}

$\bullet$ A manifold, such as a 2D surface or a 3D space, is called
``flat'' if it is possible to find a coordinate system for the manifold
with a diagonal metric tensor whose all diagonal elements are $\pm1$;
the space is called ``curved'' otherwise. More formally, an $n$D
space is described as flat \textit{iff} it is possible to find a coordinate
system for which the line element $ds$ is given by:
\begin{equation}
(ds)^{2}=\epsilon_{1}(dx^{1})^{2}+\epsilon_{2}(dx^{2})^{2}+\ldots\epsilon_{n}(dx^{n})^{2}=\sum_{i=1}^{n}\epsilon_{i}(dx^{i})^{2}\label{eqFlatSpaceCondition}
\end{equation}
where the indexed $\epsilon$ are $\pm1$. Examples of flat space
are the 3D Euclidean space which can be coordinated by a rectangular
Cartesian system whose metric tensor is diagonal with all the diagonal
elements being $+1$, and the 4D Minkowski space-time manifold whose
metric is diagonal with elements of $\pm1$. An example of curved
space is the 2D surface of a sphere or an ellipsoid.

$\bullet$ For the space to be flat, the condition given by Eq. \ref{eqFlatSpaceCondition}
should apply all over the space and not just at certain points or
regions.

$\bullet$ As discussed in \cite{SochiTC2}, a necessary and sufficient
condition for an $n$D space to be globally flat is that the Riemann-Christoffel
curvature tensor of the space vanishes identically.

$\bullet$ Due to the connection between the Gaussian curvature and
the Riemann-Christoffel curvature tensor which implies vanishing each
one of these if the other does (see Eq. \ref{eqRK1}), we see that
having an identically vanishing Gaussian curvature (see $\S$ \ref{subGaussianCurvature})
is another sufficient and necessary condition for a 2D space to be
flat.\footnote{The Gaussian curvature in differential geometry is defined for 2D
spaces although the concept may be extended to higher dimensionality
manifolds.}

$\bullet$ Curved spaces may have constant non-vanishing curvature
all over the space, or have variable curvature and hence the curvature
is position dependent. As an example of a space of constant curvature
is the surface of a sphere of radius $R$ whose curvature\footnote{This is the Gaussian curvature (refer to $\S$ \ref{subGaussianCurvature}).}
is $\frac{1}{R^{2}}$ at each point of the surface. Ellipsoid, paraboloid
and torus are simple examples of surfaces with variable curvature.\footnote{There are various characterizations and quantifications for the curvature
and hence in the present context ``curvature'' may be a generic
term. For 2D surfaces, curvature usually refers to the Gaussian curvature
(see $\S$ \ref{subGaussianCurvature}) which is strongly linked to
the Riemann curvature.}

$\bullet$ A surface with positive/negative Gaussian curvature (see
$\S$ \ref{subGaussianCurvature}) at each point is described as a
surface of positive/negative curvature.\footnote{Surfaces with constant non-zero Gaussian curvature $K$ may be described
as spherical if $K>0$ and pseudo-spherical if $K<0$.} Ellipsoids, elliptic paraboloids and hyperboloids of two sheets are
examples of surfaces of positive curvature while hyperbolic paraboloids
and hyperboloids of one sheet are examples of surfaces of negative
curvature (see $\S$ \ref{subQuadraticSurfaces}).

$\bullet$ Schur theorem related to $n$D spaces ($n>2$) of constant
curvature states that: if the Riemann-Christoffel curvature tensor
at each point of a space is a function of the coordinates only, then
the curvature is constant all over the space.\footnote{Schur theorem may also be stated as: the Riemannian curvature is constant
over an isotropic region of an $n$D ($n>2$) Riemannian space. }

$\bullet$ All 1D spaces are Euclidean and hence they cannot be curved.
So twisted curves are curved only when viewed externally from the
embedding space which they reside in (e.g. the 2D space of a surface
curve or the 3D space of a space curve).

$\bullet$ The geometry of curved spaces is usually described as the
Riemannian geometry. One approach for investigating the Riemannian
geometry of a curved manifold is to embed the manifold in a Euclidean
space of higher dimensionality. This approach is largely followed
in the present notes where the geometry of curved 2D spaces (twisted
surfaces) is investigated by immersing the surfaces in a 3D Euclidean
space and examining their properties as viewed from this external
enveloping 3D space. Such an external view is necessary for examining
the extrinsic geometry of the surface but not its intrinsic geometry.

$\bullet$ The geometric description and quantification of flat spaces
are simpler than those of curved spaces, and hence in general the
differential geometry of flat spaces is wealthier, more motivating
and less challenging than that of curved spaces.\footnote{As there is a subjective element in this type of statements, it may
not apply to everyone.}

\subsection{Homogeneous Coordinate Systems}

$\bullet$ When all the diagonal elements of a diagonal metric tensor
of a flat space are $+1$, the coordinate system is described as homogeneous.
In this case the line element of Eq. \ref{eqFlatSpaceCondition} becomes:
\begin{equation}
(ds)^{2}=dx^{i}dx^{i}\label{eqLineElemHomoSys}
\end{equation}
An example of homogeneous coordinate systems is the rectangular Cartesian
system ($x,y,z$) of a 3D Euclidean space.

$\bullet$ A homogeneous coordinate system can be transformed to another
homogeneous coordinate system only by linear transformations.

$\bullet$ Any coordinate system obtained from a homogeneous coordinate
system by an orthogonal transformation is homogeneous.

$\bullet$ As a consequence of the last points, infinitely many homogeneous
coordinate systems can be constructed in any flat space.

$\bullet$ A coordinate system of a flat space can always be homogenized
by allowing the coordinates to be imaginary. This is done by redefining
the coordinates as:
\begin{equation}
\underline{x}^{i}=\sqrt{\epsilon_{i}}x^{i}
\end{equation}
where the new coordinates $\underline{x}^{i}$ are imaginary when
$\epsilon_{i}=-1$. Consequently, the line element will be given by:
\begin{equation}
(ds)^{2}=d\underline{x}^{i}d\underline{x}^{i}
\end{equation}
which is of the same form as Eq. \ref{eqLineElemHomoSys}. An example
of a homogeneous coordinate system with some real and some imaginary
coordinates is the coordinate system of a Minkowski 4D space-time
of special relativity.

\subsection{Geodesic Coordinates\label{subGeodesicCoordinates}}

$\bullet$ It is always possible to introduce coordinates at particular
points in a multi-dimensional manifold so that the Christoffel symbols
vanish at these points. These coordinates are called geodesic coordinates.\footnote{Some authors define geodesic coordinates on a coordinate patch of
a surface as a coordinate system whose $u$ and $v$ coordinate curve
families are orthogonal with one of these families ($u$ or $v$)
being a family of geodesic curves (refer to $\S$ \ref{subGeodesicCurves}).
So, ``geodesic coordinates'' seems to have multiple usage.}

$\bullet$ Geodesic coordinates are employed as local coordinate systems
for certain advantages. In geodesic coordinates the Christoffel symbols
are made to vanish at certain allocated points described as the poles.
These systems are called geodesic for these particular points and
also described as locally Cartesian coordinates.

$\bullet$ The main reason for the use of geodesic coordinates is
that the covariant and absolute derivatives in such systems become
respectively partial and total derivatives at the poles since the
Christoffel symbol terms in the covariant and absolute derivative
expressions vanish at these points. Any tensor property can then be
easily proved in the geodesic system at the pole and consequently
generalized to other systems due to the invariance of the zero tensor
under permissible coordinate transformations. If the allocated pole
is a general point in the space, the property is then established
over the whole space.

$\bullet$ In any Riemannian space it is possible to find a coordinate
system for which the coordinates are geodesic at every point of a
given analytic curve.

$\bullet$ There is an infinite number of ways by which geodesic coordinates
can be defined over a coordinate patch.

\subsection{Christoffel Symbols for Curves and Surfaces}

$\bullet$ For 1D spaces, the Christoffel symbols are not defined.

$\bullet$ The Christoffel symbols of the first kind for a 2D surface
are given by:
\begin{equation}
\begin{aligned}\left[11,1\right] & =\frac{\partial_{u}a_{11}}{2}=\frac{E_{u}}{2}\\
\left[11,2\right] & =\partial_{u}a_{12}-\frac{\partial_{v}a_{11}}{2}=F_{u}-\frac{E_{v}}{2}\\
\left[12,1\right] & =\frac{\partial_{v}a_{11}}{2}=\frac{E_{v}}{2}=\left[21,1\right]\\
\left[12,2\right] & =\frac{\partial_{u}a_{22}}{2}=\frac{G_{u}}{2}=\left[21,2\right]\\
\left[22,1\right] & =\partial_{v}a_{12}-\frac{\partial_{u}a_{22}}{2}=F_{v}-\frac{G_{u}}{2}\\
\left[22,2\right] & =\frac{\partial_{v}a_{22}}{2}=\frac{G_{v}}{2}
\end{aligned}
\end{equation}
where the indexed $a$ are the elements of the surface covariant metric
tensor (refer to \ref{subSurfaceMetric}) and $E,F,G$ are the coefficients
of the first fundamental form (refer to \ref{subFirstFundamentalForm}).\footnote{As we will see in the forthcoming sections, $E=a_{11}$, $F=a_{12}=a_{21}$
and $G=a_{22}$.} The subscripts $u$ and $v$ which suffix the coefficients stand
for partial derivatives with respect to these variables (i.e. $\frac{\partial}{\partial u}$
and $\frac{\partial}{\partial v}$). In orthogonal coordinate systems
$F=a_{12}=a_{21}=0$ and hence these formulae will be simplified accordingly
by dropping any term involving these coefficients.

$\bullet$ The Christoffel symbols of the first kind are linked to
the surface covariant basis vectors by the following relation:

\begin{equation}
\left[\alpha\beta,\gamma\right]=\frac{\partial\mathbf{E}_{\alpha}}{\partial u^{\beta}}\cdot\mathbf{E}_{\gamma}\,\,\,\,\,\,\,\,\,\,\,\,\text{(\ensuremath{\alpha,\beta,\gamma=1,2})}\label{eqChris1}
\end{equation}
which may be written as:
\begin{equation}
\left[\alpha\beta,\gamma\right]=\mathbf{r}_{\alpha\beta}\cdot\mathbf{r}_{\gamma}\,\,\,\,\,\,\,\,\,\,\,\,\text{(\ensuremath{\alpha,\beta,\gamma=1,2})}
\end{equation}
where the subscripts represent partial derivatives with respect to
these coordinate indices. The last equation provides an easier form
to remember these formulae.

$\bullet$ Applying the index raising operator to Eq. \ref{eqChris1},
we obtain a similar expression for the Christoffel symbols of the
second kind, that is:
\begin{equation}
\Gamma_{\alpha\beta}^{\gamma}=\frac{\partial\mathbf{E}_{\alpha}}{\partial u^{\beta}}\cdot\mathbf{E}^{\gamma}\,\,\,\,\,\,\,\,\,\,\,\,\text{(\ensuremath{\alpha,\beta,\gamma=1,2})}\label{eqChris2}
\end{equation}
where $\mathbf{E}^{\gamma}$ is the contravariant form of the surface
basis vectors.

$\bullet$ The Christoffel symbols of the second kind\footnote{They are also called affine connections.}
for a 2D surface are given by:
\begin{equation}
\begin{aligned}\Gamma_{11}^{1} & =\frac{a_{22}\partial_{u}a_{11}-2a_{12}\partial_{u}a_{12}+a_{12}\partial_{v}a_{11}}{2a}=\frac{GE_{u}-2FF_{u}+FE_{v}}{2a}\\
\Gamma_{11}^{2} & =\frac{2a_{11}\partial_{u}a_{12}-a_{11}\partial_{v}a_{11}-a_{12}\partial_{u}a_{11}}{2a}=\frac{2EF_{u}-EE_{v}-FE_{u}}{2a}\\
\Gamma_{12}^{1} & =\frac{a_{22}\partial_{v}a_{11}-a_{12}\partial_{u}a_{22}}{2a}=\frac{GE_{v}-FG_{u}}{2a}=\Gamma_{21}^{1}\\
\Gamma_{12}^{2} & =\frac{a_{11}\partial_{u}a_{22}-a_{12}\partial_{v}a_{11}}{2a}=\frac{EG_{u}-FE_{v}}{2a}=\Gamma_{21}^{2}\\
\Gamma_{22}^{1} & =\frac{2a_{22}\partial_{v}a_{12}-a_{22}\partial_{u}a_{22}-a_{12}\partial_{v}a_{22}}{2a}=\frac{2GF_{v}-GG_{u}-FG_{v}}{2a}\\
\Gamma_{22}^{2} & =\frac{a_{11}\partial_{v}a_{22}-2a_{12}\partial_{v}a_{12}+a_{12}\partial_{u}a_{22}}{2a}=\frac{EG_{v}-2FF_{v}+FG_{u}}{2a}
\end{aligned}
\label{eqChristoffel2}
\end{equation}
where $a$ is the determinant of the surface covariant metric tensor,
and the other symbols are as explained in the previous points. The
formulae will also be simplified in orthogonal coordinate systems
where $F=a_{12}=a_{21}=0$ by dropping the vanishing terms.

$\bullet$ The Christoffel symbols of the second kind for a 2D surface
may also be given by:
\begin{equation}
\begin{aligned}\Gamma_{11}^{1} & =-\frac{\left(\mathbf{E}_{2}\times\partial_{1}\mathbf{E}_{1}\right)\cdot\mathbf{n}}{\sqrt{a}}\\
\Gamma_{11}^{2} & =+\frac{\left(\mathbf{E}_{1}\times\partial_{1}\mathbf{E}_{1}\right)\cdot\mathbf{n}}{\sqrt{a}}\\
\Gamma_{12}^{1} & =-\frac{\left(\mathbf{E}_{2}\times\partial_{2}\mathbf{E}_{1}\right)\cdot\mathbf{n}}{\sqrt{a}}=\Gamma_{21}^{1}\\
\Gamma_{12}^{2} & =+\frac{\left(\mathbf{E}_{1}\times\partial_{2}\mathbf{E}_{1}\right)\cdot\mathbf{n}}{\sqrt{a}}=\Gamma_{21}^{2}\\
\Gamma_{22}^{1} & =-\frac{\left(\mathbf{E}_{2}\times\partial_{2}\mathbf{E}_{2}\right)\cdot\mathbf{n}}{\sqrt{a}}\\
\Gamma_{22}^{2} & =+\frac{\left(\mathbf{E}_{1}\times\partial_{2}\mathbf{E}_{2}\right)\cdot\mathbf{n}}{\sqrt{a}}
\end{aligned}
\label{eqChristoffel2b}
\end{equation}
where the indexed $\mathbf{E}$ are the surface covariant basis vectors,
$\mathbf{n}$ is the unit vector normal to the surface and $a$ is
the determinant of the surface covariant metric tensor.

$\bullet$ Since the Christoffel symbols of both kinds are dependent
on the metric only, as can be seen from the previous points, they
represent intrinsic properties of the surface geometry and hence they
are part of its intrinsic geometry.

\subsection{Riemann-Christoffel and Ricci Curvature Tensors\label{subRiemannTensor}}

$\bullet$ The Riemann-Christoffel curvature tensor is an absolute
rank-4 tensor that characterizes important properties of spaces, including
2D surfaces, and hence it plays an important role in differential
geometry. The tensor is used, for instance, to test for the space
flatness.

$\bullet$ There are two kinds of Riemann-Christoffel curvature tensor:
first and second. The Riemann-Christoffel curvature tensor of the
first kind is a type ($0,4$) tensor while the Riemann-Christoffel
curvature tensor of the second kind is a type ($1,3$) tensor. Shifting
from one kind to the other is achieved by using the index-shifting
operator.

$\bullet$ The first and second kinds of the Riemann-Christoffel curvature
tensor are given respectively by:\footnote{In these and the following equations in this subsection, the Latin
indices do not necessarily range over $1,2,3$ as these equations
are valid for general $n$D manifolds ($n\ge2$) including surfaces
and spaces with $n>3$.}
\begin{equation}
\begin{aligned}R_{ijkl} & =\partial_{k}\left[jl,i\right]-\partial_{l}\left[jk,i\right]+\left[il,r\right]\Gamma_{jk}^{r}-\left[ik,r\right]\Gamma_{jl}^{r}\\
R_{\,jkl}^{i} & =\partial_{k}\Gamma_{jl}^{i}-\partial_{l}\Gamma_{jk}^{i}+\Gamma_{jl}^{r}\Gamma_{rk}^{i}-\Gamma_{jk}^{r}\Gamma_{rl}^{i}
\end{aligned}
\label{eqRieChrTensor}
\end{equation}

$\bullet$ The Riemann-Christoffel curvature tensor vanishes identically
\textit{iff} the space is globally flat; otherwise the space is curved.

$\bullet$ A surface is isometric to the Euclidean plane \textit{iff}
the Riemann-Christoffel curvature tensor is zero at each point on
the surface.\footnote{The same statement also applies if the Gaussian curvature of the surface
vanishes identically due to the link between the Riemann and Gaussian
curvatures (refer to $\S$ \ref{subGaussianCurvature} and see Eq.
\ref{eqRK1})}

$\bullet$ From Eq. \ref{eqRieChrTensor}, it can be seen that the
Riemann-Christoffel curvature tensor depends exclusively on the Christoffel
symbols of the first and second kind which are both dependent on the
metric (or the first fundamental form, see $\S$ \ref{subFirstFundamentalForm})
only. Hence, the Riemann-Christoffel curvature, as represented by
the Riemann-Christoffel curvature tensor, is an intrinsic property
of the manifold.

$\bullet$ Since the Riemann-Christoffel curvature tensor depends
on the metric which, in general curvilinear coordinates, is a function
of position, the Riemann-Christoffel curvature tensor follows this
dependency on position.

$\bullet$ The Ricci curvature tensor of the first kind, which is
a rank-2 symmetric tensor, is obtained by contracting the contravariant
index with the last covariant index of the Riemann-Christoffel curvature
tensor of the second kind, that is:
\begin{equation}
R_{ij}=R_{ija}^{a}=\partial_{j}\Gamma_{ia}^{a}-\partial_{a}\Gamma_{ij}^{a}+\Gamma_{bj}^{a}\Gamma_{ia}^{b}-\Gamma_{ba}^{a}\Gamma_{ij}^{b}
\end{equation}
The Ricci tensor of the second kind is obtained by raising the first
index of the Ricci tensor of the first kind using the index-raising
operator.

$\bullet$ The Ricci scalar, which is also called the curvature scalar
and the curvature invariant, is the result of contracting the indices
of the Ricci curvature tensor of the second kind, that is:
\begin{equation}
R=R_{\,i}^{i}
\end{equation}

$\bullet$ The Riemann-Christoffel curvature tensor vanishes identically
for 1D manifolds as represented by space and surface curves.

$\bullet$ As discussed in \cite{SochiTC2}, the 2D Riemann-Christoffel
curvature tensor has only one degree of freedom and hence it possesses
a single independent non-vanishing component which is represented
by $R_{1212}$. Hence for a 2D Riemannian space we have:
\begin{equation}
R_{1212}=R_{2121}=-R_{1221}=-R_{2112}
\end{equation}
while all the other components of the tensor are identically zero.
This can be expressed in a single equation as:
\begin{equation}
R_{\alpha\beta\gamma\delta}=R_{1212}\epsilon_{\alpha\beta}\epsilon_{\gamma\delta}=\frac{R_{1212}}{a}\underline{\epsilon}_{\alpha\beta}\underline{\epsilon}_{\gamma\delta}=K\underline{\epsilon}_{\alpha\beta}\underline{\epsilon}_{\gamma\delta}\label{eqRK1}
\end{equation}
where $K$ is the Gaussian curvature (see $\S$ \ref{subGaussianCurvature}).

$\bullet$ The non-vanishing component of the 2D Riemann-Christoffel
curvature tensor, $R_{1212}$, is given by:\footnote{Here, $\partial_{\alpha\beta}\equiv\frac{\partial^{2}}{\partial u^{\alpha}\partial u^{\beta}}$.}
\begin{equation}
R_{1212}=\frac{1}{2}\left(2\partial_{12}a_{12}-\partial_{22}a_{11}-\partial_{11}a_{22}\right)+a_{\alpha\beta}\left(\Gamma_{12}^{\alpha}\Gamma_{12}^{\beta}-\Gamma_{11}^{\alpha}\Gamma_{22}^{\beta}\right)\label{eqR1212}
\end{equation}
where the indexed $a$ are the coefficients of the surface covariant
metric tensor, the Christoffel symbols are based on the surface metric,
and $\alpha,\beta=1,2$.

$\bullet$ For 2D spaces, the Riemann-Christoffel curvature tensor
is related to the Ricci tensor by the following relations:
\begin{equation}
\frac{R_{1212}}{a}=-\frac{R_{11}}{a_{11}}=-\frac{R_{12}}{a_{12}}=-\frac{R_{21}}{a_{21}}=-\frac{R_{22}}{a_{22}}
\end{equation}
where $a$ is the determinant of the 2D covariant metric tensor (see
$\S$ \ref{subSurfaceMetric}) and the indexed $a$ are its elements.
Since $K=\frac{R_{1212}}{a}$, the above relations also link the Gaussian
curvature to the Ricci tensor.

$\bullet$ More details about the Riemann-Christoffel curvature tensor,
Ricci curvature tensor and Ricci scalar can be found in \cite{SochiTC2}.

\subsection{Curves\label{subCurves}}

$\bullet$ In simple terms, a space curve is a set of connected points\footnote{The points are usually assumed to be totally connected so that any
point on the curve can be reached from any other point by passing
through other curve points or at least they are piecewise connected.
We also consider mostly open curves with simple connectivity and hence
the curve does not intersect itself.} in the space such that any totally-connected subset of it can be
twisted into a straight line segment without affecting the neighborhood
of any point. More technically, a curve is defined as a differentiable
parameterized mapping between an interval of the real line and a connected
subset of the space, that is $C(t):I\rightarrow\mathbb{R}^{3}$ where
$C$ is a space curve defined on the interval $I$$\subseteq\mathbb{R}$
and parameterized by the variable $t\in I$. Hence different parameterizations
of the same ``geometric curve'' will lead to different ``mapping
curves''. The image of the mapping in $\mathbb{R}^{3}$ is known
as the trace of the curve; hence different mapping curves can have
identical traces. The curve may also be defined as a topological image
of a real interval and may be linked to the concept of Jordan arc.\footnote{Jordan arcs or Jordan curves are injective mappings with no self intersection.}

$\bullet$ Space curves can be defined symbolically in different ways;
the most common of these is parametrically where the three space coordinates
are given as functions of a real valued parameter, e.g. $x^{i}=x^{i}(t)$
where $t\in\mathbb{R}$ is the curve parameter and $i=1,2,3$. The
parameter $t$ may represent time or arc length or even an arbitrarily
defined quantity. Similarly, surface curves are defined parametrically
where the two surface coordinates are given as functions of a real
valued parameter, e.g. $u^{\alpha}=u^{\alpha}(t)$ with $\alpha=1,2$.\footnote{In the present notes, many statements like this one about space curves
apply as well to surface curves and vice versa. To be concise and
avoid repetition and unnecessary complication of text we usually talk
about one type only since the extension to the other type should be
obvious with consideration of the difference in dimensionality and
symbols.}

$\bullet$ Parameterized curves are oriented objects as they can be
traversed in one direction or the other depending on the sense of
increase of their parameter.

$\bullet$ The condition for a space curve $C(t):I\rightarrow\mathbb{R}^{3}$,
where $t\in I$ is the curve parameter and $I\subseteq\mathbb{R}$
is an interval over which the curve is defined, to be parameterized
by arc length is that: for all $t$ we have $|\frac{d\mathbf{r}}{dt}|=1$
where $\mathbf{r}(t)$ is the position vector representing the curve
in the 3D ambient space.

$\bullet$ As a consequence of the last point, parameterization by
arc length is equivalent to traversing the curve with unity speed.

$\bullet$ The parameter symbol which is used normally for parameterization
by arc length is $s$, while $t$ is used to represent a general parameter
which could be arc length or something else. This notation is followed
in the present text.

$\bullet$ For curves parameterized by arc length, the length of a
segment between two points on the curve corresponding to $s_{1}$
and $s_{2}$ is given by the simple formula: $L=\left|\int_{s_{1}}^{s_{2}}dt\right|=\left|s_{2}-s_{1}\right|$.

$\bullet$ Parameterization by arc length $s$ may be called natural
representation or natural parameterization of the curve and hence
$s$ is called natural parameter.

$\bullet$ Natural parameterization is not unique; however any other
natural parameterization $\check{s}$ is related to a given natural
parameterization $s$ by the relation $\check{s}=\pm s+c$ where $c$
is a real constant and hence the above-stated condition $|\frac{d\mathbf{r}}{dt}|=1$
remains valid.\footnote{In this formula, $t$ is a generic symbol and hence it stands for
$s$.} This may be stated in a different way by saying that natural parameterization
with arc length $s$ is unique apart from the possibility of having
a different sense of orientation and an additive constant to $s$.

$\bullet$ Natural parameterization may also be used for parameterization
which is proportional to $s$ and hence the transformation relation
between two natural parameterizations becomes $\check{s}=\pm ms+c$
where $m$ is another real constant. The two parameterizations then
differ, apart from the sense of orientation and the constant shift,
by the length scale which can be chosen arbitrarily.

$\bullet$ In a general $n$D space, the tangent vector to a space
curve, represented parametrically by the spatial representation $\mathbf{r}(t)$
where $t$ is a general parameter, is given by $\frac{d\mathbf{r}}{dt}$.

$\bullet$ The tangent line to a sufficiently smooth curve at one
of its non-singular points $P$ is a straight line passing through
$P$ but not through any point in a deleted neighborhood of $P$.
More technically, a tangent line to a curve $C$ at a point $P$ is
a straight line passing through $P$ and having the same orientation
as the derivative $\frac{d\mathbf{r}}{dt}$ where $\mathbf{r}(t)$
is the spatial representation of $C$. A vector tangent to a space
curve at $P$ is a vector oriented in either directions of the tangent
line at $P$ and hence it is a non-trivial scalar multiple of $\frac{d\mathbf{r}}{dt}$.

$\bullet$ The tangent line to a curve at a given point $P$ on the
curve may also be defined as the limit of a secant line passing through
$P$ and another neighboring point on the curve where the other point
converges, while staying on the curve, to the tangent point. All these
different definitions are equivalent as they represent the same entity.

$\bullet$ A vector tangent to a surface curve, represented parametrically
by: $C(u(t),v(t))$ where $u$ and $v$ are the surface curvilinear
coordinates and $t$ is a general parameter, is given by:

\begin{equation}
\frac{d\mathbf{r}}{dt}=\frac{\partial\mathbf{r}}{\partial u}\frac{du}{dt}+\frac{\partial\mathbf{r}}{\partial v}\frac{dv}{dt}
\end{equation}
where $\mathbf{r}(u(t),v(t))$ is the spatial representation of $C$
and all these quantities are defined and evaluated at a particular
point on the curve. The last equation, in tensor notation, becomes:
\begin{equation}
\frac{dx^{i}}{dt}=\frac{\partial x^{i}}{\partial u^{\alpha}}\frac{du^{\alpha}}{dt}=x_{\alpha}^{i}\frac{du^{\alpha}}{dt}
\end{equation}
where $i=1,2,3$, $\alpha=1,2$ and $(u^{1},u^{2})\equiv(u,v)$.

$\bullet$ A space curve $C(t):I\rightarrow\mathbb{R}^{3}$, where
$I\subseteq R$ and $t\in I$ is a parameter, is ``regular at point
$t_{0}$'' \textit{iff} $\dot{C}(t_{0})$ exists and $\dot{C}(t_{0})\ne0$
where the overdot stands for differentiation with respect to the general
parameter $t$. The curve is ``regular'' \textit{iff} it is regular
at each interior point in $I$.

$\bullet$ On a regular parameterized curve there is a neighborhood
to each point in its domain in which the curve is injective.

$\bullet$ On transforming a surface $S$ by a differentiable regular
mapping $f$ of class $C^{n}$ to a surface $\bar{S}$, a regular
curve $C$ of class $C^{n}$ on $S$ will be mapped on a regular curve
$\bar{C}$ of class $C^{n}$ on $\bar{S}$ by the same functional
mapping relation, that is $\mathbf{\bar{r}}(t)=f(\mathbf{r}(t))$
where the barred and unbarred $\mathbf{r}(t)$ are the spatial parametric
representations of the two curves on the barred and unbarred surfaces.

$\bullet$ A non-trivial vector $\mathbf{v}$ is said to be tangent
to a regular surface $S$ at a given point $P$ on $S$ if there is
a regular curve $C$ on $S$ passing through $P$ such that $\mathbf{v}=\frac{d\mathbf{r}(t)}{dt}$
where $\mathbf{r}(t)$ is the spatial representation of $C$ and $\frac{d\mathbf{r}(t)}{dt}$
is evaluated at $P$ (also see \ref{subSurfaces} for further details).\footnote{In fact any vector $\mathbf{v}=c\frac{d\mathbf{r}(t)}{dt}$, where
$c\ne0$ is a real number, is a tangent although it may not be \textit{the}
tangent.}

$\bullet$ A periodic curve $C$ is a curve that can be represented
parametrically by a continuous function of the form $\mathbf{r}(t+T)=\mathbf{r}(t)$
where $\mathbf{r}$ is the spatial representation of $C$, $t$ is
a real parameter and $T$ is a real constant called the function period.
Circles and ellipses are prominent examples of periodic curves where
they can be represented respectively by $\mathbf{r}(t)=\left(a\cos t,a\sin t\right)$
and $\mathbf{r}(t)=\left(a\cos t,b\sin t\right)$ where $a$ and $b$
are real constants, $t\in\mathbb{R}$ and $\mathbf{r}(t+2\pi)=\mathbf{r}(t)$.
Hence, circles and ellipses are periodic curves with a period of $2\pi$.

$\bullet$ A closed curve is a periodic curve defined over a minimum
of one period.\footnote{Periodicity is not a necessary requirement for the definition of closed
curves as the curves can be defined over a single period without being
considered as such.}

$\bullet$ Closed curves may be regarded as topological images of
circles.

$\bullet$ A curve is described as plane curve if it can be embedded
entirely in a plane with no distortion.

$\bullet$ Orthogonal trajectories of a given family of curves is
a family of curves that intersect the given family perpendicularly
at their intersection points.

$\bullet$ Any curve can be mapped isometrically to a straight line
segment where both are naturally parameterized by arc length.\footnote{From this statement plus the fact that isometric transformation is
an equivalence relation (see $\S$ \ref{subIsometricSurfaces}), it
can be concluded that any two curves of equal length can be mapped
isometrically to each other.}

\subsection{Surfaces\label{subSurfaces}}

$\bullet$ A 2D surface embedded in a 3D space may be defined loosely
as a set of connected points in the space such that the immediate
neighborhood of each point on the surface can be deformed continuously
to form a flat disk. Technically, a surface in a 3D manifold is a
mapping from a subset of coordinate plane to a 3D space, that is $S:\Omega\rightarrow\mathbb{R}^{3}$,
where $\Omega$ is a subset of $\mathbb{R}^{2}$ plane and $S$ is
a sufficiently smooth injective function. Other conditions may also
be imposed to ensure the existence of a tangent plane and a normal
at each point of the surface. In particular, the condition $\partial_{u}\mathbf{r}\times\partial_{v}\mathbf{r}\ne\mathbf{0}$
at all points on the surface is usually imposed to ensure regularity.
Like space and surface curves, the image of the mapping in $\mathbb{R}^{3}$
is known as the trace of the surface.\footnote{For convenience, in the present text we use curve and surface for
trace as well as for mapping; the meaning should be obvious from the
context. Also, the trace of a curve or a surface should not be confused
with the trace of a matrix which is the sum of its diagonal elements.}

$\bullet$ A 2D surface embedded in a 3D space can be defined explicitly:
$z=f(x,y)$, or implicitly: $F(x,y,z)=0$, or parametrically: $x(u^{1},u^{2}),\,y(u^{1},u^{2}),\,z(u^{1},u^{2})$
where $u^{1}$ and $u^{2}$ are independent parameters described as
the curvilinear coordinates of the surface. By substitution, elimination
and algebraic manipulation these forms can be transformed interchangeably.

$\bullet$ A coordinate patch of a surface is an injective, bicontinuous,
regular, parametric representation of a part of the surface. In more
technical terms, a coordinate patch of class $C^{n}$ ($n>0$) on
a surface $S$ is a functional mapping of an open set $\Omega$ in
the $uv$ plane onto $S$ that satisfies the following conditions:

(A) The functional mapping relation is of class $C^{n}$ over $\Omega$.

(B) The mapping is one-to-one and bicontinuous over $\Omega$.

(C) $\mathbf{E}_{1}\times\mathbf{E}_{2}\ne\mathbf{0}$ at any point
in $\Omega$.

$\bullet$ As indicated previously, a vector $\frac{d\mathbf{r}(t)}{dt}$,
where $\mathbf{r}$ is a $t$-parameterized position vector with $t\in I\subseteq\mathbb{R}$,
is described as a tangent vector to the surface $S$ at point $P$
on the surface if there is a regular curve embedded in $S$ and passing
through $P$ such that $\frac{d\mathbf{r}(t)}{dt}$ (or a non-trivial
scalar multiple of $\frac{d\mathbf{r}(t)}{dt}$) is a tangent to the
curve at $P$. The set of all tangent vectors to the surface $S$
at point $P$ forms a tangent plane to $S$ at $P$. This set is called
the tangent space of $S$ at $P$ and it is notated with $T_{P}S$.

$\bullet$ As we will see (also refer to $\S$ \ref{subBasisVectors}),
the tangent space of a regular surface at a given point on the surface
is the span of the two linearly independent basis vectors defined
as $\mathbf{E}_{1}=\frac{\partial\mathbf{r}}{\partial u^{1}}$ and
$\mathbf{E}_{2}=\frac{\partial\mathbf{r}}{\partial u^{2}}$ where
$\mathbf{r}(u^{1},u^{2})$ is the spatial representation of the coordinate
curves in a 3D coordinate system and $u^{1}$ and $u^{2}$ are the
curvilinear coordinates of the surface. The tangent space therefore
is the plane passing through $P$ and is perpendicular to the vector
$\mathbf{E}_{1}\times\mathbf{E}_{2}$.

$\bullet$ As indicated previously, every vector tangent to a regular
surface $S$ at a given point $P$ on $S$ can be expressed as a linear
combination of the surface basis vectors $\mathbf{E}_{1}$ and $\mathbf{E}_{2}$
at $P$. The reverse is also true, that is every linear combination
of $\mathbf{E}_{1}$ and $\mathbf{E}_{2}$ at $P$ is a tangent vector
to a regular curve embedded in $S$ and passing through $P$ and hence
is a tangent to $S$ at $P$.

$\bullet$ The tangent space at a specific point $P$ of a surface
is a property of the surface at $P$ and hence it is independent of
the patch that contains $P$.

$\bullet$ For any non-trivial vector $\mathbf{v}$ which is parallel
to the tangent plane of a simple and smooth surface $S$ at a given
point $P$ on $S$, there is a curve in $S$ passing through $P$
and represented parametrically by $\mathbf{r}(t)$ such that $\mathbf{v}=c\frac{d\mathbf{r}}{dt}$
where $c\ne0$ is a real constant.

$\bullet$ A non-trivial vector is parallel to the tangent plane of
a surface $S$ at a given point $P$ \textit{iff} it is tangent to
$S$ at $P$.

$\bullet$ From the previous points, we see that the tangent plane
of a surface at a given point $P$ is given by:
\begin{equation}
\mathbf{r}=\mathbf{r}_{P}+p\mathbf{E}_{1}+q\mathbf{E}_{2}
\end{equation}
where $\mathbf{r}$ is the position vector of an arbitrary point on
the tangent plane, $\mathbf{r}_{P}$ is the position vector of the
point $P$, $p,q\in(-\infty,\infty)$ are real variables, and $\mathbf{E}_{1}$
and $\mathbf{E}_{1}$ are the surface basis vectors at $P$.

$\bullet$ The straight line passing through a given point $P$ on
a surface $S$ in the direction of the normal vector $\mathbf{n}$
$\left(=\frac{\mathbf{E}_{1}\times\mathbf{E}_{2}}{\left|\mathbf{E}_{1}\times\mathbf{E}_{2}\right|}\right)$
of $S$ at $P$ is called the normal line to $S$ at $P$.\footnote{This should not be confused with the normal line of a space curve
$C$ at point $P$ which is normally called the principal normal line
(see $\S$ \ref{subSpaceCurves}). Anyway, the two should be distinguished
easily by noticing their affiliation to a surface or a curve.} The equation of this normal line is given by:
\begin{equation}
\mathbf{r}=\mathbf{r}_{P}+\lambda\mathbf{n}
\end{equation}
where $\mathbf{r}$ is the position vector of an arbitrary point on
the normal line, $\mathbf{r}_{P}$ is the position vector of the point
$P$, $\lambda\in(-\infty,\infty)$ is a real variable, and $\mathbf{n}$
is the unit normal vector of $S$ at $P$.

$\bullet$ A regular curve of class $C^{n}$ on a sufficiently smooth
surface is an image of a unique regular plane curve of class $C^{n}$
in the parameter plane $\Omega$.\footnote{Here, we are considering each connected part of the curve embedded
in a coordinate patch if there is no single patch that contains the
entire curve.}

$\bullet$ A surface is regular at a given point $P$ \textit{iff}
$\mathbf{E}_{1}\times\mathbf{E}_{2}\ne\mathbf{0}$ at $P$ where $\mathbf{E}_{1}=\partial_{u}\mathbf{r}$
and $\mathbf{E}_{2}=\partial_{v}\mathbf{r}$ are the tangent vectors
to the surface coordinate curves. A surface is regular \textit{iff}
$\mathbf{E}_{1}\times\mathbf{E}_{2}\ne\mathbf{0}$ at any point on
the surface.

$\bullet$ A surface of revolution is an axially-symmetric surface
generated by a plane curve $C$ revolving around a straight line $L$
contained in the plane of the curve but not intersecting the curve.
The curve $C$ is called the profile of the surface and the line $L$
is called the axis of revolution which is also the axis of symmetry
of the surface.

$\bullet$ Meridians of a surface of revolution are plane curves on
the surface formed by the intersection of a plane containing the axis
of revolution with the surface. Parallels of a surface of revolution
are circles generated by intersecting the surface by planes perpendicular
to the axis of revolution. Meridians and parallels intersect at right
angles.\footnote{For spheres, these curves are called meridians of longitude and parallels
of latitude.}

$\bullet$ A ``Monge patch'' is a coordinate patch in a 3D space
defined by a function in one of the following forms:
\begin{eqnarray}
\mathbf{r}(u,v) & = & \left(f(u,v),u,v\right)\nonumber \\
\mathbf{r}(u,v) & = & \left(u,f(u,v),v\right)\\
\mathbf{r}(u,v) & = & \left(u,v,f(u,v)\right)\nonumber
\end{eqnarray}
where $f$ is a differentiable function of the surface coordinates
($u,v$). When $f$ is of class $C^{n}$ then the coordinate patch
is of this class.

$\bullet$ A ``simply connected'' region on a surface means that
a closed curve contained in the region can be shrunk down continuously
onto any point in the region without leaving the region. In simple
terms, it means that the region contains no holes or gaps that separate
its parts.

$\bullet$ In simple terms, a simple surface is a continuously deformed
plane by compression, stretching and bending. Examples of simple surfaces
are cylinders, cones and paraboloids.

$\bullet$ A connected surface $S$ is a simple surface which cannot
be entirely represented by the union of two disjoint open point sets
in $\mathbb{R}^{3}$ where these sets have non-empty intersection
with $S$. Hence, for any two arbitrary points, $P_{1}$ and $P_{2}$,
on $S$ there is a regular arc which is totally-embedded in $S$ with
$P_{1}$ and $P_{2}$ being its end points. Examples of connected
surfaces are planes, ellipsoids and cylinders.

$\bullet$ A closed surface is a simple surface with no open edges;
examples of closed surfaces are spheres and ellipsoids.

$\bullet$ A bounded surface is a surface that can be contained entirely
in a sphere of a finite radius.

$\bullet$ A compact surface is a simple surface which is bounded
and closed like a torus or a Klein bottle.

$\bullet$ If $f$ is a differentiable regular mapping from a surface
$S$ to a surface $\bar{S}$, then if $S$ is compact then $\bar{S}$
is compact.

$\bullet$ If $S_{1}$ and $S_{2}$ are two simple surfaces where
$S_{1}$ is connected and $S_{2}$ is closed and contained in $S_{1}$,
then the two surfaces are equal as point sets. As a result, a simple
closed surface cannot be a proper subset of a simple connected surface.

$\bullet$ An orientable surface is a simple surface over which a
continuously-varying normal vector can be defined. Hence, spheres,
cylinders and tori are orientable surfaces while the Mobius strip
is a non-orientable surface since a normal vector moved continuously
around the strip from a given point will return to the point in the
opposite direction. An orientable surface which is connected can be
oriented in only one of two possible ways.

$\bullet$ An oriented surface is an orientable surface over which
the direction of the normal vector is determined.

$\bullet$ An elementary surface is a simple surface which possesses
a single coordinate patch basis, and hence it is an orientable surface
which can be mapped bicontinuously to an open set in the plane. Examples
of elementary surfaces are planes, cones and elliptic paraboloids.

$\bullet$ A surface that can be flattened into a plane by unfolding
without local distortion by compression or stretching is called developable
surface.\footnote{It is called developable because it can be developed into a plane
by rolling the surface out on a plane without compression or stretching.} A characteristic feature of a developable surface is that, like the
plane, its Gaussian curvature (see $\S$ \ref{subGaussianCurvature})
is zero at every point on the surface.

$\bullet$ A topological property of a surface is a property which
is invariant with respect to injective bicontinuous mappings. An example
of a topological property is compactness.

$\bullet$ A differentiable regular mapping from a surface $S$ to
a surface $\bar{S}$ is called conformal if it preserves angles between
oriented intersecting curves on the surface. The mapping is described
as direct if it preserves the sense of the angles and inverse if it
reverses it.

$\bullet$ Technically, the mapping is conformal if there is a function
$\lambda(u,v)>0$ that applies to all patches on the surface such
that $a_{\alpha\beta}=\lambda\bar{a}_{\alpha\beta}\,\,(\alpha,\beta=1,2)$
where the unbarred and barred indexed $a$ are the coefficients of
the surface covariant metric tensor in $S$ and $\bar{S}$ respectively.\footnote{This may be stated by saying that the first fundamental forms of the
two surfaces are proportional but proportionality should not be understood
to mean that $\lambda$ is constant.} An example of conformal mapping is the stereographic projection from
the Riemann sphere to a plane.

$\bullet$ An isometry or isometric mapping is a one-to-one mapping
from a surface $S$ to a surface $\bar{S}$ that preserves distances,
hence any regular arc in $S$ is mapped onto an arc in $\bar{S}$
with equal length. The two surfaces $S$ and $\bar{S}$ are described
as isometric surfaces. An example of isometric mapping is the deformation
of a rectangular plane sheet into a cylinder with no local distortion
by compression or stretching and hence the two surfaces are isometric
since all distances are preserved.

$\bullet$ Isometry is a symmetric relation and hence the inverse
of an isometric mapping is an isometric mapping, that is if $f$ is
an isometry from $S$ to $\bar{S}$, then $f^{-1}$ is an isometry
from $\bar{S}$ to $S$.

$\bullet$ An injective mapping from a surface $S$ onto a surface
$\bar{S}$ is an isometry \textit{iff} the coefficients of the first
fundamental form (see $\S$ \ref{subFirstFundamentalForm}) for any
patch on $S$ are identical to the coefficients of the first fundamental
form of its image on $\bar{S}$, that is
\begin{equation}
E=\bar{E}\,\,\,\,\,\,\,\,\,\,\,\,\,\,\,\,\,\,\,F=\bar{F}\,\,\,\,\,\,\,\,\,\,\,\,\,\,\,\,\,\,\,\,\,G=\bar{G}
\end{equation}
where the unbarred and barred $E,F,G$ are the coefficients of the
first fundamental form in the two surfaces.\footnote{As we will see, the coefficients of the first fundamental form are
the same as the coefficients of the surface covariant metric tensor,
that is:
\begin{equation}
a_{11}=E\,\,\,\,\,\,\,\,\,\,\,\,\,\,\,\,a_{12}=a_{21}=F\,\,\,\,\,\,\,\,\,\,\,\,\,\,\,\,\,a_{22}=G
\end{equation}
}

$\bullet$ \label{Llocalisometry}The mapping that preserves distances
but it is not injective is described as local isometry. The last point
also applies to local isometry.

$\bullet$ Since intrinsic properties are dependent only on the coefficients
of the first fundamental form of the surface, an intrinsic property
of the surface is invariant with respect to isometric mappings.

$\bullet$ As a consequence of the equality of corresponding lengths
of two isometric surfaces, the corresponding angles are also equal.
However, the reverse is not valid, that is a mapping that attains
the equality of corresponding angles (refer to conformal mapping above)
does not necessarily ensue the equality of corresponding lengths.

$\bullet$ Similarly, isometric mapping preserves areas of mapped
surfaces since it preserves lengths and angles.

$\bullet$ Isometric mapping is more restrictive than conformal mapping.
In fact conformal mapping can be set up between any two surfaces and
in many different ways but this is not always possible for isometric
mapping.

$\bullet$ A surface generated by the collection of all the tangent
lines to a given space curve is called the ``tangent surface'' of
the curve while the tangent lines are called the generators of the
surface.\footnote{These tangent lines are also called the rulings of the surface.}
The tangent surface of a curve may be demonstrated visually by a taut
flexible string connected to the curve where it scans the surface
while being directed tangentially at each point of the curve at its
base.\footnote{The ``tangent surface'' of a \textit{curve} should not be confused
with the aforementioned ``tangent plane'' of a \textit{surface}
at a given point. Also the taut string visualization should extend
to both directions to give the full extent of the tangent surface.}

$\bullet$ A ``branch'' of the tangent surface of a curve $C$ at
a given point $P$ on the curve refers to the tangent line of $C$
at $P$.

$\bullet$ If $C_{e}$ is a space curve with a tangent surface $S_{T}$
and $C_{i}$ is a curve embedded in $S_{T}$ and it is orthogonal
to all the tangent lines of $C_{e}$ at their intersection points,
then $C_{i}$ is called an involute of $C_{e}$ while $C_{e}$ is
called an evolute of $C_{i}$.

\pagebreak{}

\section{Curves in Space\label{subSpaceCurves}}

$\bullet$ ``Space'' in this title is general and hence it includes
surface, as explained in the Preface.

$\bullet$ Let have a space curve of class $C^{2}$ in a 3D Riemannian
manifold with a given metric $g_{ij}$ $(i,j=1,2,3)$. The curve is
parameterized by $s$ which is the distance along the curve starting
from an arbitrarily-chosen initial point on the curve, $P_{0}$.\footnote{We choose to parameterize the curve by $s$ to have simpler formulae.
Other formulae based on a more general parameterization will also
be given. We also use a mix of tensor and symbolic notations as each
has certain advantages and to familiarize the reader with both notations
as different authors use different notations.} The curve can therefore be represented by:
\begin{equation}
x^{i}=x^{i}(s)\,\,\,\,\,\,\,\,\,\,\,\,\,\,\,\,\,\,\,\,\left(i=1,2,3\right)
\end{equation}
where the indexed $x$ represent the space coordinates.\footnote{Here, we have $\mathbf{r}(s)=x^{i}(s)\mathbf{e}_{i}$ where $\mathbf{r}$
is the spatial representation of the space curve and $\mathbf{e}_{i}$
are the space basis vectors.}

\begin{figure}[h]
\centering~~~~~~~\includegraphics{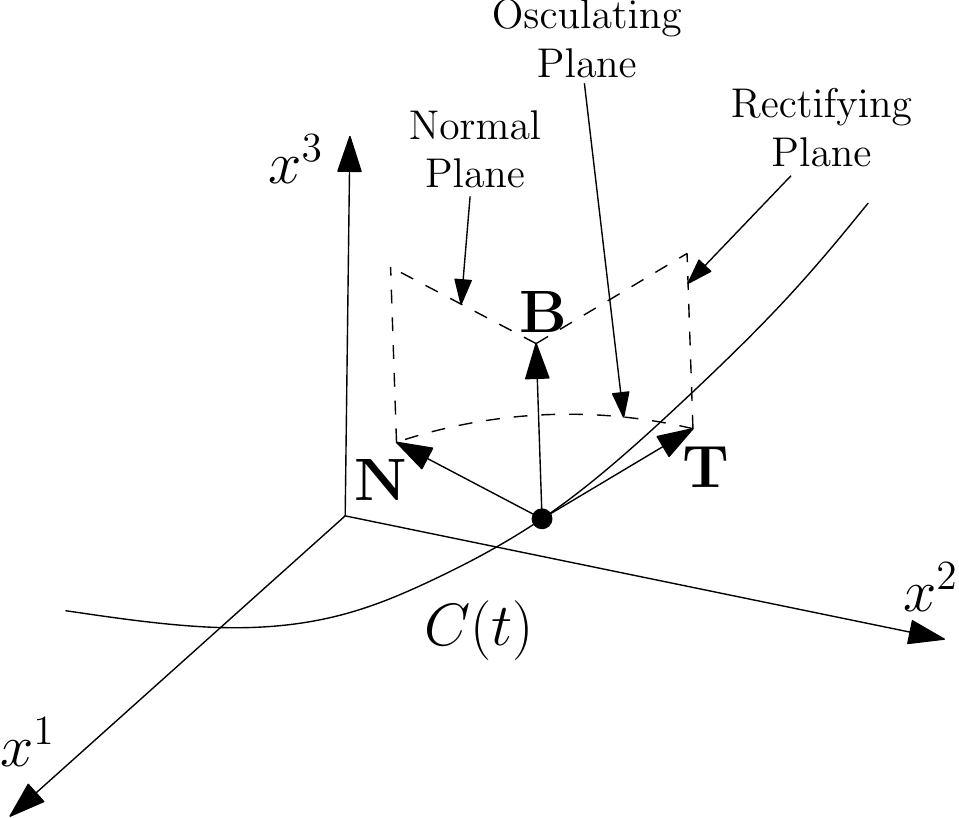}\label{FigFrenetFrame}\caption{The vectors $\mathbf{T},\mathbf{N},\mathbf{B}$ and their associated
planes.}
\end{figure}

$\bullet$ Three mutually perpendicular vectors each of unit length
can be defined at each point $P$ with non-zero curvature of the above-described
space curve: tangent $\mathbf{T}$, normal $\mathbf{N}$ and binormal
$\mathbf{B}$ (see Figure \ref{FigFrenetFrame}). These vectors can
serve as a moving coordinate system for the space.

$\bullet$ The unit vector tangent to the curve at point $P$ on the
curve is given by:\footnote{For simplicity, we employ Cartesian coordinates and hence ordinary
derivatives (i.e. $\frac{d}{ds}$ and $\frac{d}{dt}$) are used in
this and the following formulae. For general curvilinear coordinates,
these ordinary derivatives should be replaced by absolute derivatives
(i.e. $\frac{\delta}{\delta s}$ and $\frac{\delta}{\delta t}$) along
the curves. }
\begin{equation}
\left[\mathbf{T}\right]^{i}=T^{i}=\frac{dx^{i}}{ds}
\end{equation}

$\bullet$ For a $t$-parameterized curve, where $t$ is not necessarily
the arc length, the tangent vector is given by:
\begin{equation}
\mathbf{T}=\frac{\dot{\mathbf{r}}(t)}{\left|\mathbf{\dot{\mathbf{r}}}(t)\right|}
\end{equation}
where the overdot represents differentiation with respect to $t$.

$\bullet$ The unit vector normal\footnote{This is also called the unit principal normal vector. This vector
is defined only on points of the curve where the curvature $\kappa\ne0$.} to the tangent $T^{i}$, and hence to the curve, at the point $P$
is given by:
\begin{equation}
\left[\mathbf{N}\right]^{i}=N^{i}=\frac{\frac{dT^{i}}{ds}}{\left|\frac{dT^{i}}{ds}\right|}=\frac{1}{\kappa}\frac{dT^{i}}{ds}\label{eqNi}
\end{equation}
where $\kappa$ is a scalar called the ``curvature'' of the curve
at the point $P$ and is defined, according to the normalization condition,
by:\footnote{For general curvilinear coordinates, the formula becomes: $\kappa=\sqrt{g_{ij}\frac{\delta T^{i}}{\delta s}\frac{\delta T^{j}}{\delta s}}$
where $g_{ij}$ is the space covariant metric tensor.}

\begin{equation}
\kappa=\sqrt{\frac{dT^{i}}{ds}\frac{dT^{i}}{ds}}\label{eqkapdef}
\end{equation}

$\bullet$ For a $t$-parameterized curve, where $t$ is not necessarily
the arc length, the principal normal vector is given by:
\begin{equation}
\mathbf{N}=\frac{\dot{\mathbf{r}}(t)\times\left(\mathbf{\ddot{r}}(t)\times\dot{\mathbf{r}}(t)\right)}{\left|\dot{\mathbf{r}}(t)\right|\left|\mathbf{\ddot{r}}(t)\times\dot{\mathbf{r}}(t)\right|}
\end{equation}
where the overdot represents differentiation with respect to $t$.

$\bullet$ The binormal unit vector is defined as:
\begin{equation}
\left[\mathbf{B}\right]^{i}=B^{i}=\frac{1}{\tau}\left(\kappa T^{i}+\frac{dN^{i}}{ds}\right)
\end{equation}
which is a linear combination of two vectors both of which are perpendicular
to $N^{i}$ and hence it is perpendicular to $N^{i}$. In the last
equation, the normalization scalar factor $\tau$ is the ``torsion''
whose sign is chosen to make $T^{i},N^{i}$ and $B^{i}$ a right handed
triad vectors satisfying the condition:\footnote{Some authors reverse the sign in the definition of $\tau$ and this
reversal affects the signs in the forthcoming Frenet-Serret formulae
(see $\S$ \ref{subRelationCurveBasis}). The convention that we follow
in these notes has certain advantages.}
\begin{equation}
\underline{\epsilon}_{ijk}T^{i}N^{j}B^{k}=1
\end{equation}

$\bullet$ For a $t$-parameterized curve, where $t$ is not necessarily
the arc length, the binormal vector is given by:
\begin{equation}
\mathbf{B}=\frac{\dot{\mathbf{r}}(t)\times\mathbf{\ddot{r}}(t)}{\left|\dot{\mathbf{r}}(t)\times\mathbf{\ddot{r}}(t)\right|}
\end{equation}
where the overdot represents differentiation with respect to $t$.

$\bullet$ Apart from making $\mathbf{T},\mathbf{N}$ and $\mathbf{B}$
a right handed system, there is a geometric significance for the sign
of the torsion as it affects the orientation of the space curve.

$\bullet$ At any point on the space curve, the triad $T^{i},N^{i}$
and $B^{i}$ represents a mutually perpendicular right handed system
fulfilling the condition:
\begin{equation}
B^{i}=\left[\mathbf{T}\times\mathbf{N}\right]^{i}=\underline{\epsilon}^{ijk}T_{j}N_{k}
\end{equation}

$\bullet$ Since the triad $\mathbf{T},\mathbf{N}$ and $\mathbf{B}$
are mutually perpendicular, they satisfy the condition:
\begin{equation}
\mathbf{T}\cdot\mathbf{N}=\mathbf{T}\cdot\mathbf{B}=\mathbf{N}\cdot\mathbf{B}=0
\end{equation}
Also, because they are unit vectors we have:
\begin{equation}
\mathbf{T}\cdot\mathbf{T}=\mathbf{N}\cdot\mathbf{N}=\mathbf{B}\cdot\mathbf{B}=1
\end{equation}

$\bullet$ The tangent line of a curve $C$ at a given point $P$
on the curve is a straight line passing through $P$ and is parallel
to the tangent vector, $\mathbf{T}$, of $C$ at $P$. The principal
normal line of a curve $C$ at a given point $P$ on the curve is
a straight line passing through $P$ and is parallel to the principal
normal vector, $\mathbf{N}$, of $C$ at $P$. The binormal line of
a curve $C$ at a given point $P$ on the curve is a straight line
passing through $P$ and is parallel to the binormal vector, $\mathbf{B}$,
of $C$ at $P$.

$\bullet$ Following the last point, the equations of the three lines
can be given by the following generic form:
\begin{equation}
\mathbf{r}=\mathbf{r}_{P}+k\mathbf{V}_{P}\,\,\,\,\,\,\,\,\,\,\,\,\,\,\,\,\,\,\,\,\,\,\text{(\ensuremath{-\infty<k<\infty})}
\end{equation}
where $\mathbf{r}$ is the position vector of an arbitrary point on
the line, $\mathbf{r}_{P}$ is the position vector of the point $P$,
$k$ is a real variable and the vector $\mathbf{V}_{P}$ is the vector
corresponding to the particular line, that is $\mathbf{V}_{P}\equiv\mathbf{T}$
for the tangent line, $\mathbf{V}_{P}\equiv\mathbf{N}$ for the normal
line, and $\mathbf{V}_{P}\equiv\mathbf{B}$ for the binormal line.

$\bullet$ At any regular point $P$ on the space curve, the triad
$T^{i},N^{i}$ and $B^{i}$ define three mutually-perpendicular planes
where each one of these planes passes through the point $P$ and is
formed by a linear combination of two of these vectors in turn. These
planes are: the ``osculating plane'' which is the span of $T^{i}$
and $N^{i}$, the ``rectifying plane'' which is the span of $T^{i}$
and $B^{i}$, and the ``normal plane'' which is the span of $N^{i}$
and $B^{i}$ and is orthogonal to the curve at $P$ (see Figure \ref{FigFrenetFrame}).

$\bullet$ Following the last point, the equations of the three planes
can be given by the following generic form:
\begin{equation}
\left(\mathbf{r}-\mathbf{r}_{P}\right)\cdot\mathbf{V}_{P}=0
\end{equation}
where $\mathbf{r}$ is the position vector of an arbitrary point on
the plane, $\mathbf{r}_{P}$ is the position vector of the point $P$,
and where for each plane the vector $\mathbf{V}_{P}$ is the perpendicular
vector to the plane at $P$, that is $\mathbf{V}_{P}\equiv\mathbf{B}$
for the osculating plane, $\mathbf{V}_{P}\equiv\mathbf{N}$ for the
rectifying plane, and $\mathbf{V}_{P}\equiv\mathbf{T}$ for the normal
plane.

$\bullet$ Following the style of the definition of the tangent line
as a limit of a secant line (see $\S$ \ref{subCurves}), the osculating
plane may be defined as the limiting position of a plane passing through
$P$ and two other points on the curve one on each side of $P$ as
the two points converge simultaneously along the curve to $P$.

$\bullet$ The positive sense of a parameterized curve, which corresponds
to the direction in which the parameter increases and hence defines
the orientation of the curve, can be determined in two opposite ways.
While the sense of the tangent $\mathbf{T}$ and the binormal $\mathbf{B}$
is dependent on the curve orientation and hence they are in opposite
directions in these two ways, the principal normal $\mathbf{N}$ is
the same as it remains parallel to the normal plane in the direction
in which the curve is turning.

\subsection{Curvature and Torsion of Space Curves\label{subCurvatureTorsion}}

$\bullet$ The curvature and torsion of space curves may also be called
the first and second curvatures respectively, and hence a twisted
curve with non-vanishing curvature and non-vanishing torsion is described
as double-curvature curve. The expression $\sqrt{\left(ds_{\mathbf{T}}\right)^{2}+\left(ds_{\mathbf{B}}\right)^{2}}$,
where $ds_{\mathbf{T}}$ and $ds_{\mathbf{B}}$ are respectively the
line element components in the tangent and binormal directions, may
be described as the total or the third curvature of the curve.\footnote{The ``total curvature'' is also used for surfaces (see $\S$ \ref{subPrincipalCurvatures}
and \ref{subGaussBonnetTheorem}) but the meaning is different.}

$\bullet$ The equation of Lancret states that:
\begin{equation}
\left(ds_{\mathbf{N}}\right)^{2}=\left(ds_{\mathbf{T}}\right)^{2}+\left(ds_{\mathbf{B}}\right)^{2}
\end{equation}
where $ds_{\mathbf{N}}$ is the line element component in the normal
direction.

$\bullet$ According to the fundamental theorem of space curves in
differential geometry, a space curve is completely determined by its
curvature and torsion. More technically, given a real interval $I\text{\,}\text{\ensuremath{\subseteq}}\,\mathbb{R}$
and two differentiable real functions: $\kappa(s)>0$ and $\tau(s)$
where $s\,\in\,I$, there is a uniquely defined parameterized regular
space curve $C(s)$: $I\rightarrow\mathbb{R}^{3}$ of class $C^{2}$
with $\kappa(s)$ and $\tau(s)$ being the curvature and torsion of
$C$ respectively and $s$ is its arc length. Hence, any other curve
meeting these conditions will be different from $C$ only by a rigid
motion transformation (translation and rotation) which determines
its position and orientation in space.\footnote{In rigid motion transformations, which may also be called Euclidean
motion, the distance between any two points on the image is the same
as the distance between the corresponding points on the inverse image.
Hence, rigid motion transformation is a form of isometric mapping.}

$\bullet$ On the other hand, any curve with the properties given
in the last point possesses uniquely defined $\kappa(s)$ and $\tau(s)$.

$\bullet$ From the previous points, the fundamental theorem of space
curves provides the existence and uniqueness conditions for curves.

$\bullet$ The equations: $\kappa=\kappa(s)$ and $\tau=\tau(s)$,
where $s$ is the arc length, are called the intrinsic or natural
equations of the curve.

$\bullet$ The curvature and torsion are invariants of the space curve
and hence they do not depend on the employed coordinate system or
the type of parameterization.

$\bullet$ While the curvature is always non-negative ($\kappa\ge0$),
as it represents the magnitude of a vector according to the above-stated
definition (see e.g. Eqs. \ref{eqNi} and \ref{eqkapdef}), the torsion
can be negative as well as zero or positive.\footnote{Some authors define the curvature vector (see $\S$ \ref{subCurvatureVector})
and the principal normal vector of the curve in such a way that it
is possible for the curvature to be negative.}

$\bullet$ The following are some examples of the curvature and torsion
of a number of commonly-occurring simple curves:

(A) Straight line: $\kappa=0$ and $\tau=0$.

(B) Circle of radius $R$: $\kappa=\frac{1}{R}$ and $\tau=0$.\footnote{Hence, the radius of curvature (see $\S$ \ref{subCurvature}) of
a circle is its own radius}

(C) Helix parameterized by $\mathbf{r}(t)=\left(a\cos(t),a\sin(t),bt\right)$:
$\kappa=\frac{a}{a^{2}+b^{2}}$ and $\tau=\frac{b}{a^{2}+b^{2}}$.\footnote{It can be shown that a space curve of class $C^{3}$ with non-vanishing
curvature is a helix \textit{iff} the ratio of its torsion to curvature
is constant. Geometrically, the helix is characterized by having a
tangent vector that forms a constant angle with a specified direction
which is the direction defined by its axis of rotation (the circle
can be regarded as a degenerate helix).}

In these three examples, the curvature and torsion are constants along
the whole curve. However, in general the curvature and torsion of
space curves are position dependent and hence they vary from point
to point.

$\bullet$ \label{1Dinhabitant}Following the example of 2D surfaces,
a 1D inhabitant of a space curve can detect all the properties related
to the arc length. Hence, the curvature and torsion, $\kappa$ and
$\tau$, of the curve are extrinsic properties for such a 1D inhabitant.
This fact may be expressed by saying that curves are intrinsically
Euclidean, and hence their Riemann-Christoffel curvature tensor vanishes
identically and they naturally admit 1D Cartesian systems represented
by their natural parameterization of arc length.\footnote{This should be obvious when we consider that any curve can be mapped
isometrically to a straight line where both are naturally parameterized
by arc length.} Another demonstration of their intrinsic 1D nature is represented
by the Frenet-Serret formulae (see $\S$ \ref{subRelationCurveBasis}).

$\bullet$ Some authors resemble the role of $\kappa$ and $\tau$
in curve theory to the surface curvature tensor $b_{\alpha\beta}$
in surface theory and describe $\kappa$ and $\tau$ as the curve
theoretic analogues of the $b_{\alpha\beta}$ in surface theory. In
another context, $\kappa$ and $\tau$ may be contrasted with the
first and second fundamental forms of surfaces in their roles in defining
the curve and surface in the fundamental theorems of these structures.

$\bullet$ From the first and the last of the Frenet-Serret formulae
(Eq. \ref{eqFrenet}), we have:
\begin{equation}
\left|\kappa\tau\right|=\left|\mathbf{T}'\cdot\mathbf{B}'\right|
\end{equation}
where the prime stands for derivative with respect to a natural parameter
$s$ of the curve.

\subsubsection{Curvature\label{subCurvature}}

$\bullet$ The curvature $\kappa$ of a space curve is a measure of
how much the curve bends as it progresses in the tangent direction
at a particular point. The curvature represents the magnitude of the
rate of change of the direction of the tangent vector with respect
to the arc length and hence it is a measure for the departure of the
curve from the orientation of the straight line passing through that
point and oriented in the tangent direction. Consequently, the curvature
vanishes identically for straight lines.\footnote{Having an identically vanishing curvature is a necessary and sufficient
condition for a curve of class $C^{2}$ to be straight line.}

$\bullet$ From the first of the Frenet-Serret formulae (see Eq. \ref{eqFrenet})
and the fact that:
\begin{equation}
\left(\mathbf{N}\cdot\mathbf{T}\right)'=\left(0\right)'=0\,\,\,\,\,\,\,\,\,\,\,\,\,\,\,\,\,\,\,\Rightarrow\,\,\,\,\,\,\,\,\,\,\,\,\,\,\,\,\,\,\,\mathbf{N}\cdot\mathbf{T}'=-\mathbf{N}'\cdot\mathbf{T}
\end{equation}
the curvature $\kappa$ can be expressed as:
\begin{equation}
\kappa=\mathbf{N}\cdot\mathbf{T}'=-\mathbf{N}'\cdot\mathbf{T}
\end{equation}
where the prime represents differentiation with respect to the arc
length $s$ of the curve.

$\bullet$ The ``radius of curvature'', which is the radius of the
osculating circle (see $\S$ \ref{subOsculatingCircle}), is defined
at each point of the space curve for which $\kappa\ne0$ as the reciprocal
of the curvature, i.e. $R_{\kappa}=\frac{1}{\kappa}$.\footnote{A different way for introducing these concepts, which is followed
by some authors, is to define first the radius of curvature as the
reciprocal of the magnitude of the acceleration vector, that is $R_{\kappa}=\frac{1}{\left|\mathbf{r}''(s)\right|}$
where $\mathbf{r}$ is the spatial representation of an $s$-parameterized
curve; the curvature is then defined as the reciprocal of the radius
of curvature. Hence, the radius of curvature may be described as the
reciprocal of the norm of the acceleration vector where acceleration
means the second derivative of the curve.}

$\bullet$ There may be an advantage in using the concept of ``curvature''
as the principal concept instead of ``radius of curvature'', that
is the curvature is defined at all regular points while the radius
of curvature is defined only at the regular points with non-vanishing
curvature.

$\bullet$ As indicated earlier, if $C$ is a space curve of class
$C^{2}$ which is defined on a real interval $I$$\subseteq\mathbb{R}$
and is parameterized by arc length $s\in I$, that is $C(s):I\rightarrow\mathbb{R}^{3}$,
then the curvature of $C$ at a given point $P$ on the curve is defined
by: $\kappa=|\mathbf{r}''|$ where $\mathbf{r}(s)$ is the spatial
representation of the curve, the double prime represents the second
derivative with respect to $s$, and $\mathbf{r}''$ is evaluated
at $P$.\footnote{There should be no confusion in the present and forthcoming points
between the bare $C$ symbol and the superscripted $C$ symbol as
the bare $C$ symbolizes a curve while the superscripted $C$ stands
for the differentiability condition as explained earlier.}

$\bullet$ For a space curve represented parametrically by $\mathbf{r}(t)$,
where $t$ is a general parameter not necessarily a natural parameter,
we have:
\begin{equation}
\kappa=\frac{\left|\mathbf{\dot{T}}\right|}{\left|\mathbf{\dot{r}}\right|}=\frac{\left|\mathbf{\dot{r}}\times\mathbf{\ddot{r}}\right|}{\left|\mathbf{\dot{r}}\right|^{3}}=\frac{\sqrt{\left(\mathbf{\dot{r}}\cdot\mathbf{\dot{r}}\right)\left(\mathbf{\ddot{r}}\cdot\mathbf{\ddot{r}}\right)-\left(\mathbf{\dot{r}}\cdot\mathbf{\ddot{r}}\right)^{2}}}{\left(\mathbf{\dot{r}}\cdot\mathbf{\dot{r}}\right)^{3/2}}
\end{equation}
where all the quantities, which are functions of $t$, are evaluated
at a given point corresponding to a given value of $t$, and the overdot
represents derivative with respect to $t$.

$\bullet$ All surface curves passing through a point $P$ on a surface
$S$ and have the same osculating plane at $P$ have identical curvature
$\kappa$ at $P$ if the osculating plane is not tangent to $S$ at
$P$.

\subsubsection{Torsion}

$\bullet$ The torsion $\tau$ represents the rate of change of the
osculating plane, and hence it quantifies the twisting, in magnitude
and sense, of the space curve out of the plane of curvature and its
deviation from being a plane curve. The torsion therefore vanishes
identically for plane curves.\footnote{Having an identically vanishing torsion is a necessary and sufficient
condition for a curve of class $C^{2}$ to be plane curve. }

$\bullet$ If $C$ is a space curve of class $C^{2}$ which is defined
on a real interval $I$$\subseteq\mathbb{R}$ and it is parameterized
by arc length $s\in I$, that is $C(s):I\rightarrow\mathbb{R}^{3}$,
then the torsion of $C$ at a given point $P$ on the curve is given
by:\footnote{This can be obtained from the second of the Frenet-Serret formulae
(see Eq. \ref{eqFrenet}) by the dot product of both sides with $\mathbf{B}$.}
\begin{equation}
\tau=\mathbf{N}'\text{·}\mathbf{B}
\end{equation}
where $\mathbf{N}'$ and $\mathbf{B}$ are evaluated at $P$ and the
prime represents differentiation with respect to $s$.

$\bullet$ For a space curve represented parametrically by $\mathbf{r}(t)$,
where $t$ is a general parameter not necessarily a natural parameter,
we have:
\begin{equation}
\tau=\frac{\mathbf{\dot{r}}\cdot\left(\mathbf{\ddot{r}}\times\mathbf{\dddot{r}}\right)}{\left|\mathbf{\dot{r}}\times\mathbf{\ddot{r}}\right|^{2}}=\frac{\mathbf{\dot{r}}\cdot\left(\mathbf{\ddot{r}}\times\mathbf{\dddot{r}}\right)}{\left(\mathbf{\dot{r}}\times\mathbf{\ddot{r}}\right)\cdot\left(\mathbf{\dot{r}}\times\mathbf{\ddot{r}}\right)}=\frac{\mathbf{\dot{r}}\cdot\left(\mathbf{\ddot{r}}\times\mathbf{\dddot{r}}\right)}{\left(\mathbf{\dot{r}}\cdot\mathbf{\dot{r}}\right)\left(\mathbf{\ddot{r}}\cdot\mathbf{\ddot{r}}\right)-\left(\mathbf{\dot{r}}\cdot\mathbf{\ddot{r}}\right)^{2}}
\end{equation}
where all the quantities, which are functions of $t$, are evaluated
at a given point $P$ corresponding to a given value of $t$, and
the overdot represents derivative with respect to $t$. The curve
should have non-vanishing curvature $\kappa$ at $P$.

$\bullet$ For general curvilinear coordinates, the torsion of an
$s$-parameterized curve is given in tensor notation by:
\begin{equation}
\tau=\epsilon^{ijk}T_{i}N_{j}\frac{\delta N_{k}}{\delta s}
\end{equation}

$\bullet$ For rectangular Cartesian coordinates, the torsion of a
$t$-parameterized curve is given in tensor notation by:\footnote{This is based on the formula in the previous point.}

\begin{equation}
\tau=\frac{\epsilon_{ijk}\dot{x}_{i}\ddot{x}_{j}\dddot{x}_{k}}{\kappa^{2}}
\end{equation}
where $\kappa$ is the curvature of the curve as defined previously.

$\bullet$ The ``radius of torsion'' is defined at each point of
a space curve for which $\tau\ne0$ as the absolute value of the reciprocal
of the torsion, i.e. $R_{\tau}=\left|\frac{1}{\tau}\right|$.\footnote{Some authors do not take the absolute value and hence the radius of
torsion can be negative.}

$\bullet$ The value of torsion is invariant under permissible coordinate
transformations. It is also independent of the nature of its parameterization
and the orientation of the curve which is determined by the sense
of increase of its parameter.

\subsection{Geodesic Torsion\label{subGeodesicTorsion}}

$\bullet$ Geodesic torsion, which is an attribute of a curve embedded
in a surface, is also known as the relative torsion.

$\bullet$ The geodesic torsion of a surface curve $C$ at a given
point $P$ is the torsion of the geodesic curve (see $\S$ \ref{subGeodesicCurves})
passing through $P$ in the tangent direction of $C$ at $P$.\footnote{As we will see in $\S$ \ref{subGeodesicCurves}, in the neighborhood
of a point $P$ on a smooth surface and for any specified direction
there is one and only one geodesic curve passing through $P$ in that
direction. }

$\bullet$ The geodesic torsion $\tau_{g}$ of a surface curve represented
spatially by $\mathbf{r}=\mathbf{r}(s)$ is given by the following
scalar triple product:
\begin{equation}
\tau_{g}=\mathbf{n}\cdot\left(\mathbf{n}'\times\mathbf{r}'\right)
\end{equation}
where $\mathbf{n}$ is the unit normal vector to the surface, the
primes represent differentiation with respect to a natural parameter
$s$, and all these quantities are to be evaluated at a given point
on the curve corresponding to the value of $\tau_{g}$.

$\bullet$ The geodesic torsion of a curve $C$ at a non-umbilical
(see $\S$ \ref{subUmbilicPoints}) point $P$ is given in terms of
the principal curvatures (see $\S$ \ref{subPrincipalCurvatures})
by:
\begin{equation}
\tau_{g}=\left(\kappa_{1}-\kappa_{2}\right)\sin\theta\cos\theta\label{eqtaugphi}
\end{equation}
where $\theta$ is the angle between the tangent vector $\mathbf{T}$
to the curve $C$ at $P$ and the first principal direction $\mathbf{d}_{1}$
(see Darboux frame in $\S$ \ref{subPrincipalCurvatures}).

$\bullet$ The geodesic torsion of a surface curve $C$ parameterized
by arc length $s$ at a given point $P$ is given by:\footnote{The curve should not be asymptotic (see $\S$ \ref{subAsymptoticDirections}).}
\begin{equation}
\tau_{g}=\tau-\frac{d\phi}{ds}
\end{equation}
where $\tau$ is the torsion of $C$ at $P$ and $\phi$ is the angle
between the unit vector $\mathbf{n}$ normal to the surface and the
principal normal $\mathbf{N}$ of $C$ at $P$ (i.e. $\phi=\arccos\left(\mathbf{n}\cdot\mathbf{N}\right)$).
This is known as the Bonnet formula. This formula demonstrates that
when $\mathbf{n}$ and $\mathbf{N}$ are collinear (i.e. $\phi=0$),
the geodesic torsion and the torsion are equal (i.e. $\tau_{g}=\tau$).
In this case, the geodesic curvature will vanish and the curve becomes
a geodesic.\footnote{When $\mathbf{n}$ and $\mathbf{N}$ are collinear, the geodesic component
of the curvature vector will vanish (see $\S$ \ref{subGeodesicCurves}).}

$\bullet$ On a line of curvature (see $\S$ \ref{subLineofCurvature}),
the geodesic torsion vanishes identically.

$\bullet$ The geodesic torsion of a surface curve $C$ at a given
point $P$ is zero \textit{iff} $C$ is a tangent to a line of curvature
at $P$.

$\bullet$ The geodesic torsions of two orthogonal surface curves
are equal in magnitude and opposite in sign.

\subsection{Relationship between Curve Basis Vectors and their Derivatives\label{subRelationCurveBasis}}

$\bullet$ The three basis vectors $\mathbf{T}$, $\mathbf{N}$ and
$\mathbf{B}$ (see $\S$ \ref{subSpaceCurves}) are connected to their
derivatives by the Frenet-Serret formulae\footnote{These are also called Frenet formulae. The sign of the terms involving
$\tau$ depends on the convention about the torsion and hence these
equations differ between different authors.} which are given in rectangular Cartesian coordinates by:
\begin{eqnarray}
\frac{dT^{i}}{ds} & = & \kappa N^{i}\nonumber \\
\frac{dN^{i}}{ds} & = & \tau B^{i}-\kappa T^{i}\label{eqFrenet}\\
\frac{dB^{i}}{ds} & = & -\tau N^{i}\nonumber
\end{eqnarray}

$\bullet$ The Frenet-Serret formulae can be cast in the following
matrix form using symbolic notation:\footnote{As seen, the coefficient matrix is anti-symmetric.}
\begin{equation}
\left[\begin{array}{c}
\mathbf{T}^{'}\\
\mathbf{N}^{'}\\
\mathbf{B}^{'}
\end{array}\right]=\left[\begin{array}{ccc}
0 & \kappa & 0\\
-\kappa & 0 & \tau\\
0 & -\tau & 0
\end{array}\right]\left[\begin{array}{c}
\mathbf{T}\\
\mathbf{N}\\
\mathbf{B}
\end{array}\right]
\end{equation}
where all the quantities in this equation are functions of arc length
$s$ and the prime represents derivative with respect to $s$ (i.e.
$\frac{d}{ds}$).

$\bullet$ The Frenet-Serret formulae can also be given in the following
form:
\begin{eqnarray}
\mathbf{T}^{'} & = & \mathbf{d}\times\mathbf{T}\nonumber \\
\mathbf{N}^{'} & = & \mathbf{d}\times\mathbf{N}\\
\mathbf{B}^{'} & = & \mathbf{d}\times\mathbf{B}\nonumber
\end{eqnarray}
where $\mathbf{d}$ is the ``Darboux vector'' which is given by:
\begin{equation}
\mathbf{d}=\tau\mathbf{T}+\kappa\mathbf{B}
\end{equation}
$\bullet$ The three equations in the last point may be abbreviated
in a single equation as:
\begin{equation}
\left(\mathbf{T}',\mathbf{N}',\mathbf{B}'\right)=\mathbf{d}\times\left(\mathbf{T},\mathbf{N},\mathbf{B}\right)
\end{equation}

$\bullet$ In general curvilinear coordinates, the Frenet-Serret formulae
are given in terms of the absolute derivatives of the three vectors,
that is:
\begin{equation}
\begin{aligned}\frac{\delta T^{i}}{\delta s} & =\frac{dT^{i}}{ds}+\Gamma_{jk}^{i}T^{j}\frac{dx^{k}}{ds}=\kappa N^{i}\\
\frac{\delta N^{i}}{\delta s} & =\frac{dN^{i}}{ds}+\Gamma_{jk}^{i}N^{j}\frac{dx^{k}}{ds}=\tau B^{i}-\kappa T^{i}\\
\frac{\delta B^{i}}{\delta s} & =\frac{dB^{i}}{ds}+\Gamma_{jk}^{i}B^{j}\frac{dx^{k}}{ds}=-\tau N^{i}
\end{aligned}
\end{equation}
where the indexed $x$ represent general spatial coordinates and $s$
is a natural parameter; the other symbols are as given earlier.

$\bullet$ The triad $\mathbf{T}$, $\mathbf{N}$ and $\mathbf{B}$
is called the Frenet frame which forms a set of orthonormal basis
vectors for $\mathbb{R}^{3}$. This frame serves as a moving orthogonal
coordinate system on the points of the curve.\footnote{This triad may also be called the Frenet trihedron or the moving trihedron
of the curve. This frame can suffer from problems or become undefined
at inflection points where $\frac{d\mathbf{T}}{ds}=\mathbf{0}$.}

$\bullet$ According to the fundamental theorem of space curves, which
is outlined previously in $\S$ \ref{subCurvatureTorsion}, a curve
does exist and it is unique \textit{iff} its curvature and torsion
as functions of arc length are given. Now, it is natural to expect
that such a solution can be obtained from the system of differential
equations given by Frenet-Serret formulae. However, such a solution
cannot be obtained by direct integration of these equations. More
elaborate methods (e.g. methods based on the Riccati equation for
reducing a system of simultaneous differential equations to a first
order differential equation) may be used to obtain the solution. Nevertheless,
a solution can be obtained by direct integration of the Frenet-Serret
formulae for plane curves where the torsion vanishes identically.
A solution by direct integration of the Frenet-Serret formulae can
also be obtained in simple cases such as when the curvature and torsion
are constants.

\subsection{Osculating Circle and Sphere\label{subOsculatingCircle}}

$\bullet$ At any point $P$ with non-zero curvature of a smooth space
curve $C$ an ``osculating circle''\footnote{It may also be called the circle of curvature or the kissing circle.}
can be defined, where this circle is characterized by:

(A) It is tangent to $C$ at $P$ (i.e. the circle and the curve have
a common tangent vector at $P$).

(B) It lies in the osculating plane.

(C) Its radius is $\frac{1}{\kappa}$ where $\kappa$ is the curvature
of $C$ at $P$.

(D) Its center $\mathbf{r}_{C}$ is at $\mathbf{r}_{C}=\mathbf{r}_{P}+\frac{\mathbf{N}}{\kappa}$
where $\mathbf{r}_{P}$ is the position vector of $P$ and $\mathbf{N}$
is the principal normal of $C$ at $P$.\footnote{In some cases the center may be defined geometrically but not analytically
when the second derivative of the curve is not defined at the given
point.}

$\bullet$ The center of curvature of a curve at a point on the curve
is defined as the center of the osculating circle at that point, as
given in the last point.

$\bullet$ For all points on a circle, the center of curvature is
the center of the circle, so the circle is its own osculating circle.

$\bullet$ The osculating circle provides a good approximation to
the curve in the neighborhood of its regular points.

$\bullet$ Following the manner of defining the tangent line to a
curve as a limit of the secant line (see $\S$ \ref{subCurves}),
the osculating circle to a curve at a point $P$ may be defined as
the limit of the circle passing through $P$ and two other points
on the curve one on each side of $P$ as these two points converge
to $P$ while staying on the curve.

$\bullet$ The ``osculating sphere'' of a curve $C$ at a point
$P$ may be defined as the limit of a sphere passing through $P$
and three neighboring points on the curve as these points converge
to $P$. The position of the center $C_{S}$ of the osculating sphere
at $P$, which is called the center of spherical curvature of $C$
at $P$, is given by:
\begin{equation}
\mathbf{r}_{S}=\mathbf{r}_{P}+\frac{1}{\kappa}\mathbf{N}-\frac{\kappa'}{\tau\kappa^{2}}\mathbf{B}=\mathbf{r}_{P}+R_{\kappa}\mathbf{N}+R_{\tau}R_{\kappa}'\mathbf{B}\label{eqCenterSC}
\end{equation}
where $\mathbf{r}_{S}$ and $\mathbf{r}_{P}$ are the position vectors
of $C_{S}$ and $P$, $\mathbf{B}$ and $\mathbf{N}$ are the binormal
and principal normal vectors, $\kappa$ and $\tau$ are the curvature
and torsion, $R_{\kappa}$ and $R_{\tau}$ are the radii of curvature
and torsion, and the prime represents derivative with respect to a
natural parameter of $C$. All these quantities belong to $C$ at
$P$ which should have non-vanishing curvature and torsion ($\ensuremath{\kappa,\tau\ne0}$).

$\bullet$ From Eq. \ref{eqCenterSC}, the radius of the osculating
sphere is given by:
\begin{equation}
\left|\mathbf{r}_{S}-\mathbf{r}_{P}\right|=\sqrt{R_{\kappa}^{2}+\left(R_{\tau}R_{\kappa}'\right)^{2}}
\end{equation}

\subsection{Parallelism and Parallel Propagation\label{subParallelism}}

$\bullet$ In flat spaces parallelism is an absolute property as it
is defined without reference to an external object. However, in Riemannian
spaces the idea of parallelism is defined in comparison to a prescribed
curve and hence it is different from the idea of parallelism in the
Euclidean sense.

$\bullet$ A vector field $A^{\alpha}$ is described as being parallel
along the surface curve $u^{\beta}=u^{\beta}(t)$ \textit{iff} its
absolute derivative (see $\S$ \ref{secTensorDifferentiation}) along
the curve vanishes, that is:
\begin{equation}
\frac{\delta A^{\alpha}}{\delta t}\equiv A_{\,\,;\beta}^{\alpha}\frac{du^{\beta}}{dt}\equiv\frac{dA^{\alpha}}{dt}+\Gamma_{\beta\gamma}^{\alpha}A^{\gamma}\frac{du^{\beta}}{dt}=0
\end{equation}
This means that the sufficient and necessary condition for a vector
field to be parallel along a surface curve is that the covariant derivative
of the field is normal to the surface (see $\S$ \ref{secTensorDifferentiation}).

$\bullet$ All the vectors of a field of parallel vectors have the
same constant magnitude.

$\bullet$ A field of absolutely parallel unit vectors on a surface
do exist \textit{iff} there is an isometric correspondence between
the plane and the surface.

$\bullet$ When two surfaces are tangent to each other along a given
curve $C$, then a vector field which is parallel along $C$ with
respect to one of these surfaces will also be parallel along $C$
with respect to the other surface.

$\bullet$ As a consequence of the definition of parallelism in Riemannian
space we have:

(A) A surface vector field parallelly propagated along a given curve
between two points $P_{1}$ and $P_{2}$ on the curve does not necessarily
coincide with another vector field parallelly propagated along another
curve connecting $P_{1}$ and $P_{2}$.

(B) Parallel propagation is path dependent, that is: given two points
$P_{1}$ and $P_{2}$, the vector obtained at $P_{2}$ by parallel
propagation of a vector from $P_{1}$ along a given curve $C$ connecting
$P_{1}$ to $P_{2}$ depends on the curve $C$.

(C) Starting from a given point $P$ on a closed surface curve $C$
enclosing a simply connected region on the surface, parallel propagation
of a vector field around $C$ starting from $P$ does not necessarily
result in the same vector field when arriving at $P$.\footnote{The angle between the initial and final vectors is a measure of the
Gaussian curvature on the surface.}

$\bullet$ When two non-trivial vectors experience parallel propagation
along a particular curve the angle between them stays constant.

\pagebreak{}

\section{Surfaces in Space\label{secSurfacesinSpace}}

$\bullet$ Here, we examine surfaces of class $C^{2}$ embedded in
a 3D Euclidean space using a Cartesian coordinate system ($x,y,z$)
for the most parts. Some notes are based on a more general Riemannian
space with a curvilinear coordinate system.

$\bullet$ Assuming a parametric representation for the surface, where
each one of the space coordinates ($x,y,z$) on the surface is a real
differentiable function of the two surface curvilinear coordinates
($u,v$), the position vector of a point $P$ on the surface as a
function of the surface curvilinear coordinates is given by:
\begin{equation}
\mathbf{r}(u,v)=x(u,v)\mathbf{e}_{1}+y(u,v)\mathbf{e}_{2}+z(u,v)\mathbf{e}_{3}
\end{equation}
where the indexed $\mathbf{e}$ are the Cartesian orthonormal basis
vectors in the three directions. It is also assumed that $\partial_{u}\mathbf{r}$
and $\partial_{v}\mathbf{r}$ are linearly independent and hence they
are not parallel or anti-parallel, that is:\footnote{This is a sufficient and necessary condition for the surface to be
``regular'' at the given point. The point is also described as ``regular'';
otherwise it is ``singular'' if the condition is violated. The surface
is regular on $\Omega$, a closed subset of $\mathbb{R}^{2}$, if
it is regular at each interior point of $\Omega$. The regularity
condition guarantees that the surface mapping is one-to-one and possesses
continuous inverse.}
\begin{equation}
\frac{\partial\mathbf{r}}{\partial u}\times\frac{\partial\mathbf{r}}{\partial v}\ne\mathbf{0}
\end{equation}

$\bullet$ To express the position vector of $P$ in tensor notation,
we re-label the space and surface coordinates as:
\begin{equation}
(x,y,z)\equiv(x^{1},x^{2},x^{3})\,\,\,\,\,\,\,\,\,\,\,\,\,\,\,\,\,\,\,\,\&\,\,\,\,\,\,\,\,\,\,\,\,\,\,\,\,(u,v)\equiv(u^{1},u^{2})
\end{equation}
and hence the position vector becomes:
\begin{equation}
\mathbf{r}(u^{1},u^{2})=x^{i}(u^{1},u^{2})\mathbf{e}_{i}
\end{equation}

$\bullet$ To define a surface grid serving as a curvilinear positioning
system for the surface, one of the coordinate variables is held fixed
in turn while the other is varied. Hence, each one of the following
two surface functions:
\begin{equation}
\mathbf{r}(u^{1},c_{2})\,\,\,\,\,\,\,\,\,\,\,\,\,\,\,\,\,\,\,\&\,\,\,\,\,\,\,\,\,\,\,\,\,\,\,\,\,\,\,\,\,\mathbf{r}(c_{1},u^{2})
\end{equation}
where $c_{1}$ and $c_{2}$ are given real constants, defines a coordinate
curve for the surface. These two curves meet at the common surface
point ($c_{1},c_{2}$). The grid is then generated by varying $c_{1}$
and $c_{2}$ uniformly to obtain coordinate curves at regular intervals.

$\bullet$ The surface coordinate curves of the above grid are orthogonal
\textit{iff} the surface metric tensor (see $\S$ \ref{subSurfaceMetric})
is diagonal everywhere on the surface.

$\bullet$ Corresponding to each one of the surface coordinate curves
in the above order, a tangent vector to the curve at a given point
on the curve is defined by:\footnote{That is $\mathbf{E}_{1}$ is tangent to the $\mathbf{r}(u^{1},c_{2})$
curve and $\mathbf{E}_{2}$ is tangent to the $\mathbf{r}(c_{1},u^{2})$
curve.}
\begin{equation}
\mathbf{E}_{\alpha}=\frac{\partial\mathbf{r}}{\partial u^{\alpha}}=\frac{\partial x^{i}}{\partial u^{\alpha}}\mathbf{e}_{i}=x_{\alpha}^{i}\mathbf{e}_{i}
\end{equation}
where the derivatives are evaluated at that point, and $\alpha=1,2$
and $i=1,2,3$. These tangent vectors serve as a set of basis vectors
for the surface, and for each given point on the surface they generate,
by their linear combination, any vector in the surface at that point.\footnote{This should be understood in an infinitesimal sense or, equivalently,
as a vector lying in the tangent plane of the surface at the given
point.} They also define, by their linear combination, a plane tangent to
the surface at that point.\footnote{The plane generated by the linear combination of $\mathbf{E}_{1}$
and $\mathbf{E}_{2}$ is the tangent space, $T_{P}S$, to the surface
at point $P$ as described earlier, and hence $\mathbf{E}_{1}(u_{P}^{1},u_{P}^{2})$
and $\mathbf{E}_{2}(u_{P}^{1},u_{P}^{2})$ form a basis for this space
where the subscript $P$ is a reference to the point $P$.} A normal vector to the surface at that point is then defined as the
cross product of these tangent basis vectors: $\mathbf{E}_{1}\times\mathbf{E}_{2}$.
This normal vector can be scaled by its magnitude to produce a unit
vector normal to the surface at that point:
\begin{equation}
\mathbf{n}=\frac{\mathbf{E}_{1}\times\mathbf{E}_{2}}{\left|\mathbf{E}_{1}\times\mathbf{E}_{2}\right|}=\frac{\mathbf{E}_{1}\times\mathbf{E}_{2}}{\sqrt{a}}
\end{equation}
where the vectors are labeled so that $\mathbf{E}_{1},\mathbf{E}_{2}$
and $\mathbf{n}$ form a right-handed system, and $a$ is the determinant
of the surface covariant metric tensor (see $\S$ \ref{subSurfaceMetric}).

$\bullet$ The surface basis vectors, $\mathbf{E}_{\alpha}$, are
symbolized in full tensor notation by:
\begin{equation}
\left[\mathbf{E}_{\alpha}\right]^{i}\equiv E_{\alpha}^{i}=\frac{\partial x^{i}}{\partial u^{\alpha}}=x_{\alpha}^{i}\,\,\,\,\,\,\,\,\,\,\,\,\,\,\,\,\,\,\text{(\ensuremath{i=1,2,3} and \ensuremath{\alpha=1,2})}
\end{equation}
and hence they can be regarded as 3D contravariant space vectors or
as 2D covariant surface vectors (refer to $\S$ \ref{subSurfaceMetric}
for further details).

$\bullet$ It can be shown that the covariant form of the unit vector
$\mathbf{n}$ normal to the surface is given by:
\begin{equation}
n_{i}=\frac{1}{2}\underline{\epsilon}^{\alpha\beta}\underline{\epsilon}_{ijk}x_{\alpha}^{j}x_{\beta}^{k}
\end{equation}
where $x_{\alpha}^{j}=\frac{\partial x^{j}}{\partial u^{\alpha}}$
and similarly for $x_{\beta}^{k}$. The implication of this equation,
which defines $\mathbf{n}$ in terms of the surface basis vectors
$x_{\alpha}^{j}$ and $x_{\beta}^{k}$, is that $\mathbf{n}$ is a
space vector which is independent of the choice of the surface coordinates
$u^{1},u^{2}$ in support of the geometric intuition.

$\bullet$ Since $\mathbf{n}$ is normal to the surface, we have:
\begin{equation}
g_{ij}n^{i}x_{\alpha}^{j}=0
\end{equation}
which is the statement, in tensor notation, that $\mathbf{n}$ is
orthogonal to every vector in the tangent space of the surface at
the given point.

$\bullet$ Although $\mathbf{E}_{1}$ and $\mathbf{E}_{2}$ are linearly
independent they are not necessarily orthogonal or of unit length.
However, they can be orthonormalized as follow:
\begin{equation}
\underline{\mathbf{E}}_{1}=\frac{\mathbf{E}_{1}}{\left|\mathbf{E}_{1}\right|}=\frac{\mathbf{E}_{1}}{\sqrt{a_{11}}}\,\,\,\,\,\,\,\,\,\,\,\,\,\,\,\,\,\,\,\,\,\,\,\,\,\,\,\,\,\,\,\,\underline{\mathbf{E}}_{2}=\frac{a_{11}\mathbf{E}_{2}-a_{12}\mathbf{E}_{1}}{\sqrt{a_{11}a}}
\end{equation}
where $a$ is the determinant of the surface covariant metric tensor
(see $\S$ \ref{subSurfaceMetric}), the indexed $a$ are the coefficients
of this tensor, and the underlined vectors are orthonormal basis vectors,
that is:
\begin{equation}
\underline{\mathbf{E}}_{1}\cdot\underline{\mathbf{E}}_{1}=1\,\,\,\,\,\,\,\,\,\,\,\,\,\,\,\,\,\,\,\,\,\,\,\,\,\,\,\,\,\,\,\,\underline{\mathbf{E}}_{2}\cdot\underline{\mathbf{E}}_{2}=1\,\,\,\,\,\,\,\,\,\,\,\,\,\,\,\,\,\,\,\,\,\,\,\,\,\,\,\,\,\,\,\,\underline{\mathbf{E}}_{1}\cdot\underline{\mathbf{E}}_{2}=0
\end{equation}

$\bullet$ The transformation rules from one curvilinear coordinate
system of the surface to another coordinate system, notated with unbarred
($u^{1},u^{2}$) and barred ($\bar{u}^{1},\bar{u}^{2}$) symbols respectively,
where
\begin{equation}
\begin{aligned} & u^{1}=u^{1}(\bar{u}^{1},\bar{u}^{2})\,\,\,\,\,\,\,\,\,\,\,\,\,\,\,\,\,\,\,\,\&\,\,\,\,\,\,\,\,\,\,\,\,\,\,\,\,\,\,\,\,u^{2}=u^{2}(\bar{u}^{1},\bar{u}^{2})\\
 & \bar{u}^{1}=\bar{u}^{1}(u^{1},u^{2})\,\,\,\,\,\,\,\,\,\,\,\,\,\,\,\,\,\,\,\,\&\,\,\,\,\,\,\,\,\,\,\,\,\,\,\,\,\,\,\,\,\bar{u}^{2}=\bar{u}^{2}(u^{1},u^{2})
\end{aligned}
\end{equation}
are similar to the general rules outlined in \cite{SochiTC1,SochiTC2}
for the transformation between coordinate systems in a general $n$D
space.

$\bullet$ Following a transformation from the unbarred surface coordinate
system to the barred surface coordinate system, the surface becomes
a function of the barred coordinates, and a new set of basis vectors
for the surface, which are the tangents to the coordinate curves of
the barred system, are defined by the following equations:
\begin{equation}
\begin{aligned}\mathbf{\bar{E}}_{1} & =\frac{\partial\mathbf{r}}{\partial\bar{u}^{1}}=\frac{\partial\mathbf{r}}{\partial u^{1}}\frac{\partial u^{1}}{\partial\bar{u}^{1}}+\frac{\partial\mathbf{r}}{\partial u^{2}}\frac{\partial u^{2}}{\partial\bar{u}^{1}}=\mathbf{E}_{1}\frac{\partial u^{1}}{\partial\bar{u}^{1}}+\mathbf{E}_{2}\frac{\partial u^{2}}{\partial\bar{u}^{1}}\\
\mathbf{\bar{E}}_{2} & =\frac{\partial\mathbf{r}}{\partial\bar{u}^{2}}=\frac{\partial\mathbf{r}}{\partial u^{1}}\frac{\partial u^{1}}{\partial\bar{u}^{2}}+\frac{\partial\mathbf{r}}{\partial u^{2}}\frac{\partial u^{2}}{\partial\bar{u}^{2}}=\mathbf{E}_{1}\frac{\partial u^{1}}{\partial\bar{u}^{2}}+\mathbf{E}_{2}\frac{\partial u^{2}}{\partial\bar{u}^{2}}
\end{aligned}
\end{equation}
These equations, which correlate the surface basis vectors in the
barred and unbarred surface curvilinear coordinate systems, can be
compactly presented in tensor notation as:
\begin{equation}
\frac{\partial x^{i}}{\partial\bar{u}^{\alpha}}=\frac{\partial x^{i}}{\partial u^{\beta}}\frac{\partial u^{\beta}}{\partial\bar{u}^{\alpha}}\,\,\,\,\,\,\,\,\,\,\,\,\,\,\,\,\,\,\text{(\ensuremath{i=1,2,3} and \ensuremath{\alpha,\beta=1,2})}
\end{equation}

$\bullet$ A set of contravariant basis vectors for the surface may
also be defined as the gradient of the surface curvilinear coordinates,
that is:
\begin{equation}
\mathbf{E}^{\alpha}=\nabla u^{\alpha}
\end{equation}
In tensor notation, this basis is given by:
\begin{equation}
\left[\mathbf{E}^{\alpha}\right]_{i}\equiv E_{i}^{\alpha}=\frac{\partial u^{\alpha}}{\partial x^{i}}=x_{i}^{\alpha}\,\,\,\,\,\,\,\,\,\,\,\,\,\,\,\,\,\,\text{(\ensuremath{i=1,2,3} and \ensuremath{\alpha=1,2})}
\end{equation}
Hence, they can be regarded as 2D contravariant surface vectors or
as 3D covariant space vectors.

$\bullet$ The contravariant and covariant forms of the surface basis
vectors, $\mathbf{E}^{\alpha}$ and $\mathbf{E}_{\alpha}$, are obtained
from each other by the index-shifting operator for the surface, that
is:
\begin{equation}
\mathbf{E}_{\alpha}=a_{\alpha\beta}\mathbf{E}^{\beta}\,\,\,\,\,\,\,\,\,\,\,\,\,\,\,\,\,\,\,\,\,\,\,\,\,\,\,\,\,\,\,\,\,\,\,\,\mathbf{E}^{\alpha}=a^{\alpha\beta}\mathbf{E}_{\beta}\,\,\,\,\,\,\,\,\,\,\,\,\,\,\,\,\,\,\,\,\,\,(\alpha,\beta=1,2)
\end{equation}
where the indexed $a$ are the covariant and contravariant forms of
the surface metric tensor (see $\S$ \ref{subSurfaceMetric}).

$\bullet$ The contravariant and covariant forms of the surface basis
vectors, $\mathbf{E}^{\alpha}$ and $\mathbf{E}_{\alpha}$, are reciprocal
systems and hence they satisfy the following relations:
\begin{equation}
\mathbf{E}_{\alpha}\cdot\mathbf{E}^{\beta}=\delta_{\alpha}^{\beta}\equiv a_{\alpha}^{\beta}\,\,\,\,\,\,\,\,\,\,\,\,\,\,\,\,\,\,\,\,\,\,\,\,\,\,\,\,\,\,\,\mathbf{E}^{\alpha}\cdot\mathbf{E}_{\beta}=\delta_{\beta}^{\alpha}\equiv a_{\beta}^{\alpha}\,\,\,\,\,\,\,\,\,\,\,\,\,\,\,\,\,\,\,\,(\alpha,\beta=1,2)
\end{equation}

\subsection{Surface Metric Tensor\label{subSurfaceMetric}}

$\bullet$ The surface metric tensor\footnote{In differential geometry, the surface metric tensor $a_{\alpha\beta}$
may be called the first groundform.} is an absolute, rank-2, $2\times2$ symmetric tensor.\footnote{The coefficients of the metric tensor are real numbers.}

$\bullet$ Following the example of the metric in general $n$D spaces,
as explained in \cite{SochiTC1,SochiTC2}, the surface metric tensor
of a 2D surface embedded in a 3D Euclidean flat space with metric
$g_{ij}=\delta_{ij}$ is given in its covariant form by:
\begin{equation}
a_{\alpha\beta}=\mathbf{E}_{\alpha}\cdot\mathbf{E}_{\beta}=\frac{\partial\mathbf{r}}{\partial u^{\alpha}}\cdot\frac{\partial\mathbf{r}}{\partial u^{\beta}}=\frac{\partial x^{i}}{\partial u^{\alpha}}\frac{\partial x^{i}}{\partial u^{\beta}}\label{eqSurMet1}
\end{equation}
where the indexed $x$ and $u$ are the space Cartesian coordinates
and the surface curvilinear coordinates respectively, and $i=1,2,3$
and $\alpha,\beta=1,2$.

$\bullet$ The surface and space metric tensors in a general Riemannian
space with general metric $g_{ij}$ are related by:
\begin{equation}
a_{\alpha\beta}=\mathbf{E}_{\alpha}\cdot\mathbf{E}_{\beta}=\frac{\partial\mathbf{r}}{\partial u^{\alpha}}\cdot\frac{\partial\mathbf{r}}{\partial u^{\beta}}=g_{ij}\frac{\partial x^{i}}{\partial u^{\alpha}}\frac{\partial x^{j}}{\partial u^{\beta}}=g_{ij}x_{\alpha}^{i}x_{\beta}^{j}\label{eqSurMet2}
\end{equation}
where $a_{\alpha\beta}$ and $g_{ij}$ are respectively the surface
and space covariant metric tensors, the indexed $x$ and $u$ are
the curvilinear coordinates of the space and surface respectively,
and $i,j=1,2,3$ and $\alpha,\beta=1,2$. It is noteworthy that Eq.
\ref{eqSurMet1} is a special instance of Eq. \ref{eqSurMet2} for
the case of a flat space with a Cartesian system where the space metric
is the unity tensor.

$\bullet$ Eq. \ref{eqSurMet2} is the fundamental relation that provides
the crucial link between the surface and its enveloping space. As
indicated before, the partial derivatives in this relations, $\frac{\partial x^{i}}{\partial u^{\alpha}}$
and $\frac{\partial x^{j}}{\partial u^{\beta}}$, may be considered
as contravariant rank-1 3D space tensors or as covariant rank-1 2D
surface tensors. A tensor like $\frac{\partial x^{i}}{\partial u^{\alpha}}$
is usually labeled as $x_{\alpha}^{i}$ where it represents two surface
vectors which are contravariantly-transformed with respect to the
three space coordinates $x^{i}$:
\begin{equation}
x_{1}^{i}=\left(\frac{\partial x^{1}}{\partial u^{1}},\frac{\partial x^{2}}{\partial u^{1}},\frac{\partial x^{3}}{\partial u^{1}}\right)\,\,\,\,\,\,\,\,\,\,\,\,\,\,\,\,\,\,x_{2}^{i}=\left(\frac{\partial x^{1}}{\partial u^{2}},\frac{\partial x^{2}}{\partial u^{2}},\frac{\partial x^{3}}{\partial u^{2}}\right)
\end{equation}
or three space vectors which are covariantly-transformed with respect
to the two surface coordinates $u^{\alpha}$:
\begin{equation}
x_{\alpha}^{1}=\left(\frac{\partial x^{1}}{\partial u^{1}},\frac{\partial x^{1}}{\partial u^{2}}\right)\,\,\,\,\,\,\,\,\,\,\,\,\,x_{\alpha}^{2}=\left(\frac{\partial x^{2}}{\partial u^{1}},\frac{\partial x^{2}}{\partial u^{2}}\right)\,\,\,\,\,\,\,\,\,\,\,\,\,x_{\alpha}^{3}=\left(\frac{\partial x^{3}}{\partial u^{1}},\frac{\partial x^{3}}{\partial u^{2}}\right)
\end{equation}

$\bullet$ Any surface vector $A^{\alpha}$ ($\alpha=1,2$), defined
as a linear combination of the surface basis vectors $\mathbf{E}_{1}$
and $\mathbf{E}_{2}$, can also be considered as a space vector $A^{i}$
($i=1,2,3$) where the two representations are linked through the
relation:
\begin{equation}
A^{i}=\frac{\partial x^{i}}{\partial u^{\alpha}}A^{\alpha}=x_{\alpha}^{i}A^{\alpha}\,\,\,\,\,\,\,\,\,\,\,\,\text{(\ensuremath{i=1,2,3} and \ensuremath{\alpha=1,2})}\label{eqAiAalpha}
\end{equation}
Now, since we have (see Eqs. \ref{eqSurMet2} and \ref{eqAiAalpha}):
\begin{equation}
a_{\alpha\beta}A^{\alpha}A^{\beta}=g_{ij}x_{\alpha}^{i}x_{\beta}^{j}A^{\alpha}A^{\beta}=g_{ij}x_{\alpha}^{i}A^{\alpha}x_{\beta}^{j}A^{\beta}=g_{ij}A^{i}A^{j}
\end{equation}
then the two representations are equivalent, that is they define a
vector of the same magnitude and direction.

$\bullet$ The surface basis vectors in their covariant and contravariant
forms, $x_{\alpha}^{i}$ and $x_{i}^{\alpha}$, and the unit vector
normal to the surface $n_{i}$ are linked by the following relation:

\begin{equation}
x_{\alpha}^{i}=\underline{\epsilon}^{ijk}\underline{\epsilon}_{\alpha\beta}x_{j}^{\beta}n_{k}
\end{equation}
This equation means that the given product (which looks like a vector
cross product) of the surface contravariant basis vector $x_{j}^{\beta}$
and the unit normal vector $n_{k}$ produces a surface covariant basis
vector $x_{\alpha}^{i}$ and hence it is perpendicular to both.\footnote{Being a surface basis vector implies orthogonality to the unit normal
vector, while being a covariant surface basis vector implies orthogonality
to the contravariant surface basis vector.}

$\bullet$ The contravariant form of the surface metric tensor is
defined as the inverse of the surface covariant metric tensor\footnote{Since the first fundamental form is positive definite (see $\S$ \ref{subFirstFundamentalForm}),
and hence $a>0$, the existence of an inverse is guaranteed.}, that is:
\begin{equation}
a^{\alpha\gamma}\,a_{\gamma\beta}=\delta_{\,\,\beta}^{\alpha}\,\,\,\,\,\,\,\,\,\,\,\,\,\,\,\,\,\,\,a_{\alpha\gamma}\,a^{\gamma\beta}=\delta_{\alpha}^{\,\,\beta}
\end{equation}

$\bullet$ Similar to the metric tensor in general $n$D spaces, the
covariant and contravariant forms of the surface metric tensor, $a_{\alpha\beta}$
and $a^{\alpha\beta}$, are used for lowering and raising indices
and related tensor operations.

$\bullet$ The covariant form of the surface metric tensor $a_{\alpha\beta}$
is given by:
\begin{equation}
\left[a_{\alpha\beta}\right]=\left[\begin{array}{cc}
a_{11} & a_{12}\\
a_{21} & a_{22}
\end{array}\right]=\left[\begin{array}{cc}
E & F\\
F & G
\end{array}\right]
\end{equation}
where $E,F,G$ are the coefficients of the first fundamental form
(refer to $\S$ \ref{subFirstFundamentalForm}).

$\bullet$ The contravariant form of the surface metric tensor $a^{\alpha\beta}$
is the inverse of its covariant form and hence it is given by:

\begin{equation}
\left[a^{\alpha\beta}\right]=\left[\begin{array}{cc}
a^{11} & a^{12}\\
a^{21} & a^{22}
\end{array}\right]=\frac{1}{a_{11}a_{22}-a_{12}a_{21}}\left[\begin{array}{cc}
a_{22} & -a_{12}\\
-a_{21} & a_{11}
\end{array}\right]=\frac{1}{EG-F^{2}}\left[\begin{array}{cc}
G & -F\\
-F & E
\end{array}\right]
\end{equation}
where the symbols are as defined previously.

$\bullet$ The mixed form of the surface metric tensor $a_{\beta}^{\alpha}$
is the identity tensor, that is:
\begin{equation}
\left[a_{\beta}^{\alpha}\right]=\left[\delta_{\beta}^{\alpha}\right]=\left[\begin{array}{cc}
1 & 0\\
0 & 1
\end{array}\right]
\end{equation}

$\bullet$ Similar to the space metric tensor (refer to \cite{SochiTC1, SochiTC2}),
the surface metric tensor transforms between the barred and unbarred
surface coordinate systems as:
\begin{equation}
\bar{a}_{\alpha\beta}=a_{\gamma\delta}\frac{\partial u^{\gamma}}{\partial\bar{u}^{\alpha}}\frac{\partial u^{\delta}}{\partial\bar{u}^{\beta}}\label{eqSurMetTrans}
\end{equation}
where the indexed $\bar{a}$ and $a$ are the surface covariant metric
tensors in the barred and unbarred systems respectively. The contravariant
and mixed forms of the surface metric tensor also follow similar rules
to their counterparts of the space metric tensor as outlined in \cite{SochiTC1, SochiTC2}.

$\bullet$ Similar to the determinants of the space metric (refer
to \cite{SochiTC2}), the determinants of the surface metric in the
barred and unbarred coordinate systems are linked through the Jacobian
of transformation by the following relation:
\begin{equation}
\bar{a}=J^{2}a
\end{equation}
where $\bar{a}$ and $a$ are respectively the determinants of the
covariant form of the surface metric tensor in the barred and unbarred
systems respectively and $J$ $\left(=\left|\frac{\partial u}{\partial\bar{u}}\right|\right)$
is the Jacobian of the transformation between the two surface systems.
This relation can be obtained directly by taking the determinant of
the two sides of Eq. \ref{eqSurMetTrans}.

$\bullet$ The Christoffel symbols of the first kind $[\alpha\beta,\gamma]$
are linked to the basis vectors and their partial derivatives by the
following relation:
\begin{equation}
\left[\alpha\beta,\gamma\right]=\frac{\partial\mathbf{E}_{\alpha}}{\partial u^{\beta}}\cdot\mathbf{E}_{\gamma}
\end{equation}

$\bullet$ The relation between the partial derivative of the surface
metric tensor and the Christoffel symbols of the first kind is given
by:

\begin{equation}
\frac{\partial a_{\alpha\beta}}{\partial u^{\gamma}}=\frac{\partial\left(\mathbf{E}_{\alpha}\cdot\mathbf{E}_{\beta}\right)}{\partial u^{\gamma}}=\frac{\partial\mathbf{E}_{\alpha}}{\partial u^{\gamma}}\cdot\mathbf{E}_{\beta}+\mathbf{E}_{\alpha}\cdot\frac{\partial\mathbf{E}_{\beta}}{\partial u^{\gamma}}=\left[\alpha\gamma,\beta\right]+\left[\beta\gamma,\alpha\right]
\end{equation}

$\bullet$ Similar relations between the surface metric tensor and
the Christoffel symbols of the second kind can be obtained from the
equations in the previous point by using the index-raising operator
for the surface.

$\bullet$ \label{MetricScaling}Scaling a surface up or down by a
constant factor $c>0$, which is equivalent to scaling all the distances
on the surface by that factor, can be done by multiplying the surface
metric tensor by $c^{2}$.

$\bullet$ For a Monge patch of the form $\mathbf{r}(u,v)=\left(u,v,f(u,v)\right)$,
the surface covariant metric tensor $\mathbf{a}$ is given by:
\begin{equation}
\mathbf{a}\equiv\mathbf{I}_{S}=\left[\begin{array}{cc}
1+f_{u}^{2} & f_{u}f_{v}\\
f_{u}f_{v} & 1+f_{v}^{2}
\end{array}\right]\label{eqMongeMetric}
\end{equation}
where the subscripts $u$ and $v$ stand for partial derivatives with
respect to these surface coordinates, and $\mathbf{I}_{S}$ is the
tensor of the first fundamental form of the surface (see $\S$ \ref{subFirstFundamentalForm}).

$\bullet$ If $C:I\rightarrow S$ is a regular curve on a surface
$S$ defined on the interval $I\subseteq\mathbb{R}$, and $\mathbf{v}_{1}$
and $\mathbf{v}_{2}$ are parallel vector fields over $C$, then the
dot product associated with the metric tensor $\mathbf{v}_{1}\cdot\mathbf{v}_{2}$,
the norm of the vector fields $\left|\mathbf{v}_{1}\right|$ and $\left|\mathbf{v}_{2}\right|$,
and the angle between $\mathbf{v}_{1}$ and $\mathbf{v}_{2}$ are
constants.

$\bullet$ In the following subsections, we investigate arc length,
area and angle between two vectors on a surface. All these entities
depend in their definition and quantification on the metric tensor.
We will see that their geometric and tensor formulation is identical
to that given in \cite{SochiTC1,SochiTC2} for a general $n$D space
with the use of the surface metric tensor and the surface representation
of the involved vectors.

\subsubsection{Arc Length}

$\bullet$ Following the example of the length of an element of arc
of a curve embedded in a general $n$D space, the length of an element
of arc of a curve on a 2D surface is given in its general form by:
\begin{equation}
(ds)^{2}=d\mathbf{r}\cdot d\mathbf{r}=\frac{\partial\mathbf{r}}{\partial u^{\alpha}}\cdot\frac{\partial\mathbf{r}}{\partial u^{\beta}}\,du^{\alpha}du^{\beta}=\mathbf{E}_{\alpha}\cdot\mathbf{E}_{\beta}\,du^{\alpha}du^{\beta}=a_{\alpha\beta}du^{\alpha}du^{\beta}\label{eqArcLength}
\end{equation}
where $a_{\alpha\beta}$ is the covariant type of the surface metric
tensor, $\mathbf{r}$ is the spatial representation of the curve and
$\alpha,\beta=1,2$.

$\bullet$ From the above formula we have the following identity which
is valid at each point of a surface curve:
\begin{equation}
a_{\alpha\beta}\frac{du^{\alpha}}{ds}\frac{du^{\beta}}{ds}=1
\end{equation}

$\bullet$ Based on the above formula (Eq. \ref{eqArcLength}), the
length of a segment of a $t$-parameterized curve between a starting
point corresponding to $t=t_{1}$ and a terminal point corresponding
to $t=t_{2}$ is given by:
\begin{equation}
\begin{aligned}L & =\int_{I}ds\\
 & =\int_{t_{1}}^{t_{2}}\frac{ds}{dt}dt\\
 & =\int_{t_{1}}^{t_{2}}\sqrt{a_{\alpha\beta}\frac{du^{\alpha}}{dt}\frac{du^{\beta}}{dt}}\,dt\\
 & =\int_{t_{1}}^{t_{2}}\sqrt{a_{11}\left(\frac{\partial u^{1}}{\partial t}\right)^{2}+2a_{12}\frac{\partial u^{1}}{\partial t}\frac{\partial u^{2}}{\partial t}+a_{22}\left(\frac{\partial u^{2}}{\partial t}\right)^{2}}\,dt\\
 & =\int_{t_{1}}^{t_{2}}\sqrt{E\left(\frac{\partial u^{1}}{\partial t}\right)^{2}+2F\frac{\partial u^{1}}{\partial t}\frac{\partial u^{2}}{\partial t}+G\left(\frac{\partial u^{2}}{\partial t}\right)^{2}}\,dt
\end{aligned}
\label{eqL}
\end{equation}
where $I\subset\mathbb{R}$ is an interval on the real line and $E,F,G$
are the coefficients of the first fundamental form.

$\bullet$ For a Monge patch of the form $\mathbf{r}(u,v)=\left(u,v,f(u,v)\right)$,
the length of an element of arc of a curve is given by:
\begin{equation}
ds=\sqrt{\left(1+f_{u}^{2}\right)dudu+2f_{u}f_{v}dudv+\left(1+f_{v}^{2}\right)dvdv}
\end{equation}
where the subscripts $u$ and $v$ stand for partial derivatives with
respect to these surface coordinates.

$\bullet$ The length of a space curve is an intrinsic property since
it depends on the metric tensor only.

$\bullet$ The length of a space curve is invariant with respect to
the type of parameterization.

\subsubsection{Surface Area}

$\bullet$ The area of an infinitesimal element of a surface in the
neighborhood of a point $P$ on the surface is given by:\footnote{Here, we assume $du^{1}du^{2}$ is positive.}
\begin{equation}
\begin{aligned}d\sigma & =\left|d\mathbf{r}_{1}\times d\mathbf{r}_{2}\right|\\
 & =\left|\mathbf{E}_{1}\times\mathbf{E}_{2}\right|du^{1}du^{2}\\
 & =\sqrt{\left|\mathbf{E}_{1}\right|^{2}\left|\mathbf{E}_{2}\right|^{2}-\left(\mathbf{E}_{1}\cdot\mathbf{E}_{2}\right)^{2}}\,du^{1}du^{2}\\
 & =\sqrt{a_{11}a_{22}-\left(a_{12}\right)^{2}}\,du^{1}du^{2}\\
 & =\sqrt{EG-F^{2}}\,du^{1}du^{2}\\
 & =\sqrt{a}\,du^{1}du^{2}
\end{aligned}
\label{eqdsegma}
\end{equation}
where $\mathbf{E}_{1}$ and $\mathbf{E}_{2}$ are the surface covariant
basis vectors, $E,F,G$ are the coefficients of the first fundamental
form, $a$ is the determinant of the surface covariant metric tensor
and the indexed $a$ are its elements.\footnote{$a=a_{11}a_{22}-a_{12}a_{21}=a_{11}a_{22}-\left(a_{12}\right)^{2}$.}
All the quantities in these expressions belong to the point $P$.

$\bullet$ The area of a surface patch $S:\Omega\rightarrow\mathbb{R}^{3}$,
where $\Omega$ is a proper subset of the $\mathbb{R}^{2}$ plane,
is given by:\footnote{$S$ should be injective, sufficiently differentiable, and regular
on the interior of $\Omega$.}
\begin{equation}
\sigma=\int_{\Omega}d\sigma=\iint_{\Omega}\sqrt{a_{11}a_{22}-\left(a_{12}\right)^{2}}\,du^{1}du^{2}=\iint_{\Omega}\sqrt{EG-F^{2}}\,du^{1}du^{2}=\iint_{\Omega}\sqrt{a}\,du^{1}du^{2}
\end{equation}

$\bullet$ The formulae for the area are reminder of the volume formulae
(see \cite{SochiTC2}), so the area can be regraded as a \textit{volume}
in a 2D space.

$\bullet$ For a Monge patch of the form $\mathbf{r}(u,v)=\left(u,v,f(u,v)\right)$,
the surface area is given by:
\begin{equation}
\sigma=\iint_{\Omega}\sqrt{1+f_{u}^{2}+f_{v}^{2}}\,dudv
\end{equation}
where the subscripts $u$ and $v$ stand for partial derivatives with
respect to these surface coordinates.

\subsubsection{Angle Between Two Surface Curves}

$\bullet$ The angle between two sufficiently smooth surface curves
intersecting at a given point on the surface is defined as the angle
between their tangent vectors at that point. As there are two opposite
directions for each curve, corresponding to the two senses of traversing
the curve, there are two main angles $\theta_{1}$ and $\theta_{2}$
such that $\theta_{1}+\theta_{2}=\pi$. The principal angle between
the two curves is usually taken as the smaller of the two angles and
hence the directions are determined accordingly.\footnote{In fact there are still two possibilities for the directions but this
has no significance as far as the angle between the two curves is
concerned.}

$\bullet$ The angle between two surface curves passing through a
given point $P$ on the surface with tangent vectors $\mathbf{A}$
and $\mathbf{B}$ at $P$ is given by:
\begin{equation}
\cos\theta=\frac{\mathbf{A}\cdot\mathbf{B}}{\left|\mathbf{A}\right|\left|\mathbf{B}\right|}=\frac{a_{\alpha\beta}A^{\alpha}B^{\beta}}{\sqrt{a_{\gamma\delta}A^{\gamma}A^{\delta}}\sqrt{a_{\epsilon\zeta}B^{\epsilon}B^{\zeta}}}
\end{equation}
where the indexed $a$ are the elements of the surface covariant metric
tensor and the Greek indices run over $1,2$.

$\bullet$ If $\mathbf{A}$ and $\mathbf{B}$ are two unit surface
vectors with surface representations $A^{\delta}$ and $B^{\delta}$
($\delta=1,2$) and space representations $A^{k}$ and $B^{k}$ ($k=1,2,3$)
then the angle $\theta$ between the two vectors is given by (see
Eqs. \ref{eqSurMet2} and \ref{eqAiAalpha}):
\begin{equation}
\cos\theta=a_{\alpha\beta}A^{\alpha}B^{\beta}=g_{ij}x_{\alpha}^{i}x_{\beta}^{j}A^{\alpha}B^{\beta}=g_{ij}A^{i}B^{j}\,\,\,\,\,\,\,\,\,\,\,\,\,\,\text{(\ensuremath{\alpha,\beta=1,2} and \ensuremath{i,j=1,2,3})}
\end{equation}
and hence the surface and space representations of the two vectors
define the same angle.

$\bullet$ The vectors $\mathbf{A}$ and $\mathbf{B}$ in the previous
points are orthogonal \textit{iff} $a_{\alpha\beta}A^{\alpha}B^{\beta}=0$
at $P$.

$\bullet$ The coordinate curves at a given point $P$ on a surface
are orthogonal \textit{iff} $a_{12}\equiv F=0$ at $P$.

$\bullet$ The corresponding angles of two isometric surfaces, like
the corresponding lengths, are equal. However, the reverse is not
true in general, that is the equality of angles on two surfaces related
by a given mapping, as in the conformal mapping, does not lead to
the equality of the corresponding lengths on the two mapped surfaces.

$\bullet$ The sine of the angle $\theta$ ($\le\pi$) between two
surface unit vectors, $\mathbf{A}$ and $\mathbf{B}$, is given by:
\begin{equation}
\sin\theta=\underline{\epsilon}_{\alpha\beta}A^{\alpha}B^{\beta}
\end{equation}
which is numerically equal to the area of the parallelogram with sides
$\mathbf{A}$ and $\mathbf{B}$. Consequently, the sufficient and
necessary condition for $\mathbf{A}$ and $\mathbf{B}$ to be orthogonal
is that:
\begin{equation}
\left|\underline{\epsilon}_{\alpha\beta}A^{\alpha}B^{\beta}\right|=1
\end{equation}

\subsection{Surface Curvature Tensor\label{subSurfaceCurvature}}

$\bullet$ The surface curvature tensor\footnote{The surface curvature tensor $b_{\alpha\beta}$ may be called the
second groundform.}, $b_{\alpha\beta}$, is an absolute, rank-2, $2\times2$ symmetric
tensor.\footnote{The coefficients of the curvature tensor are real numbers.}

$\bullet$ The elements of the surface covariant curvature tensor
are given by:
\begin{equation}
b_{\alpha\beta}=-\frac{\partial\mathbf{r}}{\partial u^{\alpha}}\cdot\frac{\partial\mathbf{n}}{\partial u^{\beta}}=-\mathbf{E}_{\alpha}\cdot\frac{\partial\mathbf{n}}{\partial u^{\beta}}=-\frac{1}{2}\left(\frac{\partial\mathbf{r}}{\partial u^{\alpha}}\cdot\frac{\partial\mathbf{n}}{\partial u^{\beta}}+\frac{\partial\mathbf{r}}{\partial u^{\beta}}\cdot\frac{\partial\mathbf{n}}{\partial u^{\alpha}}\right)
\end{equation}
and also by:\footnote{The equality $\frac{\partial\mathbf{E}_{\alpha}}{\partial u^{\beta}}\cdot\mathbf{n}=-\mathbf{E}_{\alpha}\cdot\frac{\partial\mathbf{n}}{\partial u^{\beta}}$
is based on the equality $\frac{\partial\left(\mathbf{E}_{\alpha}\cdot\mathbf{n}\right)}{\partial u^{\beta}}=\frac{\partial\left(0\right)}{\partial u^{\beta}}=0$
and the product rule for differentiation.}
\begin{equation}
b_{\alpha\beta}=\frac{\partial^{2}\mathbf{r}}{\partial u^{\alpha}\partial u^{\beta}}\cdot\mathbf{n}=\frac{\partial\mathbf{E}_{\alpha}}{\partial u^{\beta}}\cdot\mathbf{n}=\frac{\partial\mathbf{E}_{\alpha}}{\partial u^{\beta}}\cdot\left(\frac{\mathbf{E}_{1}\times\mathbf{E}_{2}}{\sqrt{a}}\right)=\frac{\frac{\partial\mathbf{E}_{\alpha}}{\partial u^{\beta}}\cdot\left(\mathbf{E}_{1}\times\mathbf{E}_{2}\right)}{\sqrt{a}}\label{eqbab}
\end{equation}

$\bullet$ Considering Eq. \ref{eqbab}, the symmetry of the surface
curvature tensor (i.e. $b_{12}=b_{21}$) follows from the fact that:
\begin{equation}
\frac{\partial\mathbf{E}_{\alpha}}{\partial u^{\beta}}=\frac{\partial\mathbf{r}}{\partial u^{\beta}\partial u^{\alpha}}=\frac{\partial\mathbf{r}}{\partial u^{\alpha}\partial u^{\beta}}=\frac{\partial\mathbf{E}_{\beta}}{\partial u^{\alpha}}
\end{equation}

$\bullet$ In full tensor notation, the surface covariant curvature
tensor is given by:

\begin{equation}
b_{\alpha\beta}=\frac{1}{2}\underline{\epsilon}^{\gamma\delta}\underline{\epsilon}_{ijk}x_{\alpha;\beta}^{i}x_{\gamma}^{j}x_{\delta}^{k}=\frac{1}{\sqrt{a}}\underline{\epsilon}_{ijk}x_{\alpha;\beta}^{i}x_{1}^{j}x_{2}^{k}
\end{equation}
where $a$ is the determinant of the surface covariant metric tensor.
This formula will simplify to the following when the space coordinates
are rectangular Cartesian:
\begin{equation}
b_{\alpha\beta}=\frac{1}{\sqrt{a}}\epsilon_{ijk}\frac{\partial^{2}x^{i}}{\partial u^{\alpha}\partial u^{\beta}}x_{1}^{j}x_{2}^{k}
\end{equation}

$\bullet$ The surface curvature tensor obeys the same transformation
rules as the surface metric tensor. Hence, we have:
\begin{equation}
\bar{b}=J^{2}b\label{eqbJ2b}
\end{equation}
where $\bar{b}$ and $b$ are respectively the determinants of the
surface covariant curvature tensor in the barred and unbarred systems
respectively and $J$ $\left(=\left|\frac{\partial u}{\partial\bar{u}}\right|\right)$
is the Jacobian of the transformation between the two surface curvilinear
coordinate systems.

$\bullet$ The surface covariant curvature tensor $b_{\alpha\beta}$
in matrix form is given by:
\begin{equation}
\mathbf{b}\equiv\left[b_{\alpha\beta}\right]=\left[\begin{array}{cc}
b_{11} & b_{12}\\
b_{21} & b_{22}
\end{array}\right]=\left[\begin{array}{cc}
e & f\\
f & g
\end{array}\right]
\end{equation}
where $e,f,g$ are the coefficients of the second fundamental form
of the surface (refer to $\S$ \ref{subSecondFundamentalForm}).

$\bullet$ The mixed form of the surface curvature tensor $b_{\,\,\,\beta}^{\alpha}$
is given by:\footnote{The mixed form $b_{\alpha}^{\,\,\,\beta}=b_{\alpha\gamma}a^{\gamma\beta}$
is the transpose of the given form.}
\begin{equation}
\left[b_{\,\,\,\beta}^{\alpha}\right]=\left[a^{\alpha\gamma}b_{\gamma\beta}\right]=\frac{1}{a}\left[\begin{array}{cc}
G & -F\\
-F & E
\end{array}\right]\left[\begin{array}{cc}
e & f\\
f & g
\end{array}\right]=\frac{1}{a}\left[\begin{array}{cc}
eG-fF & fG-gF\\
fE-eF & gE-fF
\end{array}\right]
\end{equation}
where $E,F,G,e,f,g$ are the coefficients of the first and second
fundamental forms and $a=EG-F^{2}$ is the determinant of the surface
covariant metric tensor. As seen, the coefficients of $b_{\beta}^{\alpha}$
depend on the coefficients of both the first and second fundamental
forms.

$\bullet$The contravariant form of the surface curvature tensor is
given by:
\begin{equation}
\begin{aligned}\left[b^{\alpha\beta}\right] & =\left[a^{\alpha\gamma}b_{\gamma}^{\,\,\,\beta}\right]\\
 & =\left[b_{\,\,\,\gamma}^{\alpha}a^{\gamma\beta}\right]\\
 & =\frac{1}{a^{2}}\left[\begin{array}{cc}
eG^{2}-2fFG+gF^{2} & fEG-eFG-gEF+fF^{2}\\
fEG-eFG-gEF+fF^{2} & gE^{2}-2fEF+eF^{2}
\end{array}\right]
\end{aligned}
\end{equation}
where the symbols are as explained in the previous point. Like the
covariant form, the contravariant form is a symmetric tensor.

$\bullet$ As we will see, the trace of the surface mixed curvature
tensor $b_{\beta}^{\alpha}$ is twice the mean curvature $H$ (see
$\S$ \ref{subMeanCurvature}), while its determinant is the Gaussian
curvature $K$ (see $\S$ \ref{subGaussianCurvature}), that is:\footnote{Since the trace and the determinant of a tensor are its main two invariants
under permissible transformations, then $H$ and $K$ are invariant,
as will be established in the forthcoming notes. }
\begin{equation}
H=\frac{\mathrm{tr}(b_{\beta}^{\alpha})}{2}\,\,\,\,\,\,\,\,\,\,\,\,\,\,\,\,\,\,\,\,\,\,\,\,K=\mathrm{det}(b_{\beta}^{\alpha})
\end{equation}

$\bullet$ The surface covariant curvature tensor of a Monge patch
of the form $\mathbf{r}(u,v)=\left(u,v,f(u,v)\right)$ is given by:
\begin{equation}
\mathbf{b}\equiv\left[b_{\alpha\beta}\right]=\frac{1}{\sqrt{1+f_{u}^{2}+f_{v}^{2}}}\left[\begin{array}{cc}
f_{uu} & f_{uv}\\
f_{vu} & f_{vv}
\end{array}\right]\label{eqMongeb}
\end{equation}
where the subscripts $u$ and $v$ stand for partial derivatives with
respect to these surface coordinates. Since $f_{uv}=f_{vu}$, the
tensor is symmetric as it should be.

$\bullet$ The Riemann-Christoffel curvature tensor of the second
kind and the curvature tensor of a surface are linked by the following
relation:
\begin{equation}
R_{\,\,\alpha\beta\gamma}^{\delta}=b_{\alpha\gamma}b_{\,\beta}^{\delta}-b_{\alpha\beta}b_{\,\gamma}^{\delta}\label{eqRb}
\end{equation}
This relation may also be given in terms of the Riemann-Christoffel
curvature tensor of the first kind using the index-lowering operator
for the surface:

\begin{equation}
a_{\alpha\omega}R_{\,\,\beta\gamma\delta}^{\omega}=R_{\alpha\beta\gamma\delta}=b_{\alpha\gamma}b_{\beta\delta}-b_{\alpha\delta}b_{\beta\gamma}\label{eqRb2}
\end{equation}

$\bullet$ From Eqs. \ref{eqRieChrTensor} and \ref{eqRb}, the surface
curvature tensor and the surface Christoffel symbols of the second
kind and their derivatives are related by:
\begin{equation}
b_{\alpha\gamma}b_{\beta}^{\delta}-b_{\alpha\beta}b_{\gamma}^{\delta}=\frac{\partial\Gamma_{\alpha\gamma}^{\delta}}{\partial u^{\beta}}-\frac{\partial\Gamma_{\alpha\beta}^{\delta}}{\partial u^{\gamma}}+\Gamma_{\alpha\gamma}^{\omega}\Gamma_{\omega\beta}^{\delta}-\Gamma_{\alpha\beta}^{\omega}\Gamma_{\omega\gamma}^{\delta}\label{eqbGam}
\end{equation}

$\bullet$ Considering the fact that the 2D Riemann-Christoffel curvature
tensor has only one degree of freedom and hence it possesses a single
independent non-vanishing component which is represented by $R_{1212}$,
we see that Eq. \ref{eqRb2} has only one independent component which
is given by:
\begin{equation}
R_{1212}=b_{11}b_{22}-b_{12}b_{21}=b
\end{equation}
where $b$ is the determinant of the surface covariant curvature tensor.
This equation shows that each one of the following provisions: $R_{\alpha\beta\gamma\delta}=0$
and $b_{\alpha\beta}=0$ if satisfied identically is a sufficient
and necessary condition for having a flat 2D space, i.e. a plane surface.
Hence, for a plane surface, all the coefficients of the Riemann-Christoffel
curvature tensor and the coefficients of the surface curvature tensor
vanish identically throughout the surface.

$\bullet$ From Eqs. \ref{eqRieChrTensor} and \ref{eqRb2} we see
that the Riemann-Christoffel curvature tensor can be expressed in
terms of the coefficients of the surface curvature tensor as well
as in terms of the coefficients of the surface metric tensor where
the two sets of coefficients are linked through Eq. \ref{eqbGam}.\footnote{This does not mean that Riemann curvature is an extrinsic property
but it means that some intrinsic properties can also be defined in
terms of extrinsic parameters.}

$\bullet$ The sign of the surface curvature tensor $\mathbf{b}$
(i.e. the sign of its coefficients) is dependent on the choice of
the direction of the unit vector $\mathbf{n}$ normal to the surface.

\subsection{First Fundamental Form\label{subFirstFundamentalForm}}

$\bullet$ As indicated previously, the first fundamental form\footnote{In older books, this is also called the first fundamental quadratic
form.}, which is based on the metric, encompasses all the intrinsic information
about the surface that a 2D inhabitant of the surface can obtain from
measurements performed on the surface without appealing to an external
dimension.

$\bullet$ The first fundamental form of the length of an element
of arc of a curve on a surface is a quadratic expression given by:
\begin{equation}
\begin{aligned}I_{S} & =(ds)^{2}\\
 & =d\mathbf{r}\cdot d\mathbf{r}\\
 & =a_{\alpha\beta}du^{\alpha}du^{\beta}\\
 & =E(du^{1})^{2}+2F\,du^{1}du^{2}+G(du^{2})^{2}
\end{aligned}
\label{eqFundamentalForm1}
\end{equation}
where $E,F$ and $G$, which in general are continuous variable functions
of the surface coordinates ($u,v$), are given by:
\begin{eqnarray}
E & = & a_{11}=\mathbf{E}_{1}\cdot\mathbf{E}_{1}=\frac{\partial\mathbf{r}}{\partial u^{1}}\cdot\frac{\partial\mathbf{r}}{\partial u^{1}}=g_{ij}\frac{\partial x^{i}}{\partial u^{1}}\frac{\partial x^{j}}{\partial u^{1}}\nonumber \\
F & = & a_{12}=\mathbf{E}_{1}\cdot\mathbf{E}_{2}=\frac{\partial\mathbf{r}}{\partial u^{1}}\cdot\frac{\partial\mathbf{r}}{\partial u^{2}}=g_{ij}\frac{\partial x^{i}}{\partial u^{1}}\frac{\partial x^{j}}{\partial u^{2}}=\mathbf{E}_{2}\cdot\mathbf{E}_{1}=a_{21}\\
G & = & a_{22}=\mathbf{E}_{2}\cdot\mathbf{E}_{2}=\frac{\partial\mathbf{r}}{\partial u^{2}}\cdot\frac{\partial\mathbf{r}}{\partial u^{2}}=g_{ij}\frac{\partial x^{i}}{\partial u^{2}}\frac{\partial x^{j}}{\partial u^{2}}\nonumber
\end{eqnarray}
where the indexed $a$ are the elements of the surface covariant metric
tensor, the indexed $x$ are the curvilinear coordinates of the enveloping
space and $g_{ij}$ is its covariant metric tensor.

$\bullet$ For a flat space with a Cartesian coordinate system $x^{i}$,
the space metric is $g_{ij}=\delta_{ij}$ and hence the above equations
become:
\begin{equation}
\begin{aligned}E & =\frac{\partial x^{i}}{\partial u^{1}}\frac{\partial x^{i}}{\partial u^{1}}\\
F & =\frac{\partial x^{i}}{\partial u^{1}}\frac{\partial x^{i}}{\partial u^{2}}\\
G & =\frac{\partial x^{i}}{\partial u^{2}}\frac{\partial x^{i}}{\partial u^{2}}
\end{aligned}
\end{equation}

$\bullet$ The first fundamental form can be cast in the following
matrix form:
\begin{equation}
\begin{aligned}I_{S} & =\left[\begin{array}{cc}
du^{1} & du^{2}\end{array}\right]\left[\begin{array}{c}
\mathbf{E}_{1}\\
\mathbf{E}_{2}
\end{array}\right]\cdot\left[\begin{array}{cc}
\mathbf{E}_{1} & \mathbf{E}_{2}\end{array}\right]\left[\begin{array}{c}
du^{1}\\
du^{2}
\end{array}\right]\\
 & =\left[\begin{array}{cc}
du^{1} & du^{2}\end{array}\right]\left[\begin{array}{cc}
\mathbf{E}_{1}\cdot\mathbf{E}_{1} & \mathbf{E}_{1}\cdot\mathbf{E}_{2}\\
\mathbf{E}_{2}\cdot\mathbf{E}_{1} & \mathbf{E}_{2}\cdot\mathbf{E}_{2}
\end{array}\right]\left[\begin{array}{c}
du^{1}\\
du^{2}
\end{array}\right]\\
 & =\left[\begin{array}{cc}
du^{1} & du^{2}\end{array}\right]\left[\begin{array}{cc}
E & F\\
F & G
\end{array}\right]\left[\begin{array}{c}
du^{1}\\
du^{2}
\end{array}\right]\\
 & =\left[\begin{array}{cc}
du^{1} & du^{2}\end{array}\right]\left[\begin{array}{cc}
a_{11} & a_{12}\\
a_{21} & a_{22}
\end{array}\right]\left[\begin{array}{c}
du^{1}\\
du^{2}
\end{array}\right]\\
 & =\mathbf{v\,}\mathbf{I}_{S}\,\mathbf{v}^{T}
\end{aligned}
\end{equation}
where $\mathbf{v}$ is a direction vector, $\mathbf{v}^{T}$ is its
transpose, and $\mathbf{I}_{S}$ is the first fundamental form tensor
which is the surface covariant metric tensor. Hence, the matrix associated
with the first fundamental form is the covariant metric tensor of
the surface.

$\bullet$ \label{PlaneCylindeHaveSamerFirstForm}The first fundamental
form is not a unique characteristic of the surface and hence two geometrically
different surfaces as seen from the enveloping space, such as plane
and cylinder, can have the same first fundamental form.\footnote{Such surfaces are different extrinsically as seen from the embedding
space although they are identical intrinsically as viewed internally
from the surface by a 2D inhabitant.}

$\bullet$ The first fundamental form is positive definite at regular
points of 2D surfaces, hence its coefficients are subject to the conditions
$E>0$ and $\mathrm{det}(\mathbf{I}_{S})=EG-F^{2}>0$.\footnote{As indicated earlier, the conditions $E>0$ and $EG-F^{2}>0$ imply
$G>0$.} However, this condition may be amended to allow for metrics with
imaginary coordinates as it is the case in the coordinate systems
of relativistic mechanics.

$\bullet$ As indicated previously, the first fundamental form encompasses
the intrinsic properties of the surface geometry. Hence, as seen in
$\S$ \ref{subSurfaceMetric}, the first fundamental form is used
to define and quantify things like arc length, area and angle between
curves on a surface based on its qualification as a metric. For example,
Eq. \ref{eqL} shows that the length of a curve segment on a surface
is obtained by integrating the square root of the first fundamental
form of the surface along the segment.

$\bullet$ If a surface $S_{1}$ can be mapped isometrically (see
$\S$ \ref{subSurfaces}) onto another surface $S_{2}$ then the two
surfaces have identical first fundamental form coefficients, that
is $E_{1}=E_{2}$, $F_{1}=F_{2}$ and $G_{1}=G_{2}$ where the subscripts
are labels for the two surfaces.

$\bullet$ From a previous point (see Eq. \ref{eqMongeMetric}), for
a Monge patch of the form $\mathbf{r}(u,v)=\left(u,v,f(u,v)\right)$,
the first fundamental form is given by:
\begin{equation}
I_{S}=\left(1+f_{u}^{2}\right)dudu+2f_{u}f_{v}dudv+\left(1+f_{v}^{2}\right)dvdv
\end{equation}
where the subscripts $u$ and $v$ stand for partial derivatives with
respect to these surface coordinates.

\subsection{Second Fundamental Form\label{subSecondFundamentalForm}}

$\bullet$ The mathematical entity that characterizes the extrinsic
geometry of a surface is the normal vector to the surface. This entity
can only be observed externally from outside the surface by an observer
in a reference frame in the space that envelopes the surface. Hence,
the normal vector and all its subsidiaries are strange to a 2D inhabitant
to the surface who can only access the intrinsic attributes of the
surface as represented by and contained in the first fundamental form.

$\bullet$ As a consequence of the last point, the variation of the
normal vector as it moves around the surface can be used as an indicator
to characterize the surface shape from an external point of view and
that is how this indicator is exploited in the second fundamental
form to represent the extrinsic geometry of the surface as will be
seen from the forthcoming formulations.

$\bullet$ The following quadratic expression is called the second
fundamental form\footnote{In older books, this is also called the second fundamental quadratic
form.}, of the surface:
\begin{equation}
\begin{aligned}II_{S} & =-d\mathbf{r}\cdot d\mathbf{n}\\
 & =-\left(\frac{\partial\mathbf{r}}{\partial u^{\alpha}}du^{\alpha}\right)\cdot\left(\frac{\partial\mathbf{n}}{\partial u^{\beta}}du^{\beta}\right)\\
 & =-\left(\frac{\partial\mathbf{r}}{\partial u^{1}}du^{1}+\frac{\partial\mathbf{r}}{\partial u^{2}}du^{2}\right)\cdot\left(\frac{\partial\mathbf{n}}{\partial u^{1}}du^{1}+\frac{\partial\mathbf{n}}{\partial u^{2}}du^{2}\right)\\
 & =e(du^{1})^{2}+2f\,du^{1}du^{2}+g(du^{2})^{2}
\end{aligned}
\label{eqFundamentalForm21}
\end{equation}
where\footnote{Some authors use $L,M,N$ instead of $e,f,g$. However, the use of
$e,f,g$ is advantageous since they correspond to the coefficients
of the first fundamental form $E,F,G$ nicely making the formulae
involving the first and second fundamental forms more symmetric and
memorable. On the other hand, the use of $L,M,N$ is also advantageous
when reading formulae containing the coefficients of both fundamental
forms; moreover, it is less susceptible to errors when writing or
typing these formulae. Another point is that the coefficient $g$
should not be confused with the determinant of the space covariant
metric tensor which is used in the previous notes \cite{SochiTC2}
but we do not use it in the present text.}
\begin{eqnarray}
e & = & -\frac{\partial\mathbf{r}}{\partial u^{1}}\cdot\frac{\partial\mathbf{n}}{\partial u^{1}}=-\mathbf{E}_{1}\cdot\frac{\partial\mathbf{n}}{\partial u^{1}}\nonumber \\
f & = & -\frac{1}{2}\left(\frac{\partial\mathbf{r}}{\partial u^{1}}\cdot\frac{\partial\mathbf{n}}{\partial u^{2}}+\frac{\partial\mathbf{r}}{\partial u^{2}}\cdot\frac{\partial\mathbf{n}}{\partial u^{1}}\right)=-\frac{1}{2}\left(\mathbf{E}_{1}\cdot\frac{\partial\mathbf{n}}{\partial u^{2}}+\mathbf{E}_{2}\cdot\frac{\partial\mathbf{n}}{\partial u^{1}}\right)\\
g & = & -\frac{\partial\mathbf{r}}{\partial u^{2}}\cdot\frac{\partial\mathbf{n}}{\partial u^{2}}=-\mathbf{E}_{2}\cdot\frac{\partial\mathbf{n}}{\partial u^{2}}\nonumber
\end{eqnarray}
In the above equations, $\mathbf{r}(u^{1},u^{2})$ is the spatial
representation of the surface, $\mathbf{n}(u^{1},u^{2})$ is the unit
vector normal to the surface and $\alpha,\beta=1,2$.

$\bullet$ The second fundamental form is also given by:
\begin{equation}
\begin{aligned}II_{S} & =d^{2}\mathbf{r}\cdot\mathbf{n}\\
 & =\left(\frac{\partial^{2}\mathbf{r}}{\partial u^{\alpha}\partial u^{\beta}}du^{\alpha}du^{\beta}\right)\cdot\mathbf{n}\\
 & =\left(\frac{\partial\mathbf{E}_{1}}{\partial u^{1}}(du^{1})^{2}+2\frac{\partial\mathbf{E}_{1}}{\partial u^{2}}du^{1}du^{2}+\frac{\partial\mathbf{E}_{2}}{\partial u^{2}}(du^{2})^{2}\right)\cdot\mathbf{n}\\
 & =\frac{\partial\mathbf{E}_{1}}{\partial u^{1}}\cdot\mathbf{n}\,(du^{1})^{2}+2\frac{\partial\mathbf{E}_{1}}{\partial u^{2}}\cdot\mathbf{n}\,du^{1}du^{2}+\frac{\partial\mathbf{E}_{2}}{\partial u^{2}}\cdot\mathbf{n}\,(du^{2})^{2}
\end{aligned}
\end{equation}
where $d^{2}\mathbf{r}$ is the second order differential of the position
vector $\mathbf{r}$ of an arbitrary point on the surface in the direction
($du^{1},du^{2}$).

$\bullet$ From the last point, the coefficients of the second fundamental
form can also be given by the following alternative expressions:
\begin{eqnarray}
e & = & \mathbf{n}\cdot\frac{\partial\mathbf{E}_{1}}{\partial u^{1}}=-\frac{\partial\mathbf{n}}{\partial u^{1}}\cdot\mathbf{E}_{1}\nonumber \\
f & = & \mathbf{n}\cdot\frac{\partial\mathbf{E}_{1}}{\partial u^{2}}=\mathbf{n}\cdot\frac{\partial\mathbf{E}_{2}}{\partial u^{1}}=-\frac{\partial\mathbf{n}}{\partial u^{2}}\cdot\mathbf{E}_{1}=-\frac{\partial\mathbf{n}}{\partial u^{1}}\cdot\mathbf{E}_{2}\\
g & = & \mathbf{n}\cdot\frac{\partial\mathbf{E}_{2}}{\partial u^{2}}=-\frac{\partial\mathbf{n}}{\partial u^{2}}\cdot\mathbf{E}_{2}\nonumber
\end{eqnarray}
and also by:
\begin{eqnarray}
e & = & \frac{\left(\mathbf{E}_{1}\times\mathbf{E}_{2}\right)\cdot\frac{\partial\mathbf{E}_{1}}{\partial u^{1}}}{\sqrt{a}}\nonumber \\
f & = & \frac{\left(\mathbf{E}_{1}\times\mathbf{E}_{2}\right)\cdot\frac{\partial\mathbf{E}_{1}}{\partial u^{2}}}{\sqrt{a}}\\
g & = & \frac{\left(\mathbf{E}_{1}\times\mathbf{E}_{2}\right)\cdot\frac{\partial\mathbf{E}_{2}}{\partial u^{2}}}{\sqrt{a}}\nonumber
\end{eqnarray}
where $a=a_{11}a_{22}-a_{12}a_{21}=EG-F^{2}$ is the determinant of
the surface covariant metric tensor.

$\bullet$ The second fundamental form can be cast in the following
matrix form:
\begin{eqnarray}
II_{S} & = & \left[\begin{array}{cc}
du^{1} & du^{2}\end{array}\right]\left[\begin{array}{c}
-\mathbf{E}_{1}\\
-\mathbf{E}_{2}
\end{array}\right]\cdot\left[\begin{array}{cc}
\frac{\partial\mathbf{n}}{\partial u^{1}} & \frac{\partial\mathbf{n}}{\partial u^{2}}\end{array}\right]\left[\begin{array}{c}
du^{1}\\
du^{2}
\end{array}\right]\nonumber \\
 & = & \left[\begin{array}{cc}
du^{1} & du^{2}\end{array}\right]\left[\begin{array}{cc}
-\mathbf{E}_{1}\cdot\frac{\partial\mathbf{n}}{\partial u^{1}} & -\mathbf{E}_{1}\cdot\frac{\partial\mathbf{n}}{\partial u^{2}}\\
-\mathbf{E}_{2}\cdot\frac{\partial\mathbf{n}}{\partial u^{1}} & -\mathbf{E}_{2}\cdot\frac{\partial\mathbf{n}}{\partial u^{2}}
\end{array}\right]\left[\begin{array}{c}
du^{1}\\
du^{2}
\end{array}\right]\nonumber \\
 & = & \left[\begin{array}{cc}
du^{1} & du^{2}\end{array}\right]\left[\begin{array}{cc}
\frac{\partial\mathbf{E}_{1}}{\partial u^{1}}\cdot\mathbf{n} & \frac{\partial\mathbf{\mathbf{E}_{1}}}{\partial u^{2}}\cdot\mathbf{n}\\
\frac{\partial\mathbf{E}_{2}}{\partial u^{1}}\cdot\mathbf{n} & \frac{\partial\mathbf{E}_{2}}{\partial u^{2}}\cdot\mathbf{n}
\end{array}\right]\left[\begin{array}{c}
du^{1}\\
du^{2}
\end{array}\right]\\
 & = & \left[\begin{array}{cc}
du^{1} & du^{2}\end{array}\right]\left[\begin{array}{cc}
e & f\\
f & g
\end{array}\right]\left[\begin{array}{c}
du^{1}\\
du^{2}
\end{array}\right]\nonumber \\
 & = & \mathbf{v\,}\mathbf{II}_{S}\,\mathbf{v}^{T}\nonumber
\end{eqnarray}
where $\mathbf{v}$ is a direction vector, $\mathbf{v}^{T}$ is its
transpose, and $\mathbf{II}_{S}$ is the second fundamental form tensor.

$\bullet$ Like the coefficients of the first fundamental form, the
coefficients of the second fundamental form are, in general, continuous
variable functions of the surface coordinates ($u,v$).

$\bullet$ The coefficients of the second fundamental form tensor
satisfy the following relations:
\begin{equation}
e=b_{11}\,\,\,\,\,\,\,\,\,\,\,\,\,\,\,\,\,f=b_{12}=b_{21}\,\,\,\,\,\,\,\,\,\,\,\,\,\,\,\,\,g=b_{22}
\end{equation}
Hence, the second fundamental form tensor $\mathbf{II}_{S}$ is the
same as the surface covariant curvature tensor $\mathbf{b}$, that
is:
\begin{equation}
\mathbf{II}_{S}=\left[\begin{array}{cc}
e & f\\
f & g
\end{array}\right]=\left[\begin{array}{cc}
b_{11} & b_{12}\\
b_{21} & b_{22}
\end{array}\right]=\mathbf{b}
\end{equation}

$\bullet$ From the previous points, we see that the second fundamental
form can also be given by:

\begin{equation}
II_{S}=b_{\alpha\beta}du^{\alpha}du^{\beta}
\end{equation}
where the indexed $b$ are the elements of the surface covariant curvature
tensor and $\alpha,\beta=1,2$.

$\bullet$ The second fundamental form of the surface can be expressed
in terms of the first fundamental form $I_{S}$ and the normal curvature
$\kappa_{n}$ (see $\S$ \ref{subNormalCurvature}) of the surface
at a given point and in a given direction as:
\begin{equation}
II_{S}=\kappa_{n}I_{S}=\kappa_{n}(ds)^{2}\label{eqFundamentalForm22}
\end{equation}

$\bullet$ From a previous point (see Eq. \ref{eqMongeb}), for a
Monge patch of the form $\mathbf{r}(u,v)=\left(u,v,f(u,v)\right)$,
the second fundamental form is given by:
\begin{equation}
II_{S}=\frac{f_{uu}dudu+2f_{uv}dudv+f_{vv}dvdv}{\sqrt{1+f_{u}^{2}+f_{v}^{2}}}
\end{equation}
where the subscripts $u$ and $v$ stand for partial derivatives with
respect to these surface curvilinear coordinates.

$\bullet$ The second fundamental form is invariant under permissible
coordinate transformations that maintain the sense of the normal vector
to the surface, $\mathbf{n}$. The second fundamental form changes
its sign if the sense of $\mathbf{n}$ is reversed.

$\bullet$ As indicated earlier, while the first fundamental form
encompasses the intrinsic geometry of the surface, the second fundamental
form encompasses its extrinsic geometry.

$\bullet$ As seen, the first fundamental form is associated with
the surface covariant metric tensor, while the second fundamental
form is associated with the surface covariant curvature tensor.

$\bullet$ While the first fundamental form is positive definite,
as stated previously, the second fundamental form is not necessarily
positive or definite.

$\bullet$ Unlike space curves which are completely defined by specified
curvature and torsion, $\kappa$ and $\tau$, providing arbitrary
first and second fundamental forms is not a sufficient condition for
the existence of a surface with these forms, because the first and
second fundamental forms do not provide full identification for the
surface. In a more technical terms, defining six functions $E,F,G,e,f$
and $g$ of class $C^{3}$ on a subset of $\mathbb{R}^{2}$ where
these functions satisfy the conditions for the coefficients of the
first and second fundamental forms (in particular $E,G>0$ and $EG-F^{2}>0$)
does not guarantee the existence of a surface over the given subset
with a first fundamental form $E\left(du^{1}\right)^{2}+2Fdu^{1}du^{2}+G\left(du^{2}\right)^{2}$
and a second fundamental form $e\left(du^{1}\right)^{2}+2fdu^{1}du^{2}+g\left(du^{2}\right)^{2}$.
Further compatibility conditions relating the first and second fundamental
forms are required to fully identify the surface and secure its existence.\footnote{This may be linked to the fact that the curve conditions are established
based on the existence theorem for ordinary differential equations
where these equations generally have a solution, while the surface
conditions should be established based on the existence theorem for
partial differential equations which have solutions only when they
meet additional integrability conditions. The details should be sought
in more extensive books of differential geometry.}

$\bullet$ Following the last point, the required compatibility conditions
for the existence of a surface with predefined first and second fundamental
forms are given by the Codazzi-Mainardi equations (Eq. \ref{eqCodazziMainardi})
plus the following equation:
\begin{eqnarray}
eg-f^{2} & = & F\left[\frac{\partial\Gamma_{22}^{2}}{\partial u}-\frac{\partial\Gamma_{12}^{2}}{\partial v}+\Gamma_{22}^{1}\Gamma_{11}^{2}-\Gamma_{12}^{1}\Gamma_{12}^{2}\right]+\label{eqCompatibility2}\\
 &  & E\left[\frac{\partial\Gamma_{22}^{1}}{\partial u}-\frac{\partial\Gamma_{12}^{1}}{\partial v}+\Gamma_{22}^{1}\Gamma_{11}^{1}+\Gamma_{22}^{2}\Gamma_{12}^{1}-\Gamma_{12}^{1}\Gamma_{12}^{1}-\Gamma_{12}^{2}\Gamma_{22}^{1}\right]\nonumber
\end{eqnarray}

$\bullet$ From the last two points the fundamental theorem of surfaces
in differential geometry, which is the equivalent of the fundamental
theorem of curves (see $\S$ \ref{subCurvatureTorsion}), emerges.
The theorem states that: given six sufficiently-smooth functions $E,F,G,e,f$
and $g$ on a subset of $\mathbb{R}^{2}$ satisfying the following
conditions:

(A) $E,G>0$ and $EG-F^{2}>0$, and

(B) $E,F,G,e,f$ and $g$ satisfy Eqs. \ref{eqCompatibility2} and
\ref{eqCodazziMainardi},

then there is a unique surface with $E,F,G$ as its first fundamental
form coefficients and $e,f,g$ as its second fundamental form coefficients.
Hence, if two surfaces meet all these conditions, then they are identical
within a rigid motion transformation in space.

$\bullet$ As seen, the fundamental theorem of surfaces, like the
fundamental theorem of curves, provides the existence and uniqueness
conditions for surfaces.

$\bullet$ On the other hand, according to the theorem of Bonnet,
if two surfaces of class $C^{3}$, $S_{1}:\Omega\rightarrow\mathbb{R}^{3}$
and $S_{2}:\Omega\rightarrow\mathbb{R}^{3}$, are defined over a connected
set $\Omega\subseteq\mathbb{R}^{2}$ with identical first and second
fundamental forms, then the two surfaces can be mapped on each other
by a purely rigid motion transformation.

$\bullet$ Two surfaces having identical first fundamental forms but
different second fundamental forms may be described as applicable.
An example of applicable surfaces are plane and cylinder.

\subsubsection{Dupin Indicatrix\label{subDupinIndicatrix}}

$\bullet$ Dupin indicatrix at a given point $P$ on a sufficiently
smooth surface is a function of the coefficients of the second fundamental
form at the point and hence it is a function of the surface coordinates
at $P$.

$\bullet$ Dupin indicatrix is an indicator of the departure of the
surface from the tangent plane in the close proximity of the point
of tangency. Accordingly, the second fundamental form is used in Dupin
indicatrix to measure this departure.

$\bullet$ In quantitative terms, Dupin indicatrix is the family of
conic sections given by the following quadratic equation:
\begin{equation}
ex^{2}+2fxy+gy^{2}=\pm1
\end{equation}
where $e,f,g$ are the coefficients of the second fundamental form
at $P$.

$\bullet$ As a consequence of the previous points, Dupin indicatrix
can be used to classify the surface points with respect to the local
shape of the surface as flat, elliptic, parabolic or hyperbolic; more
details about this are given in $\S$ \ref{subLocalShape}.

\subsection{Third Fundamental Form}

$\bullet$ The third fundamental form of a space surface is defined
by:
\begin{equation}
III_{S}=d\mathbf{n}\cdot d\mathbf{n}=c_{\alpha\beta}du^{\alpha}du^{\beta}
\end{equation}
where $\mathbf{n}$ is the unit vector normal to the surface at a
given point $P$, $c_{\alpha\beta}$ are the coefficients of the third
fundamental form at $P$ and $\alpha,\beta=1,2$.

$\bullet$ The coefficients of the third fundamental form are given
by:\footnote{These coefficients are real numbers.}
\begin{equation}
c_{\alpha\beta}=g_{ij}n_{\,\,;\alpha}^{i}n_{\,\,;\beta}^{j}
\end{equation}
where $g_{ij}$ is the space covariant metric tensor and the indexed
$n$ is the unit vector normal to the surface.

$\bullet$ The coefficients of the third fundamental form are also
given by:
\begin{equation}
c_{\alpha\beta}=a^{\gamma\delta}b_{\alpha\gamma}b_{\beta\delta}
\end{equation}
where $a^{\gamma\delta}$ is the surface contravariant metric tensor
and the indexed $b$ are the coefficients of the surface covariant
curvature tensor.

$\bullet$ The first, second and third fundamental forms are linked,
through the Gaussian curvature $K$ and the mean curvature $H$ (see
$\S$ \ref{subGaussianCurvature} and \ref{subMeanCurvature}), by
the following relation:
\begin{equation}
KI_{S}-2HII_{S}+III_{S}=0
\end{equation}

$\bullet$ The coefficients of the first, second and third fundamental
forms are correlated, through the mean curvature $H$ and the Gaussian
curvature $K$, by the following relation:
\begin{equation}
Ka_{\alpha\beta}-2Hb_{\alpha\beta}+c_{\alpha\beta}=0
\end{equation}
By multiplying both sides with $a^{\alpha\beta}$ and contracting
we obtain:
\begin{equation}
Ka_{\alpha}^{\alpha}-2Hb_{\alpha}^{\alpha}+c_{\alpha}^{\alpha}=0
\end{equation}
that is:\footnote{We have: $a_{\alpha}^{\alpha}=\delta_{\alpha}^{\alpha}=\delta_{1}^{1}+\delta_{2}^{2}=2$
and $H=\frac{\mathrm{tr}\left(b_{\alpha\beta}\right)}{2}$ and hence
the formula is justified.}
\begin{equation}
\mathrm{tr}\left(c_{\alpha\beta}\right)=4H^{2}-2K
\end{equation}

\subsection{Relationship between Surface Basis Vectors and their Derivatives}

$\bullet$ The focus of this subsection is the equations of Gauss
and Weingarten which, for surfaces, are the analogue of the equations
of Frenet-Serret for curves. While the Frenet-Serret formulae express
the derivatives of $\mathbf{T},\mathbf{N},\mathbf{B}$ as combinations
of these vectors using $\kappa$ and $\tau$ as coefficients, the
equations of Gauss and Weingarten express the derivatives of $\mathbf{E}_{1},\mathbf{E}_{2},\mathbf{n}$
as combinations of these vectors with coefficients based on the first
and second fundamental forms.

$\bullet$ As shown earlier (see $\S$ \ref{subSpaceCurves}), three
unit vectors can be constructed on each point at which the curvature
does not vanish of a class $C^{2}$ space curve: the tangent $\mathbf{T}$,
the normal $\mathbf{N}$ and the binormal $\mathbf{B}$. These mutually
orthogonal vectors (i.e. $\mathbf{T}\cdot\mathbf{N}=\mathbf{T}\cdot\mathbf{B}=\mathbf{N}\cdot\mathbf{B}=0$)
can serve as a set of basis vectors. Hence, the derivatives of these
vectors with respect to the distance traversed along the curve, $s$,
can be expressed as combinations of this set as demonstrated by the
Frenet-Serret formulae (refer to $\S$ \ref{subRelationCurveBasis}).

$\bullet$ Similarly, the surface vectors: $\mathbf{E}_{1}=\frac{\partial\mathbf{r}}{\partial u^{1}}$,
$\mathbf{E}_{2}=\frac{\partial\mathbf{r}}{\partial u^{2}}$ and the
unit vector normal to the surface, $\mathbf{n}$, at each regular
point on a class $C^{2}$ surface also form a basis set and hence
their partial derivatives with respect to the surface curvilinear
coordinates, $u^{1}$ and $u^{2}$, can be expressed as combinations
of this set. The equations of Gauss and Weingarten demonstrate this
fact.

$\bullet$ The equations of Gauss express the partial derivatives
of the surface vectors, $\mathbf{E}_{1}$ and $\mathbf{E}_{2}$, with
respect to the surface curvilinear coordinates as combinations of
the surface basis set, that is:
\begin{eqnarray}
\frac{\partial\mathbf{E}_{1}}{\partial u^{1}} & = & \Gamma_{11}^{1}\mathbf{E}_{1}+\Gamma_{11}^{2}\mathbf{E}_{2}+e\mathbf{n}\nonumber \\
\frac{\partial\mathbf{E}_{1}}{\partial u^{2}} & = & \Gamma_{12}^{1}\mathbf{E}_{1}+\Gamma_{12}^{2}\mathbf{E}_{2}+f\mathbf{n}=\frac{\partial\mathbf{E}_{2}}{\partial u^{1}}\\
\frac{\partial\mathbf{E}_{2}}{\partial u^{2}} & = & \Gamma_{22}^{1}\mathbf{E}_{1}+\Gamma_{22}^{2}\mathbf{E}_{2}+g\mathbf{n}\nonumber
\end{eqnarray}
where $e,f,g$ are the coefficients of the second fundamental form.
These equations can be expressed compactly, with partial use of tensor
notation, as:
\begin{equation}
\frac{\partial\mathbf{E}_{\alpha}}{\partial u^{\beta}}=\Gamma_{\alpha\beta}^{\gamma}\mathbf{E}_{\gamma}+b_{\alpha\beta}\mathbf{n}\,\,\,\,\,\,\,\,\,\,\,\,\,\,\left(\alpha,\beta=1,2\right)\label{eqbCoe}
\end{equation}
where the Christoffel symbol $\Gamma_{\alpha\beta}^{\gamma}$ is based
on the surface metric, as given by Eq. \ref{eqChristoffel2}, and
$b_{\alpha\beta}$ is the surface covariant curvature tensor. The
last equation can be expressed in full tensor notation as:
\begin{equation}
x_{\,\,\alpha,\beta}^{i}=\Gamma_{\alpha\beta}^{\gamma}x_{\gamma}^{i}+b_{\alpha\beta}n^{i}
\end{equation}

$\bullet$ Likewise, the equations of Weingarten express the partial
derivatives of the unit vector normal to the surface, $\mathbf{n}$,
with respect to the surface curvilinear coordinates as combinations
of the surface vectors, $\mathbf{E}_{1}$ and $\mathbf{E}_{2}$, that
is:
\begin{equation}
\begin{aligned}\frac{\partial\mathbf{n}}{\partial u^{1}} & =\frac{fF-eG}{a}\mathbf{E}_{1}+\frac{eF-fE}{a}\mathbf{E}_{2}\\
\frac{\partial\mathbf{n}}{\partial u^{2}} & =\frac{gF-fG}{a}\mathbf{E}_{1}+\frac{fF-gE}{a}\mathbf{E}_{2}
\end{aligned}
\label{eqnr}
\end{equation}
where $E,F,G,e,f,g$ are the coefficients of the first and second
fundamental forms and $a=EG-F^{2}$ is the determinant of the surface
covariant metric tensor. These equations can be expressed compactly,
with partial use of tensor notation, as:
\begin{equation}
\frac{\partial\mathbf{n}}{\partial u^{\alpha}}=-b_{\alpha}^{\,\,\beta}\mathbf{E}_{\beta}
\end{equation}
where $b_{\alpha}^{\,\,\beta}$ ($=b_{\alpha\gamma}a^{\gamma\beta}$)
is the mixed type of the surface curvature tensor, $b_{\alpha\gamma}$
is the surface covariant curvature tensor and $a^{\gamma\beta}$ is
the surface contravariant metric tensor. They can also be expressed
with full use of tensor notation employing curvilinear space coordinates
as:\footnote{The vector $n_{\,\,,\alpha}^{i}$ is orthogonal to $n^{i}$ and hence
it is parallel to the tangent space of the surface, so it can be expressed
as a linear combination of the surface basis vectors $x_{\beta}^{i}$,
that is $n_{\,\,,\alpha}^{i}=d_{\alpha}^{\beta}x_{\beta}^{i}$ for
a certain set of coefficients $d_{\alpha}^{\beta}=-b_{\alpha\gamma}a^{\gamma\beta}$,
as given above.}
\begin{equation}
n_{\,\,,\alpha}^{i}=-b_{\alpha\gamma}a^{\gamma\beta}x_{\beta}^{i}=-b_{\alpha}^{\,\,\beta}x_{\beta}^{i}
\end{equation}

$\bullet$ Weingarten equations can be expressed in matrix form as:
\begin{equation}
\left[\begin{array}{c}
\frac{\partial\mathbf{n}}{\partial u^{1}}\\
\frac{\partial\mathbf{n}}{\partial u^{2}}
\end{array}\right]=-\mathbf{II}_{S}\mathbf{I}_{S}^{-1}\left[\begin{array}{c}
\mathbf{E}_{1}\\
\mathbf{E}_{2}
\end{array}\right]
\end{equation}
where $\mathbf{II}_{S}$ is the surface covariant curvature tensor
and $\mathbf{I}_{S}^{-1}$ is the surface contravariant metric tensor.

$\bullet$ The partial derivatives of the unit vector normal to the
surface, $\mathbf{n}$, with respect to the surface curvilinear coordinates
($u^{1},u^{2}$) are linked to the Gaussian curvature $K$ (see $\S$
\ref{subGaussianCurvature}) and the surface basis vectors, $\mathbf{E}_{1}$
and $\mathbf{E}_{2}$, by the following relation:
\begin{equation}
\frac{\partial\mathbf{n}}{\partial u^{1}}\times\frac{\partial\mathbf{n}}{\partial u^{2}}=\frac{eg-f^{2}}{EG-F^{2}}\left(\mathbf{E}_{1}\times\mathbf{E}_{2}\right)=K\left(\mathbf{E}_{1}\times\mathbf{E}_{2}\right)\label{eqparnparn}
\end{equation}

$\bullet$ The partial derivatives of the unit vector normal to the
surface, $\mathbf{n}$, with respect to the surface curvilinear coordinates
($u^{1},u^{2}$) are linked to the Gaussian and mean curvatures and
the coefficients of the first and second fundamental forms by the
following relations:
\begin{eqnarray}
\frac{\partial\mathbf{n}}{\partial u^{1}}\cdot\frac{\partial\mathbf{n}}{\partial u^{1}} & = & 2eH-EK\nonumber \\
\frac{\partial\mathbf{n}}{\partial u^{1}}\cdot\frac{\partial\mathbf{n}}{\partial u^{2}} & = & 2fH-FK\\
\frac{\partial\mathbf{n}}{\partial u^{2}}\cdot\frac{\partial\mathbf{n}}{\partial u^{2}} & = & 2gH-GK\nonumber
\end{eqnarray}
where $K$ is the Gaussian curvature, $H$ is the mean curvature,
and $E,F,G,e,f,g$ are the coefficients of the first and second fundamental
forms.

$\bullet$ The above equations of Weingarten (Eq. \ref{eqnr}) can
be solved for the surface basis vectors, $\mathbf{E}_{1}$ and $\mathbf{E}_{2}$,
and hence these vectors can be expressed as combinations of the partial
derivatives of the normal vector, $\mathbf{n}$, that is:
\begin{equation}
\begin{aligned}\mathbf{E}_{1} & =\frac{fF-gE}{b}\frac{\partial\mathbf{n}}{\partial u^{1}}+\frac{fE-eF}{b}\frac{\partial\mathbf{n}}{\partial u^{2}}\\
\mathbf{E}_{2} & =\frac{fG-gF}{b}\frac{\partial\mathbf{n}}{\partial u^{1}}+\frac{fF-eG}{b}\frac{\partial\mathbf{n}}{\partial u^{2}}
\end{aligned}
\end{equation}
where $E,F,G,e,f,g$ are the coefficients of the first and second
fundamental forms and $b=eg-f^{2}$ is the determinant of the surface
covariant curvature tensor.

$\bullet$ From Eq. \ref{eqbCoe}, it can be seen that the coefficients
of the surface covariant curvature tensor, $b_{\alpha\beta}$, are
the projections of the partial derivative of the surface basis vectors,
$\frac{\partial\mathbf{E}_{\alpha}}{\partial u^{\beta}}$, in the
direction of the unit vector normal to the surface, $\mathbf{n}$,
that is:
\begin{equation}
b_{\alpha\beta}=\frac{\partial\mathbf{E}_{\alpha}}{\partial u^{\beta}}\cdot\mathbf{n}
\end{equation}

$\bullet$ The main conclusion from the above equations of Gauss and
Weingarten is that the partial derivatives of $\mathbf{E}_{1}$, $\mathbf{E}_{2}$
and $\mathbf{n}$ can be represented as combinations of these vectors
with coefficients provided by the elements of the first and second
fundamental forms and their partial derivatives.

$\bullet$ For a Monge patch of the form $\mathbf{r}(u,v)=\left(u,v,f(u,v)\right)$,
the Gauss equations are given by:
\begin{eqnarray}
\frac{\partial\mathbf{E}_{1}}{\partial u} & = & \frac{1}{1+f_{u}^{2}+f_{v}^{2}}\left(f_{u}f_{uu}\mathbf{E}_{1}+f_{v}f_{uu}\mathbf{E}_{2}+f_{uu}\sqrt{1+f_{u}^{2}+f_{v}^{2}}\mathbf{n}\right)\nonumber \\
\frac{\partial\mathbf{E}_{1}}{\partial v} & = & \frac{1}{1+f_{u}^{2}+f_{v}^{2}}\left(f_{u}f_{uv}\mathbf{E}_{1}+f_{v}f_{uv}\mathbf{E}_{2}+f_{uv}\sqrt{1+f_{u}^{2}+f_{v}^{2}}\mathbf{n}\right)=\frac{\partial\mathbf{E}_{2}}{\partial u}\\
\frac{\partial\mathbf{E}_{2}}{\partial v} & = & \frac{1}{1+f_{u}^{2}+f_{v}^{2}}\left(f_{u}f_{vv}\mathbf{E}_{1}+f_{v}f_{vv}\mathbf{E}_{2}+f_{vv}\sqrt{1+f_{u}^{2}+f_{v}^{2}}\mathbf{n}\right)\nonumber
\end{eqnarray}
where the subscripts $u$ and $v$  represent partial derivatives
with respect to the surface curvilinear coordinates $u$ and $v$.

$\bullet$ For a Monge patch of the form $\mathbf{r}(u,v)=\left(u,v,f(u,v)\right)$,
the Weingarten equations are given by:
\begin{equation}
\begin{aligned}\frac{\partial\mathbf{n}}{\partial u} & =\frac{\left(f_{u}f_{v}f_{uv}-f_{uu}f_{v}^{2}-f_{uu}\right)\mathbf{E}_{1}+\left(f_{u}f_{v}f_{uu}-f_{u}^{2}f_{uv}-f_{uv}\right)\mathbf{E}_{2}}{\sqrt{\left(1+f_{u}^{2}+f_{v}^{2}\right)^{3}}}\\
\frac{\partial\mathbf{n}}{\partial v} & =\frac{\left(f_{u}f_{v}f_{vv}-f_{uv}f_{v}^{2}-f_{uv}\right)\mathbf{E}_{1}+\left(f_{u}f_{v}f_{uv}-f_{u}^{2}f_{vv}-f_{vv}\right)\mathbf{E}_{2}}{\sqrt{\left(1+f_{u}^{2}+f_{v}^{2}\right)^{3}}}
\end{aligned}
\end{equation}
where the subscripts $u$ and $v$ are as explained in the previous
point.

\subsubsection{Codazzi-Mainardi Equations}

$\bullet$ From the aforementioned equations of Gauss and Weingarten,
supported by further compatibility conditions, the following equations,
called Codazzi or Codazzi-Mainardi equations, can be derived:
\begin{equation}
\begin{aligned}\frac{\partial b_{12}}{\partial u^{1}}-\frac{\partial b_{11}}{\partial u^{2}} & =b_{22}\Gamma_{11}^{2}-b_{12}\left(\Gamma_{12}^{2}-\Gamma_{11}^{1}\right)-b_{11}\Gamma_{12}^{1}\\
\frac{\partial b_{22}}{\partial u^{1}}-\frac{\partial b_{21}}{\partial u^{2}} & =b_{22}\Gamma_{12}^{2}-b_{12}\left(\Gamma_{22}^{2}-\Gamma_{12}^{1}\right)-b_{11}\Gamma_{22}^{1}
\end{aligned}
\label{eqCodazziMainardi}
\end{equation}
where the Christoffel symbols are based on the surface metric. These
equations can be expressed compactly in tensor notation as:
\begin{equation}
\frac{\partial b_{\alpha\beta}}{\partial u^{\delta}}-\frac{\partial b_{\alpha\delta}}{\partial u^{\beta}}=b_{\gamma\beta}\Gamma_{\alpha\delta}^{\gamma}-b_{\gamma\delta}\Gamma_{\alpha\beta}^{\gamma}\label{eqCodazziMainardi2}
\end{equation}

$\bullet$ If we arrange the terms of the Codazzi-Mainardi equation
(Eq. \ref{eqCodazziMainardi2}) and subtract the term $\Gamma_{\delta\beta}^{\gamma}b_{\alpha\gamma}$
from both sides we obtain:

\begin{equation}
\frac{\partial b_{\alpha\beta}}{\partial u^{\delta}}-b_{\gamma\beta}\Gamma_{\alpha\delta}^{\gamma}-\Gamma_{\delta\beta}^{\gamma}b_{\alpha\gamma}=\frac{\partial b_{\alpha\delta}}{\partial u^{\beta}}-b_{\gamma\delta}\Gamma_{\alpha\beta}^{\gamma}-\Gamma_{\delta\beta}^{\gamma}b_{\alpha\gamma}
\end{equation}
which can be expressed compactly, using the covariant derivative notation,
as:
\begin{equation}
b_{\alpha\beta;\gamma}=b_{\alpha\gamma;\beta}\label{eqCodazzibb}
\end{equation}

$\bullet$ The Codazzi-Mainardi equations in the form given by Eq.
\ref{eqCodazzibb} reveals that there are only two independent components
for these equations because, adding to the fact that all the indices
range over 1 and 2, the covariant derivative according to Eq. \ref{eqCodazzibb}
is symmetric in its last two indices (i.e. $\beta$ and $\gamma$),
and the covariant curvature tensor is symmetric in its two indices
(i.e. $b_{\alpha\beta}=b_{\beta\alpha}$).\footnote{As a consequence of these two symmetries, the covariant derivative
of the surface covariant curvature tensor, $b_{\alpha\beta;\gamma}$,
is fully symmetric in all of its indices.} These two independent components are given by:
\begin{equation}
b_{\alpha\alpha;\beta}=b_{\alpha\beta;\alpha}
\end{equation}
where $\alpha\ne\beta$ and there is no summation over $\alpha$.
Writing these equations in full, using the covariant derivative expression\footnote{That is $b_{\alpha\beta;\gamma}=\frac{\partial b_{\alpha\beta}}{\partial u^{\gamma}}-\Gamma_{\alpha\gamma}^{\delta}b_{\delta\beta}-\Gamma_{\beta\gamma}^{\delta}b_{\alpha\delta}$.}
and noting that one term of the covariant derivative expression is
the same on both sides and hence it drops away, we have:
\begin{equation}
\frac{\partial b_{\alpha\alpha}}{\partial u^{\beta}}-\Gamma_{\alpha\beta}^{\delta}b_{\alpha\delta}=\frac{\partial b_{\alpha\beta}}{\partial u^{\alpha}}-\Gamma_{\alpha\alpha}^{\delta}b_{\delta\beta}\,\,\,\,\,\,\,\,\,\text{(\ensuremath{\alpha\ne\beta}, no sum on \ensuremath{\alpha})}
\end{equation}

$\bullet$ There is also another more general equation called the
Gauss-Codazzi equation which is given by:
\begin{equation}
R_{\,\,\gamma\alpha\beta}^{\delta}x_{\delta}^{i}=x_{\delta}^{i}b_{\alpha}^{\delta}b_{\beta\gamma}-x_{\delta}^{i}b_{\beta}^{\delta}b_{\alpha\gamma}+n^{i}b_{\alpha\gamma;\beta}-n^{i}b_{\beta\gamma;\alpha}
\end{equation}
The tangential component of this equation is the \textit{Theorema
Egregium}\footnote{In the form given by Eqs. \ref{eqRb} and \ref{eqRb2}. Refer to $\S$
\ref{subTheoremaEgregium} about the essence of \textit{Theorema Egregium}
as an expression of the fact that certain types of curvature are intrinsic
properties of the surface and hence they can be evaluated from purely
intrinsic parameters obtained from the first fundamental form.} while its normal component is the Codazzi equation.

\subsection{Gauss Mapping\label{subGaussMapping}}

$\bullet$ The Gauss mapping or sphere mapping is a correlation between
the points of a surface and the unit sphere where each point on the
surface is projected onto its unit normal as a point on the unit sphere
centered at the origin of coordinates.\footnote{This sort of mapping for surfaces is similar to the spherical indicatrix
mapping (see $\S$ \ref{subSphericalIndicatrix}) for space curves.}

$\bullet$ In technical terms, let $S$ be a surface embedded in an
$\mathbb{R}^{3}$ space and $S_{1}$ represents the origin-centered
unit sphere in this space, then Gauss mapping is given by:
\begin{equation}
\left\{ N:S\rightarrow S_{1},\,N(P)=\check{P}\right\}
\end{equation}
where the point $P(x,y,z)$ on the trace of $S$ is mapped by $N$
onto the point $\check{P}(\check{x},\check{y},\check{z})$ on the
trace of the unit sphere with $\check{x},\check{y},\check{z}$ being
the coordinates of the origin-based position vector of the normal
vector to the surface, $\mathbf{n}$, at $P$.

$\bullet$ The image $\mathfrak{\bar{S}}$ on the unit sphere of a
Gauss mapping of a region $\mathfrak{S}$ on $S$ is called the spherical
image  of $\mathfrak{S}$.

$\bullet$ To have a single-valued sphere mapping, the functional
relation representing the surface $S$ should be one-to-one.

$\bullet$ The limit of the ratio of the area of a region $\bar{\mathcal{R}}$
on the spherical image to the area of a region $\mathcal{R}$ on the
surface in the neighborhood of a given point $P$ equals the absolute
value of the Gaussian curvature $\left|K\right|$ at $P$ as $\mathcal{R}$
shrinks to the point $P$, that is:
\begin{equation}
\lim_{\mathcal{R}\rightarrow P}\frac{\sigma(\bar{\mathcal{R}})}{\sigma(\mathcal{R})}=\left|K_{P}\right|
\end{equation}
where $\sigma$ stands for area, and $K_{P}$ is the Gaussian curvature
at $P$. The tendency of $\mathcal{R}$ to $P$ should be understood
in the given sense.

$\bullet$ At a given point $P$ on a surface, where the Gaussian
curvature is non-zero, there exists a neighborhood $\mathcal{N}$
of $P$ where an injective mapping can be established between $\mathcal{N}$
and its spherical image $\mathcal{\bar{N}}$.

$\bullet$ A conformal correspondence can be established between a
surface and its spherical image \textit{iff} the surface is a sphere
or a minimal surface (see $\S$ \ref{subMinimalSurface}).

\subsection{Global Surface Theorems }

$\bullet$ In this subsection we state a few global theorems related
to surfaces.

$\bullet$ Planes are the only connected surfaces of class $C^{2}$
whose all points are flat.

$\bullet$ Spheres are the only connected closed surfaces of class
$C^{3}$ whose all points are spherical umbilical.

$\bullet$ Spheres are the only connected compact surfaces of class
$C^{3}$ with constant Gaussian curvature.

$\bullet$ Spheres are the only connected compact surfaces with constant
mean curvature and positive Gaussian curvature.

$\bullet$ Also refer to $\S$ \ref{subGaussBonnetTheorem} for the
global form of the Gauss-Bonnet theorem.

\newpage{}

\section{Curvature}

$\bullet$ ``Curvature'' is a property of both curves and surfaces
at a given point which is determined by the shape of the curve or
surface at that point.\footnote{There are also global characteristics of curvature like total curvature
$K_{t}$ (see $\S$ \ref{subGaussBonnetTheorem}) but they are based
in general on the local characterization of curvature at individual
points.} In this section we investigate this property in its general sense
and examine the main parameters used to describe and quantify it.

\subsection{Curvature Vector\label{subCurvatureVector}}

$\bullet$ At a given point $P$ on a surface $S$, a plane containing
the unit vector $\mathbf{n}$ normal to the surface at $P$ intersects
the surface in a surface curve $C$ having a tangent vector $\mathbf{t}$
at $P$. The curve $C$ is called the normal section of $S$ at $P$
in the direction of $\mathbf{t}$. This curve can be parameterized
by $s$, representing the distance traversed along the curve starting
from an arbitrary given point $P_{0}$, and hence it is defined by
the position vector $\mathbf{r}=\mathbf{r}(s)$. At point $P$, the
vector $\mathbf{T}=\frac{d\mathbf{r}}{ds}$ is a unit vector tangent
to $C$ at $P$ in the direction of increasing $s$ and contained
in the plane tangent to the surface at $P$ as defined previously
in $\S$ \ref{subSpaceCurves}.\footnote{The vector $\mathbf{T}$ can be parallel or anti-parallel to $\mathbf{t}$.
We chose to introduce $\mathbf{t}$ and define it in this way for
generality since the sense of increasing $s$ and hence the curve
orientation can be in one direction or the other.}

$\bullet$ The curvature vector of $C$ at $P$, which is orthogonal
to $\mathbf{T}$, is defined by:\footnote{$\mathbf{K}$ here is the uppercase Greek letter kappa.}
\begin{equation}
\mathbf{K}=\frac{d\mathbf{T}}{ds}\label{eqKT}
\end{equation}

$\bullet$ The curvature, $\kappa$, of $C$ at $P$ (which is defined
previously in $\S$ \ref{subSpaceCurves}) is the magnitude of the
curvature vector, that is: $\kappa=\left|\mathbf{K}\right|$, and
the radius of curvature when $\kappa\ne0$ is its reciprocal, i.e.
$R_{\kappa}=\frac{1}{\kappa}$.\footnote{This is the radius of the osculating circle of $C$ at $P$.}
The curvature vector can then be expressed as:

\begin{equation}
\mathbf{K}=\left|\mathbf{K}\right|\frac{\mathbf{K}}{\left|\mathbf{K}\right|}=\kappa\mathbf{N}\label{eqkapN}
\end{equation}
where $\mathbf{N}$ is the unit vector normal to the curve $C$ at
$P$ as defined previously in $\S$ \ref{subSpaceCurves}.

$\bullet$ The surface curvature vector is independent of the orientation
and parameterization of the surface and the curve. However, in general
it is a function of the position on the surface and the direction
and hence it depends on the point and direction.\footnote{Direction here refers to the tangential direction as represented by
a vector lying in the tangent plane of the surface at a given point
on the surface.}

$\bullet$ A point on the curve at which $\mathbf{K}=\mathbf{0}$,
and hence $\kappa=0$, is called inflection point. At such points,
the radius of curvature is infinite and the normal vector $\mathbf{N}$
and the osculating circle are not defined.\footnote{As we assume that the curve is of class $C^{2}$, the curvature vector
varies smoothly and hence at points of inflection on such a curve
$\mathbf{N}$ may be defined in such a way to ensure continuity when
this is possible, which is not always the case even for some $C^{\infty}$
curves unless the curve is analytic in which case the curve can have
a continuous normal vector in the neighborhood of inflection points.
However, this does not apply to straight lines where the curvature
is identically zero and hence $\mathbf{N}$ is not defined naturally
on any point on the curve although it can be defined as a constant
vector in an arbitrary direction over the whole line although the
above relations between the curve basis vectors will not hold.}

$\bullet$ On introducing a new unit vector, which is orthogonal simultaneously
to $\mathbf{n}$ and $\mathbf{T}$ and defined by the following cross
product:
\begin{equation}
\mathbf{u}=\mathbf{n}\times\mathbf{T}
\end{equation}
the curvature vector, which lies in the plane spanned by $\mathbf{n}$
and $\mathbf{u}$, can then be resolved in the $\mathbf{n}$ and $\mathbf{u}$
directions as:
\begin{equation}
\mathbf{K}=\mathbf{K}_{n}+\mathbf{K}_{g}=\kappa_{n}\mathbf{n}+\kappa_{g}\mathbf{u}\label{eqKnKg}
\end{equation}
where the subscripts $n$ and $g$ are labels and not indices, and
$\kappa_{n}$ and $\kappa_{g}$ are respectively the ``normal curvature''
and the ``geodesic curvature''.

$\bullet$ The normal and geodesic curvatures are given respectively
by:
\begin{equation}
\begin{aligned}\kappa_{n} & =\mathbf{n}\cdot\mathbf{K}=-\mathbf{T}\cdot\frac{d\mathbf{n}}{ds}=-\frac{d\mathbf{r}}{ds}\cdot\frac{d\mathbf{n}}{ds}\\
\kappa_{g} & =\mathbf{u}\cdot\mathbf{K}=\mathbf{u}\cdot\frac{d\mathbf{T}}{ds}=\left(\mathbf{n}\times\mathbf{T}\right)\cdot\frac{d\mathbf{T}}{ds}
\end{aligned}
\label{eqKapnKapg}
\end{equation}

$\bullet$ The vector $\mathbf{u}$, which is a unit vector normal
to the curve $C$, is called the geodesic normal vector. This vector
is the orthonormal projection of $\mathbf{K}$ onto the tangent space
and hence it is contained in the tangent space of the surface at point
$P$.

$\bullet$ While the normal curvature $\kappa_{n}$ is an extrinsic
property, since it depends on the first and second fundamental form
coefficients, the geodesic curvature $\kappa_{g}$ in an intrinsic
property as it depends only on the first fundamental form coefficients
and their derivatives.

$\bullet$ The triad ($\mathbf{n},\mathbf{T},\mathbf{u}$) is another
roaming frame in use in differential geometry in addition to the Frenet
curve-based frame ($\mathbf{T},\mathbf{N},\mathbf{B}$) and the surface-based
frame ($\mathbf{E}_{1},\mathbf{E}_{2},\mathbf{n}$).

$\bullet$ Let $C$ be a curve on a sufficiently smooth surface $S$.
If $\phi$ is the angle between the principal normal vector $N^{i}$
of $C$ at a given point $P$ and the unit vector $n_{i}$ normal
to the surface at $P$ then we have:
\begin{equation}
\cos\phi=n_{i}N^{i}
\end{equation}
The normal and geodesic curvatures, $\kappa_{n}$ and $\kappa_{g}$,
of $C$ at $P$ are then given by:
\begin{equation}
\kappa_{n}=\kappa\cos\phi\,\,\,\,\,\,\,\,\,\,\,\,\,\,\,\,\,\,\,\,\kappa_{g}=\kappa\sin\phi
\end{equation}
where $\kappa$ is the curvature of $C$ at $P$ as defined previously.
According to the theorem of Meusnier (see $\S$ \ref{subMeusnierTheorem}),
all the curves on $S$ that pass through $P$ with the same tangent
direction at $P$ have the same normal curvature.\footnote{The theorem of Meusnier may be stated in this context as: the curvature
of any surface curve at a given point $P$ on the curve equals the
curvature of the normal section which is tangent to the curve at $P$
divided by the cosine of the angle between the principal normal to
the curve at $P$ and the normal to the surface at $P$.}

\subsection{Normal Curvature\label{subNormalCurvature}}

$\bullet$ Using the first and second fundamental forms, given by
Eqs. \ref{eqFundamentalForm1} and \ref{eqFundamentalForm21}, the
normal curvature in the $\frac{du}{dv}$ direction can be expressed
as the following quotient of the second fundamental form involving
the coefficients of the surface covariant curvature tensor to the
first fundamental form involving the coefficients of the surface covariant
metric tensor:\footnote{The last part of the above equation should be interpreted as the sum
of the terms in the numerator divided by the sum of the terms in the
denominator. }
\begin{equation}
\kappa_{n}=\frac{II_{S}}{I_{S}}=\frac{e\dot{u}^{2}+2f\dot{u}\dot{v}+g\dot{v}^{2}}{E\dot{u}^{2}+2F\dot{u}\dot{v}+G\dot{v}^{2}}=\frac{b_{\alpha\beta}du^{\alpha}du^{\beta}}{a_{\gamma\delta}du^{\gamma}du^{\delta}}\label{eqKapnBA}
\end{equation}
where the overdot means differentiation with respect to a general
parameter $t$ of the curve. This equation can be obtained as follow:
\begin{equation}
\begin{aligned}\kappa_{n} & =-\mathbf{T}\cdot\frac{d\mathbf{n}}{ds}\\
 & =-\left(\frac{\partial\mathbf{r}}{\partial u^{\alpha}}\frac{du^{\alpha}}{ds}\right)\cdot\left(\frac{\partial\mathbf{n}}{\partial u^{\beta}}\frac{du^{\beta}}{ds}\right)\\
 & =-\left(\frac{\partial\mathbf{r}}{\partial u^{\alpha}}\cdot\frac{\partial\mathbf{n}}{\partial u^{\beta}}\right)\frac{du^{\alpha}}{ds}\frac{du^{\beta}}{ds}\\
 & =\frac{-\left(\frac{\partial\mathbf{r}}{\partial u^{\alpha}}\cdot\frac{\partial\mathbf{n}}{\partial u^{\beta}}\right)du^{\alpha}du^{\beta}}{ds\,ds}\\
 & =\frac{b_{\alpha\beta}du^{\alpha}du^{\beta}}{(ds)^{2}}\\
 & =\frac{b_{\alpha\beta}du^{\alpha}du^{\beta}}{a_{\gamma\delta}du^{\gamma}du^{\delta}}
\end{aligned}
\end{equation}

$\bullet$ From the previous point, it can be seen that the normal
component of the curvature vector is also given by:
\begin{equation}
\mathbf{K}_{n}=\left[e\left(\frac{du^{1}}{ds}\right)^{2}+2f\frac{du^{1}}{ds}\frac{du^{2}}{ds}+g\left(\frac{du^{2}}{ds}\right)^{2}\right]\mathbf{n}
\end{equation}

$\bullet$ From Eq. \ref{eqKapnBA}, it can be seen that the sign
of the normal curvature $\kappa_{n}$ (i.e. being greater than, less
than or equal to zero) is determined solely by the sign of the second
fundamental form since the first fundamental form is positive definite.
Now, since the sign of the second fundamental form is dependent on
its determinant $b$, this sign is determined by $b$.

$\bullet$ As given earlier, all surface curves passing through a
point $P$ on a surface and have the same tangent line at $P$ have
identical normal curvature at $P$. Hence, the normal curvature is
a property of the surface at a given point and in a given direction.

$\bullet$ The normal curvature $\kappa_{n}$ of a normal section
$C$ of a surface at a point $P$ is equal in magnitude to the curvature
$\kappa$ of $C$ at $P$, i.e. $\left|\kappa_{n}\right|=\kappa$.
This can be explained by the fact that the normal vector $\mathbf{n}$
to the surface at $P$ is collinear with the principal normal vector
$\mathbf{N}$ of $C$ at $P$ so there is only a normal component
to the curvature vector with no geodesic component.

$\bullet$ The normal curvature $\kappa_{n}$ of a surface at a given
point and in a given direction is an invariant property apart from
its sign which is dependent on the choice of the direction of the
unit vector $\mathbf{n}$ normal to the surface.

$\bullet$ At every point on a sphere and in any direction, the normal
curvature is a constant given by: $\left|\kappa_{n}\right|=\frac{1}{R}$
where $R$ is the sphere radius.

$\bullet$ At any point $P$ on a sphere, any surface curve $C$ passing
through $P$ in any direction is a normal section \textit{iff} $C$
is a great circle.\footnote{The great circles of a sphere are the plane sections formed by the
intersection of the sphere with the planes passing through the center
of the sphere.} All these great circles have constant curvature $\kappa$ and normal
curvature $\kappa_{n}$ which are equal in magnitude to $\frac{1}{R}$
where $R$ is the sphere radius.

$\bullet$ At flat points on a surface $\kappa_{n}=0$ in all directions.
At elliptic points $\kappa_{n}\ne0$ and have the same sign in all
directions. At parabolic points $\kappa_{n}$ have the same sign in
all directions except in the direction for which the second fundamental
form vanishes where $\kappa_{n}=0$. At hyperbolic points $\kappa_{n}$
is negative, positive or zero depending on the direction (for the
definition of flat, elliptic, parabolic and hyperbolic points see
$\S$ \ref{subLocalShape}).

$\bullet$ In any two orthogonal directions at a given point $P$
on a sufficiently smooth surface, the sum of the normal curvatures
corresponding to these directions at $P$ is constant.

$\bullet$ At any point $P$ of a sufficiently smooth surface $S$,
there exists a paraboloid tangent at its vertex to the tangent plane
of $S$ at $P$ such that the normal curvature of the paraboloid in
a given direction at $P$ is equal to the normal curvature of $S$
at $P$ in that direction.

$\bullet$ The normal curvature $\kappa_{n}$ is an extrinsic property,
since it depends on the second fundamental form coefficients.\footnote{That is it cannot be expressed purely in terms of the first fundamental
form coefficients.}

$\bullet$ The normal curvatures of surface curves at a given point
$P$ in the directions of the $u$ and $v$ coordinate curves are
given respectively by $\frac{b_{11}}{a_{11}}$ and $\frac{b_{22}}{a_{22}}$
where these are evaluated at $P$.

$\bullet$ As we will see (refer to $\S$ \ref{subPrincipalCurvatures}),
at each non-umbilical point $P$ of a sufficiently smooth surface
there are two perpendicular directions along which the normal curvature
of the surface at $P$ takes its maximum and minimum values of all
the normal curvature values at $P$ (for the definition of umbilical
point see $\S$ \ref{subUmbilicPoints}).

$\bullet$ The necessary and sufficient condition for a given point
$P$ on a sufficiently smooth surface $S$ to be umbilical point is
that the coefficients of the first and second fundamental forms of
the surface are proportional, that is:
\begin{equation}
\frac{e}{E}=\frac{f}{F}=\frac{g}{G}\left(=\kappa_{n}\right)
\end{equation}
where $E,F,G,e,f,g$ are the coefficients of the first and second
fundamental forms at $P$, and $\kappa_{n}$ is the normal curvature
of $S$ at $P$ in any direction.\footnote{This can be seen from Eq. \ref{eqKapnBA} where $\kappa_{n}$ in this
case becomes independent from the direction.}

\subsubsection{Meusnier Theorem\label{subMeusnierTheorem}}

$\bullet$ According to the theorem of Meusnier, all the surface curves
passing through a point $P$ on a surface and have the same tangent
direction\footnote{It should be non-asymptotic direction (see $\S$ \ref{subAsymptoticDirections}).}
at $P$ have identical normal curvatures which is the normal curvature
$\kappa_{n}$ of the normal section at $P$ in the given direction.
Moreover, the osculating circles of these curves lie on a sphere $S_{s}$
with radius $\frac{1}{\kappa}$ and with center at $\mathbf{r}_{C}=\mathbf{r}_{P}+\frac{\mathbf{N}}{\kappa}$
where $\kappa$ (which is equal to $\kappa_{n}$) is the curvature
of the normal section at $P$, $\mathbf{N}$ is the principal normal
vector of the normal section at $P$ and $\mathbf{r}_{P}$ is the
position vector of $P$.

$\bullet$ Following the last point, we have:

(A) The center of the sphere $S_{s}$ is the center of curvature of
the normal section at $P$ in the given direction (see $\S$ \ref{subOsculatingCircle}).

(B) These curves are characterized by being tangent to the normal
section at $P$ in the given direction and by being plane sections
of the surfaces with shared tangent direction at $P$.

(C) The osculating circles of these curves are the intersection of
the sphere $S_{s}$ with the osculating planes of these curves.

(D) The sphere $S_{s}$ is tangent to the tangent plane of the surface
at $P$.

$\bullet$ The theorem of Meusnier may also be stated as follow: the
center of curvature of a surface curve at a given point $P$ on the
curve is obtained by orthogonal projection of the center of curvature
of the normal section, which is tangent to the curve at $P$, on the
osculating plane of the curve.

\subsection{Geodesic Curvature\label{subGeodesicCurvature}}

$\bullet$ As described earlier, the curvature vector $\mathbf{K}$
of a surface curve lies in a plane perpendicular to $\mathbf{T}$
and it can be resolved into a normal component $\mathbf{K}_{n}=\kappa_{n}\mathbf{n}$
and a geodesic component $\mathbf{K}_{g}=\kappa_{g}\mathbf{u}$ where
the normal and geodesic curvatures, $\kappa_{n}$ and $\kappa_{g}$,
are given by Eq. \ref{eqKapnKapg}.

$\bullet$ The geodesic component $\mathbf{K}_{g}$ of the curvature
vector $\mathbf{K}$ of a surface curve at a given point $P$ is the
projection of $\mathbf{K}$ onto the tangent space $T_{P}S$ of the
surface at $P$.

$\bullet$ The geodesic component of the curvature vector of a surface
curve is given by:

\begin{equation}
\mathbf{K}_{g}=\kappa_{g}\mathbf{u}=\left(\frac{d^{2}u^{1}}{ds^{2}}+\Gamma_{\alpha\beta}^{1}\frac{du^{\alpha}}{ds}\frac{du^{\beta}}{ds}\right)\mathbf{E}_{1}+\left(\frac{d^{2}u^{2}}{ds^{2}}+\Gamma_{\alpha\beta}^{2}\frac{du^{\alpha}}{ds}\frac{du^{\beta}}{ds}\right)\mathbf{E}_{2}\label{eqKapgu}
\end{equation}
where the Christoffel symbols are derived from the surface metric.

$\bullet$ While the curvature vector $\mathbf{K}$ is an extrinsic
property of the surface, the geodesic curvature vector $\mathbf{K}_{g}$
is an intrinsic property. This can be seen, for example, from Eq.
\ref{eqKapgu} where all the elements in the equation (Christoffel
symbols and surface basis vectors) are solely dependent on the first
fundamental form.

$\bullet$ The geodesic component of the curvature vector is also
given by:
\begin{equation}
\mathbf{K}_{g}=\left[\mathbf{n}\times\left(\frac{\partial\mathbf{E}_{1}}{\partial u^{1}}\left(\frac{du^{1}}{ds}\right)^{2}+2\frac{\partial\mathbf{E}_{1}}{\partial u^{2}}\frac{du^{1}}{ds}\frac{du^{2}}{ds}+\frac{\partial\mathbf{E}_{2}}{\partial u^{2}}\left(\frac{du^{2}}{ds}\right)^{2}\right)\right]\times\mathbf{n}+\mathbf{E}_{1}\frac{d^{2}u^{1}}{ds^{2}}+\mathbf{E}_{2}\frac{d^{2}u^{2}}{ds^{2}}
\end{equation}

$\bullet$ It can be shown that $\kappa_{g}$ can also be given by:
\begin{equation}
\begin{aligned}\kappa_{g}= & \sqrt{a}\biggl[\Gamma_{11}^{2}\left(\frac{du^{1}}{ds}\right)^{3}+\left(2\Gamma_{12}^{2}-\Gamma_{11}^{1}\right)\left(\frac{du^{1}}{ds}\right)^{2}\frac{du^{2}}{ds}+\\
 & \left(\Gamma_{22}^{2}-2\Gamma_{12}^{1}\right)\frac{du^{1}}{ds}\left(\frac{du^{2}}{ds}\right)^{2}-\Gamma_{22}^{1}\left(\frac{du^{2}}{ds}\right)^{3}+\frac{du^{1}}{ds}\frac{d^{2}u^{2}}{ds^{2}}-\frac{d^{2}u^{1}}{ds^{2}}\frac{du^{2}}{ds}\biggr]
\end{aligned}
\label{eqkapg}
\end{equation}
 where the Christoffel symbols are derived from the surface metric
and $a=EG-F^{2}$ is the determinant of the surface covariant metric
tensor.

$\bullet$ On the $u^{1}$ coordinate curves (i.e. in the orientation
of $\mathbf{E}_{1}$), $\frac{du^{1}}{ds}=\frac{1}{\sqrt{E}}$ and
$\frac{du^{2}}{ds}=0$; hence Eq. \ref{eqkapg} will simplify to:
\begin{equation}
\left(\kappa_{g}\right)_{\mathbf{E}_{1}}=\sqrt{a}\,\Gamma_{11}^{2}\left(\frac{du^{1}}{ds}\right)^{3}=\frac{\sqrt{a}}{E^{3/2}}\,\Gamma_{11}^{2}\label{eqkapgE1}
\end{equation}
The last formula will be simplified further if the surface coordinate
curves $u^{1},u^{2}$ are orthogonal, since in this case $F=0$ and
$\Gamma_{11}^{2}=-\frac{E_{v}}{2G}$ (see Eq. \ref{eqChristoffel2}),
and the formula will become:
\begin{equation}
\left(\kappa_{g}\right)_{\mathbf{E}_{1}}=-\frac{E_{v}}{2E\sqrt{G}}\label{eqkapgE1b}
\end{equation}

$\bullet$ On the $u^{2}$ coordinate curves (i.e. in the orientation
of $\mathbf{E}_{2}$), $\frac{du^{1}}{ds}=0$ and $\frac{du^{2}}{ds}=\frac{1}{\sqrt{G}}$;
hence Eq. \ref{eqkapg} will simplify to:
\begin{equation}
\left(\kappa_{g}\right)_{\mathbf{E}_{2}}=-\sqrt{a}\,\Gamma_{22}^{1}\left(\frac{du^{2}}{ds}\right)^{3}=-\frac{\sqrt{a}}{G^{3/2}}\,\Gamma_{22}^{1}\label{eqkapgE2}
\end{equation}
The last formula will be simplified further if the coordinate curves
$u^{1},u^{2}$ are orthogonal, since in this case $F=0$ and $\Gamma_{22}^{1}=-\frac{G_{u}}{2E}$
(see Eq. \ref{eqChristoffel2}), and the formula will become:
\begin{equation}
\left(\kappa_{g}\right)_{\mathbf{E}_{2}}=\frac{G_{u}}{2G\sqrt{E}}\label{eqkapgE2b}
\end{equation}

$\bullet$ Geodesic curvature can take any real value: positive, negative
or zero. Since $\mathbf{u}=\mathbf{n}\times\mathbf{T}$, the sense
of the geodesic curvature vector depends on the orientation of the
surface and the orientation of the curve.

$\bullet$ Geodesic curvature can also be calculated extrinsically
by:

\begin{equation}
\kappa_{g}=\frac{\mathbf{\ddot{r}}\cdot\left(\mathbf{n}\times\mathbf{\dot{r}}\right)}{\mathbf{\dot{r}}\cdot\mathbf{\dot{r}}}
\end{equation}
where the overdot means differentiation with respect to a general
parameter $t$ for the curve.\footnote{As indicated previously, some intrinsic properties can also be defined
in terms of extrinsic parameters.}

$\bullet$ For any surface curve, the curvature $\kappa$ and the
geodesic curvature $\kappa_{g}$ of the curve at a given point $P$
on the curve are related by:
\begin{equation}
\kappa_{g}=\kappa\sin\theta
\end{equation}
where $\theta$ is the angle between the principal normal $\mathbf{N}$
to the curve at $P$ and the unit normal $\mathbf{n}$ to the surface.

$\bullet$ On a surface patch of class $C^{2}$ with orthogonal coordinate
curves, an $s$-parameterized curve $C$ of class $C^{2}$ given by
$\mathbf{r}\left(u(s),v(s)\right)$ has a geodesic curvature given
by:\footnote{This is known as Liouville formula.}
\begin{equation}
\kappa_{g}=\frac{d\phi}{ds}+\kappa_{u}\cos\phi+\kappa_{v}\sin\phi
\end{equation}
where $\kappa_{u}$ and $\kappa_{v}$ are the geodesic curvatures
of the $u$ and $v$ coordinate curves, $\phi$ is the angle such
that $\mathbf{T}=\frac{\mathbf{E}_{1}}{\left|\mathbf{E}_{1}\right|}\cos\phi+\frac{\mathbf{E}_{2}}{\left|\mathbf{E}_{2}\right|}\sin\phi$
where $\mathbf{T}$ is the tangent unit vector of $C$, and $\mathbf{E}_{1}$
and $\mathbf{E}_{2}$ are the surface basis vectors in the $u$ and
$v$ directions. All the given quantities are evaluated at a given
point on the patch.

\subsection{Principal Curvatures and Directions\label{subPrincipalCurvatures}}

$\bullet$ On rotating the plane containing $\mathbf{n}$, i.e. the
unit vector normal to the surface at a given point $P$ on the surface,
around $\mathbf{n}$ the normal section and hence its curvature $\kappa$
at $P$ will vary in general.\footnote{The plane containing $\mathbf{n}$ is orthogonal to the tangent plane
to the surface at $P$.} The normal curvature $\kappa_{n}$ (which is equal to the curvature
$\kappa$ for a normal section) of the surface at $P$ in a given
direction $\lambda$ is given by:
\begin{equation}
\begin{aligned}\kappa_{n} & =\frac{e+2f\lambda+g\lambda^{2}}{E+2F\lambda+G\lambda^{2}}\\
 & =\frac{f+g\lambda}{F+G\lambda}\\
 & =\frac{e+f\lambda}{E+F\lambda}
\end{aligned}
\label{eqKappan}
\end{equation}
where $E,F,G,e,f,g$ are the coefficients of the first and second
fundamental forms and $\lambda=\frac{du^{2}}{du^{1}}$.\footnote{The directions represented by $\frac{du^{2}}{du^{1}}$ are the directions
of the tangents to the normal sections at $P$.} The two principal curvatures of the surface at $P$, $\kappa_{1}$
and $\kappa_{2}$, which represent the maximum and minimum values
of the normal curvature $\kappa_{n}$ of the surface at $P$ as given
by the last equation, correspond to the two $\lambda$ roots of the
following quadratic equation:
\begin{equation}
(gF-fG)\lambda^{2}+(gE-eG)\lambda+(fE-eF)=0\,\,\,\,\,\,\,\,\,\,\,\,\,\,\,\,\,(gF-fG)\ne0\label{eqlam2}
\end{equation}
The last equation is obtained by equating the derivative of $\kappa_{n}$
(as given by Eq. \ref{eqKappan}) with respect to $\lambda$ to zero
to obtain the extremum values.

$\bullet$ Eq. \ref{eqlam2} possesses two roots, $\lambda_{1}$ and
$\lambda_{2}$, which are linked by the following relations:
\begin{equation}
\lambda_{1}+\lambda_{2}=\frac{gE-eG}{gF-fG}\,\,\,\,\,\,\,\,\,\,\,\,\,\,\,\,\,\,\,\,\,\,\,\,\,\,\,\,\lambda_{1}\lambda_{2}=\frac{fE-eF}{gF-fG}\,\,\,\,\,\,\,\,\,\,\,\,\,\,\,\,\,\text{(\ensuremath{gF-fG\ne0})}
\end{equation}
These roots represent the two directions corresponding to the two
principal curvatures, $\kappa_{1}$ and $\kappa_{2}$, of the surface.

$\bullet$ On an oriented and sufficiently smooth surface, the principal
curvatures are continuous functions of the surface coordinates.

$\bullet$ The following two vectors on the surface, which are defined
in terms of the two roots of $\lambda$ as given above, represent
the directions of the principal curvatures:
\begin{equation}
\begin{aligned}\left(\frac{d\mathbf{r}}{du^{1}}\right)_{1} & =\frac{\partial\mathbf{r}}{\partial u^{1}}+\lambda_{1}\frac{\partial\mathbf{r}}{\partial u^{2}}=\mathbf{E}_{1}+\lambda_{1}\mathbf{E}_{2}\\
\left(\frac{d\mathbf{r}}{du^{1}}\right)_{2} & =\frac{\partial\mathbf{r}}{\partial u^{1}}+\lambda_{2}\frac{\partial\mathbf{r}}{\partial u^{2}}=\mathbf{E}_{1}+\lambda_{2}\mathbf{E}_{2}
\end{aligned}
\end{equation}
These directions, which are called the principal directions\footnote{They may also be called curvature directions.}
of the surface at point $P$, are orthogonal at non-umbilical points
where $\kappa_{1}\ne\kappa_{2}$.\footnote{At umbilical points (see $\S$ \ref{subUmbilicPoints}) the normal
curvature is the same in all directions and hence there are no principal
directions to be orthogonal (or every direction is a principal direction
and hence there is no sensible meaning for being orthogonal).}

$\bullet$ On each non-umbilical (including non-flat) point $P$ of
a smooth surface $S$ an orthonormal moving ``Darboux frame'' can
be defined. This frame consists of the vector triad ($\mathbf{d}_{1},\mathbf{d}_{2},\mathbf{n}$)
where $\mathbf{d}_{1}$ and $\mathbf{d}_{2}$ are the unit vectors
corresponding to the principal directions at $P$, and $\mathbf{n}=\mathbf{d}_{1}\times\mathbf{d}_{2}$
is the unit vector normal to the surface at $P$.\footnote{This is another moving frame in use in differential geometry to be
added to the three previously-described frames: the ($\mathbf{T}$,$\mathbf{N}$,$\mathbf{B}$)
frame, the ($\mathbf{E}_{1},\mathbf{E}_{2},\mathbf{n}$) frame and
the ($\mathbf{n},\mathbf{T},\mathbf{u}$) frame. The first of these
frames is associated with curves while the last three are associated
with surfaces.}

$\bullet$ The centers of curvature of the normal sections corresponding
to the two principal curvatures\footnote{These two normal sections, which are the normal sections in the principal
directions, may be called the principal normal sections.} at a given point $P$ on a surface $S$ are given in tensor notation
by (see $\S$ \ref{subOsculatingCircle}):
\begin{equation}
\begin{aligned}x_{1}^{i} & =x_{P}^{i}+\frac{N^{i}}{\kappa_{1}}\\
x_{2}^{i} & =x_{P}^{i}+\frac{N^{i}}{\kappa_{2}}
\end{aligned}
\end{equation}
where $x_{1}^{i}$ and $x_{2}^{i}$ are the spatial coordinates of
the first and second center of curvature, $x_{P}^{i}$ are the spatial
coordinates of $P$, $N^{i}$ is the principal normal vector of the
normal sections at $P$, $\kappa_{1}$ and $\kappa_{2}$ are the principal
curvatures of $S$ at $P$, and $i=1,2,3$.\footnote{The principal normal vector $\mathbf{N}$ of a normal section at a
point $P$ on the surface is collinear with the unit vector $\mathbf{n}$
normal to the surface at $P$ and hence the principal normal vector
is the same for all the normal sections at $P$.}

$\bullet$ The centers of curvature of the normal sections of the
principal curvatures at $P$ are described as the principal centers
of curvature of the surface at $P$.

$\bullet$ According to one of the Euler theorems, the normal curvature
at a given point $P$ on a surface of class $C^{2}$ in a given direction
can be expressed as a combination of the principal curvatures, $\kappa_{1}$
and $\kappa_{2}$, at $P$ as:
\begin{equation}
\kappa_{n}=\kappa_{1}\cos^{2}\theta+\kappa_{2}\sin^{2}\theta\label{eqEuler}
\end{equation}
where $\theta$ is the angle between the principal direction of $\kappa_{1}$
at $P$ and the given direction.\footnote{Since the principal directions at non-umbilical points are orthogonal,
$\theta$ could represent the angle with the other principal direction
but with relabeling of the two $\kappa$.}

$\bullet$ A number of invariant parameters of the surface at a given
point $P$ on the surface are defined in terms of the principal curvatures;
these include:

(A) The principal radii of curvature: $R_{1}=\left|\frac{1}{\kappa_{1}}\right|$
and $R_{2}=\left|\frac{1}{\kappa_{2}}\right|$.

(B) The Gaussian curvature:\footnote{Some authors use ``total curvature'' for the ``Gaussian curvature''
and hence they are synonym, while others use total curvature for the
area integral $\int Kd\sigma$ as, for example, in the Gauss-Bonnet
theorem (refer to $\S$ \ref{subGaussBonnetTheorem}). In the present
text we follow the latter convention and we use total curvature strictly
for the integral. Hence, we label the Gaussian curvature with $K$
and the total curvature with $K_{t}$.} $K=\kappa_{1}\kappa_{2}$.

(C) The mean curvature: $H=\frac{\kappa_{1}+\kappa_{2}}{2}$.\footnote{Some authors define $H$ as the sum of $\kappa_{1}$ and $\kappa_{2}$,
that is: $H=\kappa_{1}+\kappa_{2}$. Each convention has its merit.
In the present notes we define $H$ as the average, not the sum, of
the two principal curvatures.}

$\bullet$ Table \ref{tabCurvature} shows the restricting conditions
on the principal curvatures, $\kappa_{1}$ and $\kappa_{2}$, for
a number of common surfaces with simple geometric shapes and the effect
on the Gaussian curvature $K$ and the mean curvature $H$.

\begin{table}
\centering{}\caption{The limiting conditions on the principal curvatures, $\kappa_{1}$
and $\kappa_{2}$, for a number of surfaces of simple geometric shapes
alongside the corresponding mean curvature $H$ and Gaussian curvature
$K$. Apart from the plane, the unit vector normal to the surface,
$\mathbf{n}$, is assumed to be in the outside direction.\vspace{0.2cm}\label{tabCurvature}}
\begin{tabular*}{16cm}{@{\extracolsep{\fill}}|l|c|c|c|c|}
\hline
 & $\kappa_{1}$ & $\kappa_{2}$ & $H$ & $K$\tabularnewline
\hline
Plane & 0 & 0 & 0 & 0\tabularnewline
\hline
Cylinder & $\kappa_{1}=0$ & $\kappa_{2}<0$ & $H<0$ & 0\tabularnewline
\hline
Sphere & $\kappa_{1}=\kappa_{2}<0$ & $\kappa_{2}=\kappa_{1}<0$ & $H<0$ & $K>0$\tabularnewline
\hline
Ellipsoid & $\kappa_{1}<0$ & $\kappa_{2}<0$ & $H<0$ & $K>0$\tabularnewline
\hline
Hyperboloid of one sheet (see $\S$ \ref{HyperboloidOneSheet}) & $\kappa_{1}>0$ & $\kappa_{2}<0$ & --- & $K<0$\tabularnewline
\hline
\end{tabular*}
\end{table}

$\bullet$ Let $P$ be a point on a sufficiently smooth surface $S$
embedded in a 3D space coordinated by a rectangular Cartesian system
($x,y,z$) with $P$ being above the origin, the tangent plane of
$S(x,y)$ at $P$ being parallel to the $xy$ plane, and the principal
directions being along the $x$ and $y$ coordinate lines. The equation
of $S$ in the neighborhood of $P$ can then be expressed, up and
including the quadratic terms, in the following form:
\begin{equation}
S(x,y)\simeq S(0,0)+\frac{\kappa_{1}x^{2}}{2}+\frac{\kappa_{2}y^{2}}{2}
\end{equation}
where $\kappa_{1}$ and $\kappa_{2}$ are the principal curvatures
of $S$ at $P$. This means that in the immediate neighborhood of
$P$, $S$ resembles a quadratic surface of the given form.\footnote{This includes umbilical points where the principal directions can
be arbitrarily chosen as the directions of the $x$ and $y$ coordinate
lines.}

$\bullet$ The necessary and sufficient condition for a number $\kappa\in\mathbb{R}$
to be a principal curvature of a surface $S$ at a given point $P$
and in a given direction $\frac{du}{dv}$, where $\left(du\right)^{2}+\left(dv\right)^{2}\ne0$,
is that the following equations are satisfied:\footnote{For simplicity in notation, we use $\kappa$ in the following equations
to represent principal curvature. This use should not be confused
with the curve curvature which is also symbolized by $\kappa$. However,
for this case $\kappa$ is equal to the principal curvature since
the latter is the curvature of a normal section and hence the use
of $\kappa$ is justified.}
\begin{equation}
\begin{aligned}(e-\kappa E)du+(f-\kappa F)dv & =0\\
(f-\kappa F)du+(g-\kappa G)dv & =0
\end{aligned}
\end{equation}
where $E,F,G,e,f,g$ are the coefficients of the first and second
fundamental forms at $P$.

$\bullet$ The above equations can be cast in a matrix form as:
\begin{equation}
\left[\begin{array}{cc}
e-\kappa E & f-\kappa F\\
f-\kappa F & g-\kappa G
\end{array}\right]\left[\begin{array}{c}
du\\
dv
\end{array}\right]=\left[\begin{array}{c}
0\\
0
\end{array}\right]
\end{equation}
This system of homogeneous linear equations has a non-trivial solution,
$du$ and $dv$, \textit{iff} the determinant of the coefficient matrix
is zero, that is:
\begin{equation}
\left|\begin{array}{cc}
e-\kappa E & f-\kappa F\\
f-\kappa F & g-\kappa G
\end{array}\right|=\left(EG-F^{2}\right)\kappa^{2}-\left(gE-2fF+eG\right)\kappa+\left(eg-f^{2}\right)=0\label{eqCurvEq}
\end{equation}
The above quadratic equation in $\kappa$ has a non-negative discriminant
and hence it possesses either two distinct real roots or a repeated
real root. In the former case there are two distinct principal curvatures
at $P$ corresponding to two orthogonal principal directions, while
in the latter case the point is umbilical and hence there is no specific
principal direction as each direction can be a principal direction.
So in brief, a given real number $\kappa$ is a principal curvature
of $S$ at $P$ \textit{iff} it is a solution of Eq. \ref{eqCurvEq}.

$\bullet$ From Eq. \ref{eqCurvEq}, the principal curvatures of a
surface at a given point $P$ are the solutions of the above quadratic
equation and hence they are given by:

\begin{equation}
\kappa_{1,2}=\frac{gE-2fF+eG\pm\sqrt{\left(gE-2fF+eG\right)^{2}-4\left(EG-F^{2}\right)\left(eg-f^{2}\right)}}{2\left(EG-F^{2}\right)}\label{eqkap12}
\end{equation}
where $E,F,G,e,f,g$ are the coefficients of the first and second
fundamental forms at $P$.

$\bullet$ On dividing Eq. \ref{eqCurvEq} by $a=EG-F^{2}$ we obtain:
\begin{equation}
\kappa^{2}-2H\kappa+K=0\label{eqkapHK}
\end{equation}
where $H$ and $K$ are the mean and Gaussian curvatures (see Eqs.
\ref{eqKRa} and \ref{eqHtrW}). Hence, Eq. \ref{eqkap12} can be
expressed compactly as:\footnote{This formula can be obtained directly from Eq. \ref{eqkapHK} using
the quadratic formula.}

\begin{equation}
\kappa_{1,2}=H\pm\sqrt{H^{2}-K}
\end{equation}

$\bullet$ The principal directions are invariant with respect to
permissible coordinate transformations and parameterizations.

$\bullet$ When the $u$ and $v$ coordinate curves of a surface at
a point $P$ are aligned along the principal directions at $P$, the
principal curvatures at $P$ will be given by:
\begin{equation}
\kappa_{1}=\frac{b_{11}}{a_{11}}=\frac{e}{E}\,\,\,\,\,\,\,\,\,\,\,\,\,\,\,\,\,\,\,\,\,\,\,\,\,\kappa_{2}=\frac{b_{22}}{a_{22}}=\frac{g}{G}
\end{equation}
where $E,G,e,g$ are the coefficients of the first and second fundamental
forms at $P$, and the indexed $a$ and $b$ are the coefficients
of the surface covariant metric and covariant curvature tensors.

$\bullet$ For a non-umbilical point $P$ on a sufficiently smooth
surface $S$, the direction $\frac{du}{dv}$ is a principal direction
of $S$ at $P$ \textit{iff} the following condition is true:\footnote{This equation is based on Eq. \ref{eqlam2}.}
\begin{equation}
\left(fE-eF\right)dudu+\left(gE-eG\right)dudv+\left(gF-fG\right)dvdv=0\label{eqPrincipalDirection}
\end{equation}
The last equation can be factored into two linear equations each of
the form $Adu+Bdv=0$ (with $A$ and $B$ being real parameters) where
these equations represent the two orthogonal principal directions.

$\bullet$ At a given point $P$ on a sufficiently smooth surface
$S$, a direction $\frac{du}{dv}$ is a principal direction \textit{iff}
for a real number $\kappa$ the following relation holds true:
\begin{equation}
d\mathbf{n}=-\kappa d\mathbf{r}\label{eqRodrigues}
\end{equation}
where $d\mathbf{n}=\frac{\partial\mathbf{n}}{\partial u}du+\frac{\partial\mathbf{n}}{\partial v}dv$
and $d\mathbf{r}=\frac{\partial\mathbf{r}}{\partial u}du+\frac{\partial\mathbf{r}}{\partial v}dv$.
If these conditions are satisfied, then $\kappa$ is the principal
curvature of $S$ at $P$ corresponding to the principal direction
$\frac{du}{dv}$. Eq. \ref{eqRodrigues} is called the Rodrigues curvature
formula.

$\bullet$ From the Rodrigues curvature formula (Eq. \ref{eqRodrigues})
we have:

\begin{equation}
\frac{\partial\mathbf{n}}{\partial u}=-\kappa\mathbf{E}_{1}\,\,\,\,\,\,\,\,\,\,\,\,\,\,\,\,\,\,\,\,\,\,\,\,\,\,\,\,\,\,\,\,\,\,\,\,\,\,\frac{\partial\mathbf{n}}{\partial v}=-\kappa\mathbf{E}_{2}
\end{equation}

$\bullet$ At any point on a plane surface, all the directions are
principal directions.\footnote{Or, alternatively, there is no principal direction depending on allowing
more than two principal directions or not.}

$\bullet$ At any point on a sphere, all the directions are principal
directions (or there is no principal direction).

$\bullet$ The principal curvatures, $\kappa_{1}$ and $\kappa_{2}$,
are the eigenvalues of the mixed type surface curvature tensor $b_{\beta}^{\alpha}$.

\subsubsection{Gaussian Curvature\label{subGaussianCurvature}}

$\bullet$ Gaussian curvature represents a generalization of curve
curvature to surfaces since it is the product of two curvatures of
curves belonging to the surface and hence in a sense it is a ``2D
curvature''.

$\bullet$ The Gaussian curvature $K$ of a surface\footnote{Gaussian curvature may also be called the Riemannian curvature of
the surface.} at a given point $P$ on the surface is given by:
\begin{equation}
K=\frac{eg-f^{2}}{EG-F^{2}}=\frac{b}{a}=\frac{R_{1212}}{a}\label{eqKRa}
\end{equation}
where $E,F,G$ and $e,f,g$ are the coefficients of the first and
second fundamental forms at $P$, $a$ and $b$ are the determinants
of the surface covariant metric and curvature tensors, and $R_{1212}$
is the component of the 2D covariant Riemann-Christoffel curvature
tensor.\footnote{As discussed in \cite{SochiTC2}, the 2D Riemann-Christoffel curvature
tensor has only one independent non-vanishing component which is represented
by $R_{1212}$. From Eq. \ref{eqRb2} we get $R_{1212}=b_{11}b_{22}-b_{12}b_{21}=eg-f^{2}=b$
where the indexed $b$ are the coefficients of the surface covariant
curvature tensor and $b$ is its determinant.}

$\bullet$ The above formulae (Eq. \ref{eqKRa}) are based on the
fact that the Gaussian curvature $K$ is the determinant of the mixed
curvature tensor $b_{\beta}^{\alpha}$ of the surface, that is:\footnote{Being the determinant of a tensor establishes the status of $K$ as
an invariant under permissible transformations. Also, the chain of
formulae in Eq. \ref{eqKba} may be taken in the opposite direction
starting primarily from $K=\frac{b}{a}$ as a definition or as a derived
result from other arguments, and hence the statement $K=\mathrm{det}(b_{\beta}^{\alpha})$
will be obtained as a secondary result, as done by some authors.}
\begin{equation}
K=\mathrm{det}(b_{\beta}^{\alpha})=\mathrm{det}(a^{\alpha\gamma}b_{\gamma\beta})=\mathrm{det}(a^{\alpha\gamma})\mathrm{det}(b_{\gamma\beta})=\frac{\mathrm{det}(b_{\gamma\beta})}{\mathrm{det}(a_{\alpha\gamma})}=\frac{b}{a}\label{eqKba}
\end{equation}

$\bullet$ From Eq. \ref{eqKRa}, we see that the sign of $K$ (i.e.
$K>0$, $K<0$ or $K=0$) is the sign of $b$ since $a$ is positive
definite.

$\bullet$ Since both $R_{1212}$ (see Eq. \ref{eqRieChrTensor})
and $a$ depend exclusively on the surface metric tensor, Eq. \ref{eqKRa}
reveals that $K$ depends only on the first fundamental form coefficients
and hence it is an intrinsic property of the surface (refer to $\S$
\ref{subTheoremaEgregium}). The dependence of $K$ on the second
fundamental form coefficients in Eq. \ref{eqKba} does not affect
its qualification as an intrinsic property since this dependency is
not indispensable as $K$ can be expressed in terms of the first fundamental
form coefficients exclusively.

$\bullet$ The Gaussian curvature is an invariant with respect to
permissible coordinate transformations in 2D manifolds and hence:
\begin{equation}
K=\frac{R_{1212}}{a}=\frac{\bar{R}_{1212}}{\bar{a}}
\end{equation}
where the barred and unbarred symbols represent the quantities in
the barred and unbarred surface coordinate systems.

$\bullet$ The Gaussian curvature is also invariant with respect to
the type of representation and parameterization of the surface. In
particular, the Gaussian curvature is independent of the orientation
of the surface which is based on the choice of the direction of the
normal vector $\mathbf{n}$ to the surface.

$\bullet$ From Table \ref{tabCurvature} we see that the Gaussian
curvature of planes and cylinders are both identically zero. At the
root of this is the fact that the Gaussian curvature is an intrinsic
property and the cylinder is a developable surface obtained by wrapping
a plane with no localized distortion by stretching or compression.
Hence, the planes and cylinders possess identical first fundamental
forms, as indicated previously in $\S$ \ref{subFirstFundamentalForm},
and consequently identical Gaussian curvature (also see $\S$ \ref{subTheoremaEgregium}).

$\bullet$ Since the magnitude of the normal curvature of a sphere
of radius $R$ is $\left|\kappa_{n}\right|=\frac{1}{R}$ at any point
on its surface and for any normal section in any direction, its Gaussian
curvature is a constant given by $K=\frac{1}{R^{2}}$.

$\bullet$ For a Monge patch of the form $\mathbf{r}(u,v)=\left(u,v,f(u,v)\right)$,
the Gaussian curvature is given by:
\begin{equation}
K=\frac{f_{uu}f_{vv}-f_{uv}^{2}}{\left(1+f_{u}^{2}+f_{v}^{2}\right)^{2}}
\end{equation}
where the subscripts $u$ and $v$ stand for partial derivatives with
respect to these surface coordinates.

$\bullet$ The Gaussian curvature of a surface of revolution generated
by revolving a plane curve of class $C^{2}$ having the form $y=f(x)$
around the $x$-axis is given by:
\begin{equation}
K=-\frac{f_{xx}}{f\left(1+f_{x}^{2}\right)^{2}}
\end{equation}
where the subscript $x$ represents the partial derivative of $f$
with respect to this variable.

$\bullet$ At any point on a sufficiently smooth surface the Gaussian
curvature satisfies the following relation:

\begin{equation}
\partial_{u}\mathbf{n}\times\partial_{v}\mathbf{n}=K\left(\mathbf{E}_{1}\times\mathbf{E}_{2}\right)
\end{equation}
On dot producting both sides by $\mathbf{n}$ we obtain:

\begin{equation}
\mathbf{n}\cdot\left(\partial_{u}\mathbf{n}\times\partial_{v}\mathbf{n}\right)=K\,\mathbf{n}\cdot\left(\mathbf{E}_{1}\times\mathbf{E}_{2}\right)=K\sqrt{a}
\end{equation}
Hence:

\begin{equation}
K=\frac{\mathbf{n}\cdot\left(\partial_{u}\mathbf{n}\times\partial_{v}\mathbf{n}\right)}{\sqrt{a}}
\end{equation}

$\bullet$ On a 2D surface, the Gaussian curvature $K$ is related
to the Ricci curvature scalar $R$ (see $\S$ \ref{subRiemannTensor})
by the following relation:
\begin{equation}
K=\frac{R}{2}
\end{equation}

$\bullet$ There are surfaces with constant zero Gaussian curvature
$K=0$ (e.g. planes, cylinders and cones excluding the apex), surfaces
with constant positive Gaussian curvature $K>0$ (e.g. spheres with
$K=\frac{1}{R^{2}}$ where $R$ is the sphere radius), and surfaces
with constant negative Gaussian curvature $K<0$ (e.g. Beltrami pseudo-spheres
with $K=-\frac{1}{\rho^{2}}$ where $\rho$ is the pseudo-radius of
the pseudo-sphere).\footnote{The tractrix is a plane curve in the $xy$ plane starting from point
($\rho,0$) on the $x$-axis ($\rho>0$) with the property that the
length of the line segment of its tangent between the tangency point
and the point of intersection with the $y$-axis is equal $\rho$.
Hence, the tractrix is a solution of the following differential equation:
\begin{equation}
\frac{dy}{dx}=-\frac{\sqrt{\rho^{2}-x^{2}}}{x}
\end{equation}
with the condition $y(\rho)=0$. The real parameter $\rho$ is called
pseudo-radius. The Beltrami pseudo-sphere is a surface of revolution
generated by revolving a tractrix around its asymptote which is the
$y$-axis in the above formulation.} However, in general the Gaussian curvature is a variable function
in sign and magnitude of the surface coordinates and a single surface
can have Gaussian curvature of different magnitudes and signs.

$\bullet$ On scaling a surface up or down by a constant factor $c>0$
(refer to $\S$ \ref{MetricScaling}), the Gaussian curvature $K$
will scale by a factor of $\frac{1}{c^{2}}$.

$\bullet$ The Gaussian curvature is invariant with respect to all
isometric transformations, and hence two isometric surfaces have identical
Gaussian curvature at each pair of their corresponding points. However,
two mapped surfaces with equal Gaussian curvature at their corresponding
points are not necessarily isometric. Yes in the case of two sufficiently
smooth surfaces with equal constant Gaussian curvature the two surfaces
have local isometry i.e. they are isometric in the immediate neighborhood
of their corresponding points.

$\bullet$ In 3D manifolds, there is no compact surface of class $C^{2}$
with non-positive Gaussian curvature, i.e. $K\le0$.

$\bullet$ The sphere is the only connected, compact and sufficiently
smooth surface with constant Gaussian curvature.

$\bullet$ According to the Hilbert lemma, if $P$ is a point on a
sufficiently smooth surface $S$ with $\kappa_{1}$ and $\kappa_{2}$
being the principal curvatures of $S$ at $P$ such that: $\kappa_{1}>\kappa_{2}$,
$\kappa_{1}$ is a local maximum, and $\kappa_{2}$ is a local minimum,
then the Gaussian curvature of $S$ at $P$ is non-positive, that
is $K\le0$.

$\bullet$ At a given point $P$ on a spherically-mapped (see $\S$
\ref{subGaussMapping}) and sufficiently smooth surface $S$, the
ratio of the area of the spherical image $\mathfrak{\bar{S}}$ of
a region surrounding $P$ to the area of the mapped region $\mathfrak{S}$
on $S$ converges to the value of the Gaussian curvature at $P$ as
$\mathfrak{S}$ shrinks to $P$.

$\bullet$ From the Gauss-Bonnet theorem (see $\S$ \ref{subGaussBonnetTheorem}),
it can be shown that a surface will have identically-vanishing Gaussian
curvature if at any point $P$ on the surface there are two families
of geodesic curves (see $\S$ \ref{subGeodesicCurves}) in the neighborhood
of $P$ intersecting at a constant angle.

$\bullet$ From Eqs. \ref{eqR1212} and \ref{eqKRa}, the Gaussian
curvature $K$ of a sufficiently smooth surface represented by $\mathbf{r}=\mathbf{r}(u,v)=\mathbf{r}(u^{1},u^{2})$
can also be given by:
\begin{equation}
K=\frac{1}{a}\left[F_{uv}-\frac{1}{2}E_{vv}-\frac{1}{2}G_{uu}+a_{\alpha\beta}\left(\Gamma_{12}^{\alpha}\Gamma_{12}^{\beta}-\Gamma_{11}^{\alpha}\Gamma_{22}^{\beta}\right)\right]\,\,\,\,\,\,\,\,\,\,\,\text{(\ensuremath{\alpha},\ensuremath{\beta}=1,2)}
\end{equation}
where $E,F,G$ are the coefficients of the first fundamental form,
the subscripts $u$ and $v$ stand for partial derivatives with respect
to these surface coordinates, $a$ is the determinant of the surface
covariant metric tensor and the indexed $a$ are its coefficients.
The Christoffel symbols are based on the surface metric.

$\bullet$ The Gaussian curvature of a sufficiently smooth surface
is also given by:
\begin{equation}
K=\frac{1}{2\sqrt{a}}\left[\partial_{u}\left(\frac{FE_{v}}{E\sqrt{a}}-\frac{G_{u}}{\sqrt{a}}\right)+\partial_{v}\left(\frac{2F_{u}}{\sqrt{a}}-\frac{E_{v}}{\sqrt{a}}-\frac{FE_{u}}{E\sqrt{a}}\right)\right]
\end{equation}
where $E,F,G$ are the coefficients of the first fundamental form,
$a=EG-F^{2}$ is the determinant of the surface covariant metric tensor,
and the subscripts $u$ and $v$ stand for partial derivatives with
respect to these surface coordinates.

$\bullet$ The Gaussian curvature $K$ of a surface of class $C^{3}$
represented by $\mathbf{r}(u,v)$ with orthogonal surface coordinate
curves is given by:
\begin{equation}
K=-\frac{1}{2\sqrt{EG}}\left[\partial_{u}\left(\frac{G_{u}}{\sqrt{EG}}\right)+\partial_{v}\left(\frac{E_{v}}{\sqrt{EG}}\right)\right]
\end{equation}
where $E,G$ are the coefficients of the first fundamental form, and
the subscripts $u$ and $v$ stand for partial derivatives with respect
to these surface coordinates. This formula is obtained from the previous
formula by setting $F=0$.

$\bullet$ The last formula will simplify to:
\begin{equation}
K=-\frac{\partial_{uu}\sqrt{G}}{\sqrt{G}}
\end{equation}
when the surface $\mathbf{r}(u,v)$ is represented by geodesic coordinates
(see $\S$ \ref{subGeodesicCoordinates}) with the $u$ coordinate
curves being geodesics and $u$ is a natural parameter.\footnote{Here, geodesic coordinates means a coordinate system on a coordinate
patch of a surface whose $u$ and $v$ coordinate curve families are
orthogonal with one of these families ($u$ or $v$) being a family
of geodesic curves (refer to $\S$ \ref{subGeodesicCoordinates} for
more details).}

$\bullet$ The sign of the Gaussian curvature is independent of the
choice of the direction of the unit vector $\mathbf{n}$ normal to
the surface. This is because a change in the direction ensuing a change
in the sign of the principal curvatures will change the sign of both
of these curvatures and hence their product will not be affected.

$\bullet$ The Gaussian curvature $K$ can be expressed in terms of
the mean curvature $H$ as:
\begin{equation}
K=\kappa_{1}\kappa_{2}=\left(H+C\right)\left(H-C\right)=H^{2}-C^{2}
\end{equation}
where
\begin{equation}
C=\frac{\sqrt{\left\{ e(EG-2F^{2})+2fEF-gE^{2}\right\} ^{2}+4\left(fE-eF\right)^{2}}}{2E\left(EG-F^{2}\right)}\label{eqC}
\end{equation}
and $E,F,G,e,f,g$ are the coefficients of the first and second fundamental
forms.

$\bullet$ The Gaussian curvature of a surface $S$ at a given point
$P$ on the surface is positive if all the surface points in a deleted
neighborhood of $P$ on $S$ are on the same side of the plane tangent
to $S$ at $P$. The Gaussian curvature is negative if for all deleted
neighborhoods of $P$ on $S$ some points are on one side and some
are on the other. The Gaussian curvature is zero if, in a deleted
neighborhood, either all the points lie in the tangent plane or all
the points are on one side except some which lie on a curve in the
tangent plane. Hence,

(A) A sphere has a positive Gaussian curvature at all points.

(B) A hyperbolic paraboloid has a negative Gaussian curvature at all
points.

(C) A plane has a zero Gaussian curvature at all points.

(D) A cylinder has a zero Gaussian curvature at all points.

(E) A torus has points with positive Gaussian curvature (outer half),
points with zero Gaussian curvature (top and bottom circles) and points
with negative Gaussian curvature (inner half).

$\bullet$ The Gaussian curvature of a developable surface is identically
zero. Hence, beside the plane, there are other surfaces with constant
zero Gaussian curvature. Examples are: cones, cylinders and tangent
surfaces of curves (refer to $\S$ \ref{subTangentSurfaces}).

$\bullet$ Examples of the Gaussian curvature, $K$, for a number
of simple surfaces:

(A) Plane: $K=0$.

(B) Sphere of radius $R$: $K=\frac{1}{R^{2}}$.

(C) Torus parameterized by $x=\left(R+r\sin\phi\right)\cos\theta$,
$y=\left(R+r\sin\phi\right)\sin\theta$ and $z=r\cos\phi$: $K=\frac{\sin\phi}{r\left(R+r\sin\phi\right)}$.\footnote{In these parameterizations, $R$ is the torus radius (i.e. the distance
between the center of symmetry of the torus and the center of the
tube), $r$ is the radius of the generating circle ($r<R$), $\phi\in[0,2\pi)$
is the angle of variation of $r$, and $\theta\in[0,2\pi)$ is the
angle of variation of $R$.}

$\bullet$ The total curvature $K_{t}$ is the area integral of the
Gaussian curvature $K$ over a surface $S$ or a patch of a surface,
that is:
\begin{equation}
K_{t}=\iint_{S}Kd\sigma
\end{equation}
where $d\sigma$ symbolizes infinitesimal area element on the surface.

$\bullet$ From Eqs. \ref{eqdsegma} and \ref{eqparnparn}, the total
curvature $K_{t}$ may be given by:
\begin{equation}
K_{t}\equiv\iint_{S}Kd\sigma=\iint_{S}K\left|\mathbf{E}_{1}\times\mathbf{E}_{2}\right|dudv=\iint_{S}\mathrm{sgn}\left(K\right)\left|\partial_{u}\mathbf{n}\times\partial_{v}\mathbf{n}\right|dudv
\end{equation}
where $\mathrm{sgn}\left(K\right)$ is the sign function of $K$ as
a function of the surface coordinates, $u$ and $v$.

$\bullet$ The Riemann-Christoffel curvature tensor is related to
the Gaussian curvature through the alternating absolute tensor of
the surface by the following relation:
\begin{equation}
R_{\alpha\beta\gamma\delta}=K\underline{\epsilon}_{\alpha\beta}\underline{\epsilon}_{\gamma\delta}\label{eqRK}
\end{equation}
On multiplying both sides of the last equation by $\underline{\epsilon}^{\alpha\beta}\underline{\epsilon}^{\gamma\delta}$
we get:
\begin{equation}
\underline{\epsilon}^{\alpha\beta}\underline{\epsilon}^{\gamma\delta}R_{\alpha\beta\gamma\delta}=K\underline{\epsilon}^{\alpha\beta}\underline{\epsilon}^{\gamma\delta}\underline{\epsilon}_{\alpha\beta}\underline{\epsilon}_{\gamma\delta}
\end{equation}
Now, since $\underline{\epsilon}^{\alpha\beta}\underline{\epsilon}_{\alpha\beta}=\underline{\epsilon}^{\gamma\delta}\underline{\epsilon}_{\gamma\delta}=2$,
the last equation becomes:\footnote{Hence, $K$ is a rank-0 tensor or scalar.}
\begin{equation}
K=\frac{1}{4}\underline{\epsilon}^{\alpha\beta}\underline{\epsilon}^{\gamma\delta}R_{\alpha\beta\gamma\delta}=\frac{1}{4}\underline{\epsilon}^{\alpha\beta}\underline{\epsilon}^{\gamma\delta}\left(b_{\alpha\gamma}b_{\beta\delta}-b_{\alpha\delta}b_{\beta\gamma}\right)\label{eqKR}
\end{equation}

$\bullet$ The Riemann-Christoffel curvature tensor is also linked
to the Gaussian curvature through the surface metric tensor by the
following relation:
\begin{equation}
R_{\alpha\beta\gamma\delta}=K\left(a_{\alpha\gamma}a_{\beta\delta}-a_{\alpha\delta}a_{\beta\gamma}\right)
\end{equation}

$\bullet$ The Gaussian curvature $K$ may also be given by:
\begin{equation}
K=\frac{1}{2}\underline{\epsilon}^{\alpha\beta}\underline{\epsilon}^{\gamma\delta}b_{\gamma\alpha}b_{\delta\beta}
\end{equation}
where the indexed $\underline{\epsilon}$ are the 2D absolute alternating
tensor, the indexed $b$ are the components of the surface covariant
curvature tensor and all the indices range over 1 and 2.

$\bullet$ From Eqs. \ref{eqRb2} and \ref{eqRK}, the Gaussian curvature
and the surface curvature tensor are related by:
\begin{equation}
K\underline{\epsilon}_{\alpha\beta}\underline{\epsilon}_{\gamma\delta}=b_{\alpha\gamma}b_{\beta\delta}-b_{\alpha\delta}b_{\beta\gamma}
\end{equation}

$\bullet$ Another formula for the Gaussian curvature (in terms of
the surface basis vectors, their derivatives and the coefficients
of the first fundamental form) can be obtained from the formula $K=\frac{b}{a}$,
that is:
\begin{equation}
\begin{aligned}b & =eg-f^{2}\\
 & =\frac{\left(\partial_{u}\mathbf{E}_{1}\cdot\mathbf{E}_{1}\times\mathbf{E}_{2}\right)\left(\partial_{v}\mathbf{E}_{2}\cdot\mathbf{E}_{1}\times\mathbf{E}_{2}\right)-\left(\partial_{v}\mathbf{E}_{1}\cdot\mathbf{E}_{1}\times\mathbf{E}_{2}\right)^{2}}{\left|\mathbf{E}_{1}\times\mathbf{E}_{2}\right|^{2}}\\
 & =\frac{\left(\partial_{u}\mathbf{E}_{1}\cdot\mathbf{E}_{1}\times\mathbf{E}_{2}\right)\left(\partial_{v}\mathbf{E}_{2}\cdot\mathbf{E}_{1}\times\mathbf{E}_{2}\right)-\left(\partial_{v}\mathbf{E}_{1}\cdot\mathbf{E}_{1}\times\mathbf{E}_{2}\right)^{2}}{a}\\
 & =\frac{\left(\partial_{u}\mathbf{E}_{1}\cdot\mathbf{E}_{1}\times\mathbf{E}_{2}\right)\left(\partial_{v}\mathbf{E}_{2}\cdot\mathbf{E}_{1}\times\mathbf{E}_{2}\right)-\left(\partial_{v}\mathbf{E}_{1}\cdot\mathbf{E}_{1}\times\mathbf{E}_{2}\right)^{2}}{EG-F^{2}}
\end{aligned}
\end{equation}
Hence:
\begin{equation}
\begin{aligned}K & =\frac{b}{a}\\
 & =\frac{\left(\partial_{u}\mathbf{E}_{1}\cdot\mathbf{E}_{1}\times\mathbf{E}_{2}\right)\left(\partial_{v}\mathbf{E}_{2}\cdot\mathbf{E}_{1}\times\mathbf{E}_{2}\right)-\left(\partial_{v}\mathbf{E}_{1}\cdot\mathbf{E}_{1}\times\mathbf{E}_{2}\right)^{2}}{a^{2}}\\
 & =\frac{\left(\partial_{u}\mathbf{E}_{1}\cdot\mathbf{E}_{1}\times\mathbf{E}_{2}\right)\left(\partial_{v}\mathbf{E}_{2}\cdot\mathbf{E}_{1}\times\mathbf{E}_{2}\right)-\left(\partial_{v}\mathbf{E}_{1}\cdot\mathbf{E}_{1}\times\mathbf{E}_{2}\right)^{2}}{\left(EG-F^{2}\right)^{2}}
\end{aligned}
\end{equation}
This formula is another confirmation that the Gaussian curvature is
an intrinsic property of the surface.

\subsubsection{Mean Curvature\label{subMeanCurvature}}

$\bullet$ As given earlier, the mean curvature $H$ is the average\footnote{Or the sum depending on the authors although it will not be a mean
anymore.} of the two principal curvatures, that is:
\begin{equation}
H=\frac{\kappa_{1}+\kappa_{2}}{2}
\end{equation}

$\bullet$ The mean curvature $H$ is also given by the following
formula:
\begin{equation}
H=\frac{eG-2fF+gE}{2\left(EG-F^{2}\right)}=\frac{\mathrm{tr}\left(b_{\alpha}^{\beta}\right)}{2}=\frac{b_{\alpha}^{\alpha}}{2}\,\,\,\,\,\,\,\,\,\,\,\,\,\,\,\,\text{(\ensuremath{\alpha,\beta=1,2})}\label{eqHtrW}
\end{equation}
where $E,F,G,e,f,g$ are the coefficients of the first and second
fundamental forms, the indexed $b$ are the coefficients of the surface
mixed curvature tensor, and $\mathrm{tr}$ stands for the trace of
matrix.\footnote{Being half the trace of a tensor establishes the status of $H$ as
an invariant under permissible transformations.}

$\bullet$ The sign of the mean curvature is dependent on the choice
of the direction of the unit vector $\mathbf{n}$ normal to the surface.

$\bullet$ Like the Gaussian curvature, the mean curvature is invariant
under permissible coordinate transformations and representations.

$\bullet$ The mean curvature of a surface at a given point $P$ is
a measure of the rate of change of area of the surface elements in
the neighborhood of $P$.

$\bullet$ Examples of the mean curvature, $H$, for a number of simple
surfaces:

(A) Plane: $H=0$.

(B) Sphere of radius $R$: $\left|H\right|=\frac{1}{R}$.

(C) Torus parameterized by $x=\left(R+r\sin\phi\right)\cos\theta$,
$y=\left(R+r\sin\phi\right)\sin\theta$ and $z=r\cos\phi$: $H=\frac{R+2r\sin\phi}{2r\left(R+r\sin\phi\right)}$.

$\bullet$ For a Monge patch of the form $\mathbf{r}(u,v)=\left(u,v,f(u,v)\right)$,
the mean curvature is given by:

\begin{equation}
H=\frac{\left(1+f_{v}^{2}\right)f_{uu}-2f_{u}f_{v}f_{uv}+\left(1+f_{u}^{2}\right)f_{vv}}{2\left(1+f_{u}^{2}+f_{v}^{2}\right)^{3/2}}
\end{equation}
where the subscripts $u$ and $v$ stand for partial derivatives with
respect to these surface coordinates.

$\bullet$ The mean curvature may be considered as the 2D equivalent
of the geodesic curvature in 1D.\footnote{The equivalence should be obvious since the mean curvature is a measure
for minimizing the surface area while the geodesic curvature is a
measure for minimizing the curve length.} Accordingly, the 2D minimal surfaces (see $\S$ \ref{subMinimalSurface})
correspond to the 1D geodesic curves (see $\S$ \ref{subGeodesicCurves}).

\subsubsection{Theorema Egregium\label{subTheoremaEgregium}}

$\bullet$ The essence of Gauss \textit{Theorema Egregium} or \textit{Remarkable
Theorem} is that the Gaussian curvature $K$ of a surface is an intrinsic
property of the surface and hence it can be expressed as a function
of the coefficients of the first fundamental form and their partial
derivatives only with no involvement of the coefficients of the second
fundamental form. This can be guessed for example from the last part\footnote{In fact even the first part can be used in this argument since $b$
can be expressed purely in terms of $E,F,G$ and their derivatives
as:
\begin{eqnarray}
b & = & eg-f^{2}=F\left(\frac{\partial\Gamma_{22}^{2}}{\partial u}-\frac{\partial\Gamma_{12}^{2}}{\partial v}+\Gamma_{22}^{1}\Gamma_{11}^{2}-\Gamma_{12}^{1}\Gamma_{12}^{2}\right)+\nonumber \\
 &  & E\left(\frac{\partial\Gamma_{22}^{1}}{\partial u}-\frac{\partial\Gamma_{12}^{1}}{\partial v}+\Gamma_{22}^{1}\Gamma_{11}^{1}+\Gamma_{22}^{2}\Gamma_{12}^{1}-\Gamma_{12}^{1}\Gamma_{12}^{1}-\Gamma_{12}^{2}\Gamma_{22}^{1}\right)
\end{eqnarray}
} of Eq. \ref{eqKRa}.\footnote{In several places of these notes we see that some quantities can be
expressed once in terms of the coefficients of the first fundamental
form exclusively and once in terms of expressions involving the coefficients
of the second fundamental form as well. In this regard, a quantity
is classified as intrinsic if it can be expressed as a function of
the first fundamental form only even if it can also be expressed in
terms involving the second fundamental form.}

$\bullet$ The essence of \textit{Theorema Egregium}, as a statement
of the fact that certain types of curvature are intrinsic to the surface,
is contained in several forms and equations; some of which are indicated
in these notes when they occur. For example, Eq. \ref{eqRb2} which
links the surface curvature tensor to the Riemann-Christoffel curvature
tensor (which is an intrinsic property of the surface and is related
to the Gaussian curvature by Eqs. \ref{eqRK1} and \ref{eqRK} for
instance) can be regarded as a statement of \textit{Theorema Egregium}
since it expresses a form of surface curvature in terms of a combination
of purely intrinsic surface parameters.

$\bullet$ An example may be given to demonstrate the significance
of \textit{Theorema Egregium} that is, if a piece of plane is rolled
into a cylinder of radius $R$, then $\kappa_{1},\kappa_{2},H$ will
change from $0,0,0$ to $\frac{1}{R},0,\frac{1}{2R}$.\footnote{For the cylinder, we are assuming a normal unit vector in the inner
direction.} However, as a consequence of \textit{Theorema Egregium,} $K$ will
not change since $K$ is dependent exclusively on the first fundamental
form which is the same for planes and cylinders as indicated previously
in $\S$ \ref{PlaneCylindeHaveSamerFirstForm}.

$\bullet$ According to \textit{Theorema Egregium,} the Gaussian curvature
of a sufficiently smooth surface of class $C^{3}$ at a given point
$P$ can be represented by the following function of the coefficients
of the first fundamental form and their partial derivatives at $P$:
\begin{eqnarray}
K & =\frac{1}{\left(EG-F^{2}\right)^{2}} & \left\{ \left|\begin{array}{ccc}
C & F_{v}-\frac{1}{2}G_{u} & \frac{1}{2}G_{v}\\
\frac{1}{2}E_{u} & E & F\\
F_{u}-\frac{1}{2}E_{v} & F & G
\end{array}\right|-\left|\begin{array}{ccc}
0 & \frac{1}{2}E_{v} & \frac{1}{2}G_{u}\\
\frac{1}{2}E_{v} & E & F\\
\frac{1}{2}G_{u} & F & G
\end{array}\right|\right\} \\
 & =\frac{1}{\left(EG-F^{2}\right)^{2}} & \left\{ \left|\begin{array}{ccc}
C & \left[22,1\right] & \left[22,2\right]\\
\left[11,1\right] & a_{11} & a_{12}\\
\left[11,2\right] & a_{21} & a_{22}
\end{array}\right|-\left|\begin{array}{ccc}
0 & \left[21,1\right] & \left[21,2\right]\\
\left[21,1\right] & a_{11} & a_{12}\\
\left[21,2\right] & a_{21} & a_{22}
\end{array}\right|\right\} \nonumber
\end{eqnarray}
where $C=\frac{1}{2}\left(-E_{vv}+2F_{uv}-G_{uu}\right)$ and the
subscripts $u$ and $v$ stand for the partial derivatives with respect
to these surface coordinates.

\subsubsection{Gauss-Bonnet Theorem\label{subGaussBonnetTheorem} }

$\bullet$ This theorem ties the geometry of surfaces to their topology.
There are several variants of this theorem; some of which are local
while others are global. Due to the importance and subtlety of this
theorem we give different variants of the theorem and several examples
for both plane and curved surfaces.

$\bullet$ According to the Gauss-Bonnet theorem, if $\mathfrak{S}$
is a simply connected region on a surface of class $C^{3}$ where
$\mathfrak{S}$ is bordered by a finite number of piecewise regular
curves $C_{i}$ that meet in corners\footnote{These corners can be defined as the points of discontinuity of the
tangents of the boundary curves. The angles of these corners are therefore
defined as the angles between the tangent vectors at the points of
discontinuity when traversing the boundary curves in a predefined
sense. Hence, these angles are exterior to the curves and region.
Sometimes, artificial ``corners'' at regular points are introduced,
for convenience to establish an argument, in which case the exterior
angle is zero.} then we have:\footnote{The geodesic and Gaussian curvatures in this theorem should be continuous
and finite over their domain.}
\begin{equation}
\sum_{s}\int_{C_{i}}\kappa_{g}+\sum_{c}\phi_{j}+\iint_{\mathfrak{S}}Kd\sigma=2\pi\label{eqGaussBonnet1}
\end{equation}
where the first sum is over the sides while the second sum is over
the corners, $\kappa_{g}$ is the geodesic curvature of the curves
$C_{i}$ as a function of their coordinates, $\phi_{j}$ are the exterior
angles of the corners and $K$ is the Gaussian curvature of $\mathfrak{S}$
as a function of the coordinates over $\mathfrak{S}$.\footnote{This form of the Gauss-Bonnet theorem may be labeled as a local variant
of the theorem although its locality may not be obvious. However,
these labels are not of crucial importance.}

$\bullet$ As indicated previously, the term $\iint_{\mathfrak{S}}Kd\sigma$,
which is the area integral of the Gaussian curvature over the region
$\mathfrak{S}$ of the above-described surface, is called the total
curvature $K_{t}$ of $\mathfrak{S}$.

$\bullet$ Some examples for the application of the Gauss-Bonnet theorem
are given below:

(A) A disc in a plane with radius $R$ where Eq. \ref{eqGaussBonnet1}
becomes $\frac{1}{R}2\pi R+0+0=2\pi+0+0=2\pi$ which is an identity.

(B) A semi-circular disc in a plane with radius $R$ where Eq. \ref{eqGaussBonnet1}
becomes $\left(\frac{1}{R}\pi R+0\times2R\right)+2\left(\frac{\pi}{2}\right)+0=\pi+\pi+0=2\pi$
which is an identity again.

(C) A spherical triangle on a sphere of radius $R$ whose sides are
two half meridians connecting a pole to the equator and one quarter
of an equatorial parallel and all of its three corners are right angles
($=\frac{\pi}{2}$) where Eq. \ref{eqGaussBonnet1} becomes $\left[0\left(3\times\frac{\pi R}{2}\right)\right]+3\left(\frac{\pi}{2}\right)+\frac{1}{R^{2}}\frac{4\pi R^{2}}{8}=0+\frac{3\pi}{2}+\frac{\pi}{2}=2\pi$.

(D) The upper half of a sphere (or hemisphere in general) of radius
$R$ where Eq. \ref{eqGaussBonnet1} becomes $0\left(2\pi R\right)+0+\frac{1}{R^{2}}2\pi R^{2}=0+0+2\pi=2\pi$.

$\bullet$ The fact that the sum of the interior angles of a planar
triangle is $\pi$ can be regarded as an instance of the Gauss-Bonnet
theorem since for the planar triangle Eq. \ref{eqGaussBonnet1} becomes:
\begin{eqnarray*}
0+\sum_{i=1}^{3}\left(\pi-\theta_{i}\right)+0 & =\\
3\pi-\sum_{i=1}^{3}\theta_{i} & = & 2\pi
\end{eqnarray*}
where $\theta_{i}$ are the interior angles of the triangle and hence
$\sum_{i=1}^{3}\theta_{i}=\pi$. By a similar argument, we can obtain
the sum of the interior angles of a planar polygon of $n$ sides using
the Gauss-Bonnet theorem, that is:
\begin{eqnarray*}
0+\sum_{i=1}^{n}\left(\pi-\theta_{i}\right)+0 & =\\
n\pi-\sum_{i=1}^{n}\theta_{i} & = & 2\pi
\end{eqnarray*}
where $\theta_{i}$ are the interior angles of the $n$-polygon and
hence $\sum_{i=1}^{n}\theta_{i}=\left(n-2\right)\pi$.

$\bullet$ The fact that the perimeter of a planar circle of radius
$R$ is $2\pi R$ can be regarded as another instance of the Gauss-Bonnet
theorem since for such a planar circle Eq. \ref{eqGaussBonnet1} becomes:
\begin{equation}
\frac{1}{R}L+0+0=2\pi
\end{equation}
where $L$ is the length of the circle perimeter and hence $L=2\pi R$.

$\bullet$ \label{LGeodesicTriangle}As a result of the Gauss-Bonnet
theorem, the sum $\theta_{s}$ of the interior angles of a geodesic
triangle\footnote{A geodesic triangle is a triangle with geodesic sides. Also, ``triangle''
here is more general than a three-side planar polygon with three straight
segments as it can be on a curved surface with curved non-planar sides.} on a surface with constant (i.e. have the same sign) Gaussian curvature
is:

(A) $\theta_{s}<\pi$ \textit{iff }$K<0$.

(B) $\theta_{s}=\pi$ \textit{iff }$K=0$.

(C) $\theta_{s}>\pi$ \textit{iff }$K>0$.

This, in the light of Eq. \ref{eqGaussBonnet1}, shows that the total
curvature provides the excess over $\pi$ for the sum when $K>0$
on the surface and the deficit when $K<0$. The vanishing total curvature
in the case of $K=0$ is the intermediate case where the total curvature
has no contribution to the sum.

$\bullet$ As a consequence of the last point, two geodesic curves
on a simply connected patch of a surface with negative Gaussian curvature
cannot intersect at two points because on introducing an artificial
vertex at a regular point on one curve we will have a new corner with
zero exterior angle and hence $\pi$ interior angle. We will then
have a geodesic triangle with $\theta_{s}>\pi$ on a surface over
which $K<0$, in violation of the above-stated condition.

$\bullet$ By a similar argument to the above, the area of a geodesic
polygon\footnote{A geodesic polygon is a polygon with geodesic sides. Again, ``polygon''
here is general and hence it includes that on curved surface with
non-planar curved sides.} on a surface with constant non-zero Gaussian curvature is determined
by the polygon interior angles.

$\bullet$ Because the geodesic curvature is an intrinsic property,
as discussed in $\S$ \ref{subGeodesicCurvature}, the Gauss-Bonnet
theorem is another indication to the fact that the Gaussian curvature
is an intrinsic property and hence it is another demonstration of
\textit{Theorema Egregium} (see $\S$ \ref{subTheoremaEgregium}).

$\bullet$ The Euler characteristic, or the Euler-Poincare characteristic,
is a topological parameter of surfaces which, for polyhedral surfaces,
is given by $\chi=\mathcal{V}+\mathcal{F}-\mathcal{E}$ where $\mathcal{V},\mathcal{F},\mathcal{E}$
are the numbers of vertices, faces and edges of the polyhedron. The
Euler characteristic is also defined for more general surfaces.

$\bullet$ The Euler characteristic of a compact orientable non-polyhedral
surface, like sphere and torus, can be obtained by a polygonal decomposition
by dividing the entire surface into a finite number of non-overlapping
curvilinear polygons which share at most edges or vertices.

$\bullet$ The Gauss-Bonnet theorem has also a global variant which
links the Euler characteristic $\chi$, which is a topological invariant
of the surface, to the Gaussian curvature $K$, which is a geometric
invariant of the surface. This global form of the Gauss-Bonnet theorem
states that: on a compact orientable surface $S$ of class $C^{3}$
these two invariants are linked through the following equation:
\begin{equation}
\iint_{S}Kd\sigma=2\pi\chi\label{eqGaussBonnetGlobal}
\end{equation}
Now, since $\chi$ is a topological invariant of the surface, Eq.
\ref{eqGaussBonnetGlobal} reveals that the total curvature is also
a topological invariant of the surface.

$\bullet$ The global Gauss-Bonnet theorem can be used to determine
the total curvature $K_{t}$ of a surface. For example, the Euler
characteristic of a sphere is 2 and hence from Eq. \ref{eqGaussBonnetGlobal}
its total curvature is $K_{t}=4\pi$ with no need for evaluation of
the area integral. Similarly, the Euler characteristic of a torus
is 0 and hence it can be concluded immediately that its total curvature
is $K_{t}=0$ with no need for evaluating the integral.

$\bullet$ The Euler characteristics of the sphere and torus in the
above examples can be obtained easily by a polygonal decomposition,
as described above. For example, the Euler characteristic of the sphere
can be calculated by dividing the surface of the sphere to 4 curved
polygonal faces with 4 vertices and 6 edges and hence the Euler characteristic
is $\chi=4+4-6=2$, as given in the last point.

$\bullet$ The global Gauss-Bonnet theorem can be used in the opposite
direction, that is it may be used for determining the Euler characteristic
of a surface knowing its Gaussian, and hence total, curvature although
this may be practically of little use. For instance, the Gaussian
curvature of a sphere of radius $R$ is $\frac{1}{R^{2}}$ at every
point on the sphere and hence its total curvature is $K\sigma=\frac{1}{R^{2}}4\pi R^{2}=4\pi$,
therefore from Eq. \ref{eqGaussBonnetGlobal} its Euler characteristic
is $\chi=\frac{4\pi}{2\pi}=2$.

$\bullet$ The Gauss-Bonnet theorem can also be used to find the total
curvature of a surface which is topologically equivalent (technically
described as homeomorphic) to another surface with known total curvature.
For example, the ellipsoid is homeomorphic to the sphere and hence
they have the same Euler characteristic, therefore they have the same
total curvature, according to Eq. \ref{eqGaussBonnetGlobal}, which
is $4\pi$ as known from the aforementioned sphere example.

$\bullet$ The total curvature $K_{t}$ of a surface with a complex
shape can be obtained from the Gauss-Bonnet theorem by reducing the
surface to a topologically-equivalent simpler surface whose total
curvature is known or can be computed more easily. For example, for
an orientable surface of genus\footnote{In simple terms, the genus (in topology) of a surface is the number
of handles or topological holes on the surface.} $\mathfrak{g}$ the Euler characteristic is $\chi=2\left(1-\mathfrak{g}\right)$
and hence its total curvature is given by $K_{t}=2\pi\chi=4\pi\left(1-\mathfrak{g}\right)$.
So, for a compact orientable complexly-shaped surface which is topologically-equivalent
to a sphere with 2 handles the total curvature is $K_{t}=-4\pi$.
Similarly, the genus of a torus is $\mathfrak{g}=1$ and hence its
total curvature is $K_{t}=4\pi\left(1-1\right)=0$, as found earlier
by another method.

$\bullet$ The obvious implication of the global variant of the Gauss-Bonnet
theorem is that the total curvature of a closed surface is dependent
on its genus and not on its shape and hence it is a topological parameter
of the surface as indicated previously.

\subsubsection{Local Shape of Surface\label{subLocalShape}}

$\bullet$ Using the principal curvatures, $\kappa_{1}$ and $\kappa_{2}$,
a point $P$ on a surface is classified according to the shape of
the surface in the close proximity of $P$ as:

(A) Flat when $\kappa_{1}=\kappa_{2}=0$, and hence $K=H=0$.

(B) Elliptic when either $\kappa_{1}>0$ and $\kappa_{2}>0$ or $\kappa_{1}<0$
and $\kappa_{2}<0$, and hence $K>0$.

(C) Hyperbolic when $\kappa_{1}>0$ and $\kappa_{2}<0$, and hence
$K<0$.

(D) Parabolic when either $\kappa_{1}=0$ and $\kappa_{2}\ne0$ or
$\kappa_{2}=0$ and $\kappa_{1}\ne0$, and hence $K=0$ and $H\ne0$.

These constraints on $K$ and $H$ are sufficient and necessary conditions
for determining the type of the surface point as described above.

$\bullet$ The following are some examples for the above classification:

(A) The points of a plane are flat.

(B) The points of an ellipsoid are elliptic.

(C) The points of a catenoid are hyperbolic.

(D) The points of a cone (excluding the apex) and the points of a
cylinder are parabolic.

$\bullet$ Surfaces normally contain points of different shapes. For
example, the torus has elliptic points on its outside rim, parabolic
points on its top and bottom parallels\footnote{These parallels correspond to the two circles contacting its two tangent
planes at the top and bottom which are perpendicular to its axis of
symmetry.}, and hyperbolic points on its inside rim.\footnote{In fact it is elliptic on the outside half and hyperbolic on the inside
half.} However, there are some types of surfaces whose all points are of
the same shape; e.g. all points of planes are flat, all points of
spheres are elliptic, all points of catenoids are hyperbolic, and
all points of cylinders are parabolic.

$\bullet$ The above classification can also be based on the coefficients
of the second fundamental form of the surface where:\footnote{We have $eg-f^{2}=b$ where $b$ is the determinant of the surface
covariant curvature tensor.}

(A) $eg-f^{2}=0$ and $e=f=g=0$ for flat points.

(B) $eg-f^{2}=0$ and $e^{2}+f^{2}+g^{2}\ne0$ for parabolic points.

(C) $eg-f^{2}>0$ for elliptic points.

(D) $eg-f^{2}<0$ for hyperbolic points.

The reader is referred to $\S$ \ref{subNormalCurvature} for the
significance of the sign of $b$.

$\bullet$ The above classification of the shape of a surface in the
immediate neighborhood of a point on a surface (i.e. being flat, elliptic,
hyperbolic or parabolic) is an invariant property with respect to
permissible coordinate transformations. This can be concluded from
the dependence of the classification on the sign of $b$ as explained
above, plus Eq. \ref{eqbJ2b} where the Jacobian (which is real) is
squared and hence the sign of $b$ and $\bar{b}$ is the same. The
classification is also independent of the representation and parameterization
of the surface since these types are real geometric properties of
the surface in their local definitions.

$\bullet$ The invariance of the shape type of the surface points,
as explained in the previous bullet point, holds true even for the
transformations that reverse the direction of the vector $\mathbf{n}$
normal to the surface, because the classification depends on the Gaussian
curvature which is invariant even under this type of transformations
(refer to $\S$ \ref{subGaussianCurvature}). As for the distinction
between the flat and parabolic points which involves $H$ as well,
the distinction is not affected since it depends on the magnitude
of $H$ (i.e. being zero or not) and not on its sign and the magnitude
is not affected by such transformations.

$\bullet$ In the immediate neighborhood of an elliptic point $P$
of a surface $S$, the surface lies completely on one side of the
tangent plane to $S$ at $P$, while at a hyperbolic point the tangent
plane cuts through $S$ and hence some parts of $S$ are on one side
while other parts are on the other side. In the neighborhood of a
parabolic point, the surface lies entirely on one side except for
some points on a curve which lies in the tangent plane itself. As
for planar points, the neighborhood of the point lies in the tangent
plane.

$\bullet$ The surface points can also be classified according to
the geometric shape of Dupin indicatrix (refer to $\S$ \ref{subDupinIndicatrix})
as follow:

(A) If $eg-f^{2}>0$ then $\kappa_{1}$ and $\kappa_{2}$ have the
same sign; hence the point is elliptic and the Dupin indicatrix is
an ellipse or circle.

(B) If $eg-f^{2}=0$ and $e^{2}+f^{2}+g^{2}>0$, then either $\kappa_{1}=0$
and $\kappa_{2}\ne0$ or $\kappa_{2}=0$ and $\kappa_{1}\ne0$; hence
the point is parabolic and the Dupin indicatrix becomes two parallel
lines.\footnote{These lines are characterized by having identically vanishing normal
curvature.}

(C) If $eg-f^{2}<0$ then $\kappa_{1}$ and $\kappa_{2}$ have opposite
signs; hence the point is hyperbolic and the Dupin indicatrix becomes
two conjugate hyperbolas.\footnote{The normal curvature is positive along one of these hyperbolas and
negative along the other, while along the common asymptotes of these
hyperbolas the normal curvature is zero.}

(D) If $e=f=g=0$, and hence $eg-f^{2}=0$, then the point is flat
and the Dupin indicatrix is not defined.\footnote{Hence, having undefined Dupin indicatrix is a characteristic for planar
points.}

Consequently, because of these correlations between the type of point
and its Dupin indicatrix, the Dupin indicatrix can be used to classify
the point as flat, parabolic, elliptic, or hyperbolic.

$\bullet$ In the immediate neighborhood of a point on a surface,
the surface may be approximated by:

(A) An elliptic paraboloid for an elliptic point.

(B) A hyperbolic paraboloid for a hyperbolic point.

(C) A parabolic cylinder for a parabolic point.

(D) A plane for a flat point.

$\bullet$ In the neighborhood of a parabolic point $P$ on a surface
$S$, the tangent plane of $S$ at $P$ meets $S$ in a single line
passing through $P$.

$\bullet$ In the neighborhood of a hyperbolic point $P$ on a surface
$S$, the tangent plane meets $S$ in two lines intersecting at $P$
where these two lines divide $S$ alternatively into regions above
the tangent plane and regions below the tangent plane.

$\bullet$ The following function:
\begin{equation}
\frac{II_{S}}{2}=\frac{e\,dudu+2f\,dudv+g\,dvdv}{2}
\end{equation}
evaluated at a given point $P$ of a class $C^{2}$ surface is called
the osculating paraboloid of $P$. This function is used to determine
the shape of the surface at $P$.

\subsubsection{Umbilical Point\label{subUmbilicPoints}}

$\bullet$ A point on a surface is called ``umbilical'' or ``umbilic''
or ``navel'' if all the normal sections of the surface at the point
have the same curvature $\kappa$.\footnote{For normal sections, the normal curvature is equal to the curvature.}

$\bullet$ The following are some examples of umbilical points on
common surfaces:

(A) All points of planes are umbilical.\footnote{Some authors impose the condition $K>0$ at umbilical points and hence
the points of plane are not umbilical according to these authors.}

(A) All points of spheres are umbilical.\footnote{Hence, umbilical points may be called spherical points.}

(B) The vertex of a paraboloid of revolution is an umbilical point.

(C) The two vertices of an ellipsoid of revolution are umbilical points.

$\bullet$ If all points of a surface of class $C^{3}$ are umbilical
then the surface must be a sphere (the plane is a special case of
sphere as it can be considered a sphere with infinite radius).

$\bullet$ The sufficient and necessary condition for a point to be
umbilical is that the curvature tensor $b_{\alpha\beta}$ is proportional
to the metric tensor $a_{\alpha\beta}$, that is:
\begin{equation}
b_{\alpha\beta}=\lambda(u^{1},u^{2})a_{\alpha\beta}\,\,\,\,\,\,\,\,\,\,\,\,\,\,\,\alpha,\beta=1,2
\end{equation}
where $\lambda$ is independent of the direction of the tangent to
the normal section at the umbilical point. As a result, the determinants
of the two tensors, $a$ and $b$, satisfy the relation:
\begin{equation}
b=\lambda^{2}a\label{eqblam2a}
\end{equation}

$\bullet$ Since the first fundamental form is positive definite,
$a>0$. Therefore, if at the umbilical point $\lambda=0$ then $b=0$
according to Eq. \ref{eqblam2a} and the point is a flat or parabolic
umbilic; otherwise $b>0$ (since $\lambda$ is real) and the point
is an elliptic umbilic.\footnote{By definition, a hyperbolic point cannot be umbilical point.}

$\bullet$ On a plane surface all points are planar umbilic, while
on a sphere all points are elliptic umbilic.\footnote{They may also be called spherical umbilic.}

$\bullet$ At umbilical points $\kappa_{1}=\kappa_{2}$ and hence
we have:
\begin{equation}
K=\kappa_{1}\kappa_{2}=\kappa_{1}\kappa_{1}=\left(\kappa_{1}\right)^{2}=\left(\frac{2\kappa_{1}}{2}\right)^{2}=\left(\frac{\kappa_{1}+\kappa_{2}}{2}\right)^{2}=H^{2}\label{eqKH2}
\end{equation}
where $K$ and $H$ are the Gaussian and mean curvatures at the point.
This can also be obtained from Eq. \ref{eqkapHK} where the discriminant
of this quadratic equation is zero, i.e. $4H^{2}-4K=0$, since at
an umbilical point the two roots are equal and hence $H^{2}=K$.\footnote{The provision $H^{2}=K$ is a sufficient and necessary condition for
the point to be umbilic.}

$\bullet$ The relation between $K$ and $H$ at umbilical points,
as expressed in Eq. \ref{eqKH2}, may be stated by some authors in
the following disguised form:
\begin{equation}
\left(a^{\alpha\beta}b_{\alpha\beta}\right)^{2}=\frac{4}{a}\left(b_{11}b_{22}-b_{12}^{2}\right)\,\,\,\,\,\,\,\,\,\,\,\,\,\,\,\,\,\,\,(\alpha,\beta=1,2)
\end{equation}

\pagebreak{}

\section{Special Types of Curve}

$\bullet$ There are many classifications of space curves depending
on their properties and relations; a few of these categories are given
in the following subsections.

\subsection{Straight Lines}

$\bullet$ A necessary and sufficient condition for a curve of class
$C^{2}$ to be a straight line is that its curvature is zero at every
point on the curve.

$\bullet$ A criterion for a curve to be a straight line is that all
the tangents of the curve are parallel.\footnote{``Parallel'' here is in its absolute Euclidean sense.}

$\bullet$ Another criterion for a curve $C(t):I\rightarrow\mathbb{R}^{3}$
where $t\in I\subseteq\mathbb{R}$ to be a straight line is that for
all the points $t$ in the domain of the curve, $\dot{\mathbf{r}}$
and $\ddot{\mathbf{r}}$ are linearly dependent where the overdot
represents derivative with respect to the general parameter $t$ of
the curve.

$\bullet$ A straight line lying in a surface has the same tangent
plane at each of its points, and hence the line is contained in this
tangent plane.

\subsection{Plane Curves}

$\bullet$ A curve is described as plane curve if the whole curve
can be contained in a plane with no distortion.

$\bullet$ A necessary and sufficient condition for a curve parameterized
by a general parameter $t$ to be a plane curve is that the relation
$\mathbf{\dot{r}}\cdot\left(\mathbf{\ddot{r}}\times\mathbf{\dddot{r}}\right)=0$
holds identically where the overdots represent differentiation with
respect to $t$.

$\bullet$ Plane curves are characterized by having identically vanishing
torsion.\footnote{An identically vanishing torsion is a necessary and sufficient condition
for a regular curve of class $C^{2}$ to be plane curve.}

$\bullet$ For plane curves, the osculating plane at each regular
point contains the entire curve. Therefore, the plane curve may be
characterized by having a common intersection point for all of its
osculating planes.\footnote{It also implies that the curve have the same osculating plane at all
of its points.}

$\bullet$ Two curves are plane curves if they have the same binormal
lines at each pair of corresponding points on the two curves.

$\bullet$ The locus of the centers of curvature of a curve $C$ is
an evolute of $C$ \textit{iff} $C$ is a plane curve.

$\bullet$ On a smooth surface, a geodesic curve (see $\S$ \ref{subGeodesicCurves})
which is also a line of curvature (see $\S$ \ref{subLineofCurvature})
is a plane curve.

\subsection{Involutes and Evolutes}

$\bullet$ If $C_{e}$ is a space curve with a tangent surface $S_{T}$
(see $\S$ \ref{subTangentSurfaces}) and $C_{i}$ is a curve embedded
in $S_{T}$ and it is orthogonal to all the tangent lines of $C_{e}$
at their intersection points, then $C_{i}$ is called an involute
of $C_{e}$ while $C_{e}$ is called an evolute of $C_{i}$. Hence,
the involute is an orthogonal trajectory of the generators of the
tangent surface of its evolute.

$\bullet$ Accordingly, the equation of an involute $C_{i}$ to a
curve $C_{e}$ is given by:
\begin{equation}
\mathbf{r}_{i}=\mathbf{r}_{e}+\left(c-s\right)\mathbf{T}_{e}\label{eqInvolute}
\end{equation}
where $\mathbf{r}_{i}$ is an arbitrary point on the involute, $\mathbf{r}_{e}$
is the point on the curve $C_{e}$ corresponding to $\mathbf{r}_{i}$,
$c$ is a given constant, $s$ is a natural parameter of $C_{e}$
and $\mathbf{T}_{e}$ is the unit vector tangent to $C_{e}$ at $\mathbf{r}_{e}$.

$\bullet$ A visual demonstration of how to generate an involute $C_{i}$
of a curve $C_{e}$, when $\left(c-s\right)$ in Eq. \ref{eqInvolute}
is positive, may be given by unrolling a taut string wrapped around
$C_{e}$ where the string is kept in the tangent direction as it is
unrolled. A fixed point $P$ on the string, where the distance between
$P$ and the point of contact of the string with $C_{e}$ represents
a natural parameter of $C_{e}$, then traces an involute of $C_{e}$.

$\bullet$ A curve has infinitely many involutes corresponding to
different values of $c$ and/or $s$ in Eq. \ref{eqInvolute}.

$\bullet$ For any tangent of a given curve, the length of the line
segment confined between two given involutes is constant which is
the difference between the two $c$ in Eq. \ref{eqInvolute} of the
two involutes.

$\bullet$ If $C_{e}$ is an evolute of $C_{i}$, then for a given
point $P_{e}$ on $C_{e}$ and the corresponding point $P_{i}$ on
$C_{i}$ the principal normal line of $C_{e}$ at $P_{e}$ is parallel
to the tangent line of $C_{i}$ at $P_{i}$.

$\bullet$ A curve $C_{i}$ is a plane curve \textit{iff} the locus
of the centers of curvature of $C_{i}$ is an evolute of $C_{i}$.

$\bullet$ The involutes of a circle are congruent.

$\bullet$ The evolutes of plane curves are helices.

\subsection{Bertrand Curves}

$\bullet$ Bertrand curves are two associated space curves with common
principal normal lines at their continuously-varying one-to-one corresponding
points.\footnote{Some authors characterize Bertrand curves by having a common principal
normal vector.}

$\bullet$ Associated Bertrand curves are characterized by the following
properties:

(A) The product of the torsions of their corresponding points is constant;
i.e. $\tau_{1}(s_{o})\tau_{2}(s_{o})=\mathrm{constant}$ where $\tau_{1}$
and $\tau_{2}$ are the torsions of the two curves and $s_{o}$ is
a given value of their common parameter.

(B) The distance between their corresponding points is constant.

(C) The angle between their corresponding tangent lines is constant.

$\bullet$ For a plane curve $C_{1}$, there is always a curve $C_{2}$
such that $C_{1}$ and $C_{2}$ are associated Bertrand curves.\footnote{This may also be stated as: a plane curve $C_{1}$ has always a Bertrand
curve associate $C_{2}$ (i.e. $C_{1}$ and $C_{2}$ are associated
Bertrand curves).}

$\bullet$ If $C_{1}$ is a curve with non-vanishing torsion such
that $C_{1}$ has more than one Bertrand curve associate, then $C_{1}$
is a circular helix. The reverse is also true.

$\bullet$ If $C_{1}$ is a curve with non-vanishing torsion then
a necessary and sufficient condition for $C_{1}$ to be a Bertrand
curve (i.e. it possesses an associate curve $C_{2}$ such that $C_{1}$
and $C_{2}$ are Bertrand curves) is that there are two constants
$c_{1}$ and $c_{2}$ such that:
\begin{equation}
\kappa=c_{1}\tau+c_{2}
\end{equation}
where $\kappa$ and $\tau$ are the curvature and torsion of the curve
$C_{1}$.

$\bullet$ If a plane curve $C$ has two involutes $C_{1}$ and $C_{2}$,
then $C_{1}$ and $C_{2}$ are Bertrand curves.

\subsection{Spherical Indicatrix\label{subSphericalIndicatrix}}

$\bullet$ A spherical indicatrix of a continuously-varying unit vector
is a continuous curve $\bar{C}$ on the origin-based unit sphere generated
by mapping the unit vector (e.g. $\mathbf{T}$ or $\mathbf{N}$ or
$\mathbf{B}$) of a particular space curve $C$ on an equal unit vector
represented by a point on the origin-based unit sphere. Hence, we
have $\bar{C_{\mathbf{T}}}$, $\bar{C_{\mathbf{N}}}$ and $\bar{C_{\mathbf{B}}}$
as the spherical indicatrices of $C$ corresponding respectively to
the tangent, principal normal and binormal vectors of $C$.

$\bullet$ When $C$ is a naturally parameterized curve then $\bar{\mathbf{r}}(t)=\mathbf{T}(t)$
where $\bar{\mathbf{r}}$ is the origin-based position vector of $\bar{C_{\mathbf{T}}}$.

$\bullet$ A necessary and sufficient condition for $\bar{\mathbf{r}}(t)$
to be a natural parameterization of $\bar{C_{\mathbf{T}}}$ is that
$\kappa(t)=1$ identically.

$\bullet$ The tangent to the curve $\bar{C_{\mathbf{T}}}$ of a curve
$C$ is parallel to the normal vector $\mathbf{N}$ of $C$ at the
corresponding points of the two curves.

$\bullet$ The curvature of the curve $\bar{C_{\mathbf{T}}}$ of a
curve $C$ is related to the curvature and torsion of $C$ by:
\begin{equation}
\kappa_{\mathbf{T}}^{2}=\frac{\kappa^{2}+\tau^{2}}{\kappa^{2}}
\end{equation}
where $\kappa_{\mathbf{T}}$ is the curvature of $\bar{C_{\mathbf{T}}}$
while $\kappa$ and $\tau$ are the curvature and torsion of $C$
respectively.

$\bullet$ The torsion of the curve $\bar{C_{\mathbf{T}}}$ of a curve
$C$ is given by:
\begin{equation}
\tau_{\mathbf{T}}=\frac{\kappa'\tau-\kappa\tau'}{\kappa\left(\kappa^{2}+\tau^{2}\right)}
\end{equation}
where $\tau_{\mathbf{T}}$ is the torsion of $\bar{C_{\mathbf{T}}}$
, $\kappa$ and $\tau$ are the curvature and torsion of $C$ respectively,
and the prime stands for the derivative with respect to the natural
parameter $s$.

$\bullet$ The curvature of the curve $\bar{C_{\mathbf{B}}}$ of a
curve $C$ is given by:
\begin{equation}
\kappa_{\mathbf{B}}=\frac{\kappa^{2}+\tau^{2}}{\kappa^{2}}
\end{equation}
where $\kappa_{\mathbf{B}}$ is the curvature of $\bar{C_{\mathbf{B}}}$
while $\kappa$ and $\tau$ are the curvature and torsion of $C$
respectively.

$\bullet$ The torsion of the curve $\bar{C_{\mathbf{B}}}$ of a curve
$C$ is given by:
\begin{equation}
\tau_{\mathbf{B}}=\frac{\kappa'\tau-\kappa\tau'}{\tau\left(\kappa^{2}+\tau^{2}\right)}
\end{equation}
where $\tau_{\mathbf{B}}$ is the torsion of $\bar{C_{\mathbf{B}}}$
, $\kappa$ and $\tau$ are the curvature and torsion of $C$ respectively,
and the prime stands for the derivative with respect to the natural
parameter $s$.

$\bullet$ The necessary and sufficient condition for the curve $\bar{C_{\mathbf{T}}}$
of a curve $C$ to be a circle is that $C$ is a helix.

$\bullet$ The tangent to the curve $\bar{C_{\mathbf{T}}}$ of a curve
$C$ is parallel to the tangent to the $\bar{C_{\mathbf{B}}}$ of
$C$ at the corresponding points of the two curves.

\subsection{Spherical Curves}

$\bullet$ Spherical curve is a curve that lies completely on the
surface of a sphere. Spherical indicatrices are common examples of
spherical curves (see $\S$ \ref{subSphericalIndicatrix}).

$\bullet$ Circles are the only spherical curves with constant curvature.

$\bullet$ At all points of a spherical curve, the normal plane passes
through the center of the embedding sphere. Conversely, if all the
normal planes of a curve meet in a common point, then the curve is
spherical with the common point being the center of the sphere that
envelopes the curve.

$\bullet$ The sufficient and necessary condition that should be satisfied
by a spherical curve is given by:
\begin{equation}
\frac{R_{\kappa}}{R_{\tau}}+\frac{d}{ds}\left(R_{\tau}\frac{dR_{\kappa}}{ds}\right)=0
\end{equation}
where $R_{\kappa}$ and $R_{\tau}$ are the radii of curvature and
torsion and $s$ is a natural parameter of the curve.

$\bullet$ The center of curvature of a twisted spherical curve $C$
at a given point $P$ on $C$ is the perpendicular projection of the
center of the enveloping sphere on the osculating plane of $C$ at
$P$.

\subsection{Geodesic Curves\label{subGeodesicCurves}}

$\bullet$ The characteristic feature of a geodesic curve is that
it has vanishing geodesic curvature $\kappa_{g}$ at every point on
the curve. This is a necessary and sufficient condition for a curve
to be geodesic.

$\bullet$ In more technical terms, let $S:\Omega\rightarrow\mathbb{R}^{3}$
be a surface defined on a set $\Omega\subseteq\mathbb{R}^{2}$ and
let $C(t):I\rightarrow\mathbb{R}^{3}$, where $I\subseteq\mathbb{R}$,
be a regular curve on $S$, then $C$ is a geodesic curve \textit{iff}
$\kappa_{g}(t)=0$ on all points $t\in I$ in its domain.

$\bullet$ The path of the shortest distance connecting two points
in a Riemannian space is a geodesic. The length of arc, as given by
Eq. \ref{eqArcLength}, is used in the definition of geodesic in this
sense. The geodesic is a straight line in a Euclidean space, but it
is a generalized curved path in a general Riemannian space.

$\bullet$ Examples of geodesic curves on simple surfaces are arcs
of great circles on spheres and arcs of helices on cylinders.\footnote{The generating straight lines and the circles on cylinders are also
geodesics (they may be considered as degenerate helices).} The meridians of surfaces of revolution are also geodesics. The arcs
of parallel circles on a surface of revolution corresponding to stationary
points on the generating curve of the surface are also geodesic curves.

$\bullet$ All straight lines on any surface are geodesic curves.\footnote{For plane surfaces in particular, being a straight line on a plane
is a sufficient and necessary condition for being geodesic.}

$\bullet$ If a geodesic surface curve is not a straight line then
its principal normal vector $\mathbf{N}$ is collinear with the normal
vector $\mathbf{n}$ to the surface at each point on the curve with
non-vanishing curvature.\footnote{Collinearity is equivalent to the condition that $\mathbf{n}$ lies
in the osculating plane of the curve at the given point.} The opposite is also true.

$\bullet$ Being an arc of a great circle on a sphere is a sufficient
and necessary condition for being geodesic.

$\bullet$ Intrinsically, the geodesic curves are straight lines in
the sense that a 2D inhabitant will see them straight since he cannot
measure their curvature. Any deviation from such ``straight lines''
within the surface is therefore a geodesic curvature and can be detected
intrinsically by a 2D inhabitant.

$\bullet$ Although a geodesic curve is usually the curve of the shortest
distance between two points on the surface it is not necessarily so.
For instance, the largest of the two arcs forming a great circle on
a sphere is a geodesic curve but it is not the curve of the shortest
distance on the sphere between its two end points; in fact it is the
curve of the longest distance.\footnote{A constraint my be imposed to make the geodesic minimal by stating
that geodesics minimize distance locally but not necessarily globally
where an infinitesimal element of arc length is considered in this
constraint.}

$\bullet$ If on a surface $S$ there is exactly one geodesic curve
connecting two given points, $P_{1}$ and $P_{2}$, then the length
of the geodesic segment between $P_{1}$ and $P_{2}$ is the shortest
distance on $S$ between these points.

$\bullet$ Being a shortest path is a sufficient but not necessary
condition for being a geodesic, that is all shortest paths connecting
two given points are geodesics but not all geodesics are shortest
paths, as explained in the previous points.

$\bullet$ A sufficient and necessary condition for a curve to be
a geodesic is that the first variation\footnote{The first variation of a functional $F(x)$ may be defined by the
Gateaux derivative of the functional as:
\begin{equation}
\delta F(x,h)=\lim_{\alpha\to0}\frac{F(x+\alpha h)-F(x)}{\alpha}
\end{equation}
where $x$ and $h$ are variable functions and $\alpha$ is a scalar
parameter (refer to the Calculus of Variations for details).} of its length is zero.\footnote{In fact this may be taken as the basis for the definition of geodesic
as the curve connecting two fixed points, $P_{1}$ and $P_{2}$, whose
length possesses a stationary value with regard to small variations,
that is: $\delta\int_{P_{1}}^{P_{2}}ds=0$.}

$\bullet$ It can be shown that a geodesic curve satisfies the Euler-Lagrange
variational principle which is a necessary and sufficient condition
for extremizing the arc length.\footnote{This principle in its generic, simple and most common form is given
by: $\frac{\partial f}{\partial y}-\frac{d}{dx}\left(\frac{\partial f}{\partial y'}\right)=0$.
The reader is advised to consult textbooks on the Calculus of Variations
for details. }

$\bullet$ A physical interpretation may also be given to the geodesic
curve that a free particle restricted to move on the surface will
follow a geodesic path.

$\bullet$ The geodesic, even in its restricted sense as the curve
of the shortest distance, is not necessarily unique; for example all
semicircular meridians of longitude connecting the two poles (or any
antipodal points) of a sphere are geodesics and there is an infinite
number of them.

$\bullet$ In fact even the existence, let alone uniqueness, of a
geodesic connecting two points on a surface is not guaranteed. An
example is the $xy$ plane excluding the origin with two points on
a straight line lying in the plane and passing through the origin
where there is a lower limit for the length of any curve connecting
the two points. This limit is the straight segment connecting the
two points but this segment cannot be a geodesic since it includes
the origin which is not on the surface. Any curve $C$ (other than
the straight segment) connecting the two points cannot be a geodesic
since there is always another curve on the plane connecting the two
points which is shorter than $C$.

$\bullet$ In the neighborhood of a given point $P$ on a surface
and for any specified direction, there is exactly one geodesic curve
that passes through $P$ in that direction. More technically, for
any specific point $P$ on a surface $S$ of class $C^{3}$, and for
any tangent vector $\mathbf{v}$ in the tangent space of $S$ at $P$,
there exists a geodesic curve on the surface in the direction of $\mathbf{v}$
in the immediate neighborhood of $P$.\footnote{This is based on the existence of a unique solution to the geodesic
differential equations (Eq. \ref{eqGeodesicCurve2}) for given initial
values for a point on the curve and its derivative (which represents
its tangent direction) at that point.}

$\bullet$ The obvious examples of the last point is the plane, where
a straight line passes through any point and in any direction, and
the sphere where a great circle passes through any point and in any
direction. A less obvious example is the cylinder where a helix (including
the straight line generators and the circles which can be regarded
as degenerate forms of helix) passes through any point and in any
direction.

$\bullet$ Similarly, there is exactly one geodesic passing through
two sufficiently-close points on a smooth surface.

$\bullet$ Geodesics in curved spaces are the equivalent of straight
lines in flat spaces.

$\bullet$ For planes (or in fact any Euclidean $n$D manifold) there
is a unique geodesic passing between any two points (close or not)
which is the straight line segment connecting the two points.

$\bullet$ The necessary and sufficient conditions that should be
satisfied by a naturally-parameterized curve on a surface, both of
class $C^{2}$, to be a geodesic curve are given by the following
second order non-linear differential equations:\footnote{Due to their non-linearity, closed form explicit solutions of these
equations are rare.}
\begin{equation}
\begin{aligned}\frac{d^{2}u^{1}}{ds^{2}}+\Gamma_{11}^{1}\left(\frac{du^{1}}{ds}\right)^{2}+2\Gamma_{12}^{1}\frac{du^{1}}{ds}\frac{du^{2}}{ds}+\Gamma_{22}^{1}\left(\frac{du^{2}}{ds}\right)^{2} & =0\\
\frac{d^{2}u^{2}}{ds^{2}}+\Gamma_{11}^{2}\left(\frac{du^{1}}{ds}\right)^{2}+2\Gamma_{12}^{2}\frac{du^{1}}{ds}\frac{du^{2}}{ds}+\Gamma_{22}^{2}\left(\frac{du^{2}}{ds}\right)^{2} & =0
\end{aligned}
\label{eqGeodesicCurve2}
\end{equation}
where $s$ is the arc length, and the Christoffel symbols are derived
from the surface metric. The last equations can be merged in a single
equation using tensor notation:
\begin{equation}
\frac{\delta}{\delta s}\left(\frac{du^{\alpha}}{ds}\right)\equiv\frac{d^{2}u^{\alpha}}{ds^{2}}+\Gamma_{\beta\gamma}^{\alpha}\frac{du^{\beta}}{ds}\frac{du^{\gamma}}{ds}=0\label{eqGeodesicCurve}
\end{equation}
where $\alpha,\beta,\gamma=1,2$. These equations, which can be obtained
from Eq. \ref{eqKapgu} by setting the two components of the geodesic
curvature to zero, have no closed form solutions in general. Similar
equations are used to find the geodesic curves between two given points
in general $n$D spaces.

$\bullet$ From Eq. \ref{eqGeodesicCurve}, it can be seen that being
a geodesic is an intrinsic property since it depends exclusively on
the Christoffel symbols which depend only on the coefficients of the
first fundamental form and their partial derivatives.

$\bullet$ From Eq. \ref{eqGeodesicCurve}, it can be seen that for
planes (or indeed for any Euclidean $n$D manifold) the geodesic is
a straight line since in this case the Christoffel symbols vanish
identically and Eq. \ref{eqGeodesicCurve} is reduced to $\frac{d^{2}u^{\alpha}}{ds^{2}}=0$
which has a straight line solution.

$\bullet$ From Eq. \ref{eqkapgE1}, we see that the $u^{1}$ coordinate
curves on a sufficiently smooth surface are geodesics \textit{iff}
$\Gamma_{11}^{2}=0$. Similarly, from Eq. \ref{eqkapgE2}, we see
that the $u^{2}$ coordinate curves are geodesics \textit{iff} $\Gamma_{22}^{1}=0$.

$\bullet$ We also see from Eqs. \ref{eqkapgE1b} and \ref{eqkapgE2b}
that for orthogonal coordinate systems the coordinate curves are geodesics
\textit{iff} $E$ is independent of $v$ and $G$ is independent of
$u$.

$\bullet$ For a Monge patch of the form $\mathbf{r}(u,v)=\left(u,v,f(u,v)\right)$,
the geodesic differential equations are given by:
\begin{equation}
\begin{aligned}\left(1+f_{u}^{2}+f_{v}^{2}\right)u''+f_{u}f_{uu}(u')^{2}+2f_{u}f_{uv}u'v'+f_{u}f_{vv}(v')^{2} & =0\\
\left(1+f_{u}^{2}+f_{v}^{2}\right)v''+f_{v}f_{uu}(u')^{2}+2f_{v}f_{uv}u'v'+f_{v}f_{vv}(v')^{2} & =0
\end{aligned}
\end{equation}
where the subscripts $u$ and $v$ represent partial derivatives with
respect to the surface curvilinear coordinates $u$ and $v$, and
the prime represents derivatives with respect to a natural parameter.

$\bullet$ Each one of the following provisions is a necessary and
sufficient condition for a curve $C$ on a surface $S$ to be a geodesic
curve:

(A) The geodesic component of the curvature vector is zero at each
point on the curve, that is $\mathbf{K}_{g}=\mathbf{0}$ identically.

(B) The osculating plane of the curve at each point of the curve is
orthogonal to the tangent plane of $S$ at that point.\footnote{On geodesic curves $\kappa_{g}=0$ and hence $\mathbf{n}$ and $\mathbf{N}$
are in the same orientation (parallel or anti-parallel) on all points
along the curve (refer to Eqs. \ref{eqkapN} and \ref{eqKnKg}) and
hence $\mathbf{n}$ lies in the osculating plane.}

(C) The normal vector $\mathbf{n}$ to the surface at any point on
the curve lies in the osculating plane.\footnote{This is because for a geodesic curve, $\mathbf{K}_{g}$ vanishes identically
and hence $\mathbf{K}=\kappa\mathbf{N}=\kappa_{n}\mathbf{n}=\mathbf{K}_{n}$.}

(D) The principal normal vector $\mathbf{N}$ of $C$ is normal to
the surface at each point on $C$.

(E) The curvature vector $\mathbf{K}$ is normal to the tangent plane
at each point on the curve.

$\bullet$ Being a geodesic is independent of the choice of the coordinate
system and hence it is invariant under permissible transformations.
It is also independent of the type of representation and parameterization
and hence it is invariant in this sense.\footnote{This may be regarded as another reason for non-uniqueness of geodesic
because even if the trace of the curve is unique, its parameterization
is not unique in general due to different mappings.}

$\bullet$ Geodesics can be open or closed curves and may be self-intersecting.\footnote{From the previous points, we see that ``geodesic curves'' have two
common uses: (a) curves with identically vanishing geodesic curvature
and (b) curves of shortest distance between two points, where (b)
is a subset of (a). Here, ``geodesics'' is used mainly in the first
sense. }

$\bullet$ As a result of the Gauss-Bonnet theorem (see $\S$ \ref{subGaussBonnetTheorem}),
a surface with negative Gaussian curvature cannot have a geodesic
that intersects itself.\footnote{Because the first term on the LHS of Eq. \ref{eqGaussBonnet1} will
vanish and the third term will be negative and hence the sum of the
angles (which is the angle of intersection) in the second term will
not satisfy the relation unless it is greater than $2\pi$ which is
not possible since $K<0$ over the surface (see $\S$ \ref{LGeodesicTriangle}). } Also on such a surface, two geodesics cannot intersect at more than
one point if the geodesics enclose a simply-connected region.\footnote{Similarly, on introducing an artificial vertex at a regular point
on one of these curves we will have a new corner with $\pi$ interior
angle and hence the sum of the geodesic triangle will exceed $\pi$.
Also, on introducing an artificial vertex at a regular point on each
one of these curves we will have a curvilinear quadrilateral whose
internal angles sum is greater than $2\pi$ on a surface with $K<0$
which is not possible (see $\S$ \ref{LGeodesicTriangle}).}

$\bullet$ The lines of curvature (see $\S$ \ref{subLineofCurvature})
are geodesics.

$\bullet$ On a patch on a surface of class $C^{2}$ with orthogonal
coordinate curves and with first fundamental form coefficients being
dependent on only the $u$-coordinate variable\footnote{That is: $E=E(u)$, $F=0$ and $G=G(u)$. The case of dependence on
only the $v$-coordinate variable can be obtained by re-labeling the
coordinate variables and coefficients.}, the following statements apply:

(A) The $u$-coordinate curves with constant $v$ are geodesics at
a point \textit{iff} $\partial_{u}G=0$.

(B) The $v$-coordinate curves with constant $u$ are geodesics.

(C) A curve $C$ represented by $\mathbf{r}=\mathbf{r}\left(u,v(u)\right)$
is a geodesic \textit{iff}:
\begin{equation}
v=\pm\int_{C}\frac{\alpha\sqrt{E}}{\sqrt{G\left(G-\alpha^{2}\right)}}du
\end{equation}
where $\alpha$ is a constant.

$\bullet$ A vector attained by parallel propagation (see $\S$ \ref{subParallelism})
of a tangent vector to a geodesic curve stays always tangent to the
geodesic curve.

$\bullet$ A sufficient and necessary condition for a surface curve
to be geodesic is being a tangent to a parallel vector field.

$\bullet$ As a consequence of the previous points, a vector field
attained by parallel propagation along a geodesic makes a constant
angle with the geodesic.

$\bullet$ Geodesic in curved spaces is a generalization of straight
lines in flat spaces.\footnote{Hence, geodesics may be described as the straightest curves in the
space.}

$\bullet$ Geodesic curves on a developable surface become straight
lines when the surface is developed into a plane by unrolling.

$\bullet$ Isometric surfaces possess identical geodesic equations.

$\bullet$ On a surface with orthogonal coordinate curves, the curves
of constant $u^{\alpha}$ are geodesics \textit{iff} $a_{\beta\beta}$
($\beta\ne\alpha$) is a function of $u^{\beta}$ only.

\subsection{Line of Curvature\label{subLineofCurvature}}

$\bullet$ A ``line of curvature'' is a curve $C$ on a surface
$S$ defined on an interval $I\subseteq\mathbb{R}$ as $C:I\rightarrow S$
with the condition that the tangent of $C$ at each point on $C$
is collinear with one of the principal directions (see $\S$ \ref{subPrincipalCurvatures})
of the surface at that point.

$\bullet$ Since the definition of the line of curvature is based
on the existence of distinct principal directions, the line of curvature
should not include umbilical (including flat) points (see $\S$ \ref{subUmbilicPoints})
due to the absence of distinct principal directions at these points.

$\bullet$ Referring to Eq. \ref{eqtaugphi}, on a line of curvature
either $\sin\theta=0$ or $\cos\theta=0$ and hence the lines of curvature
are characterized by having identically vanishing geodesic torsion
(i.e. $\tau_{g}=0$).

$\bullet$ The lines of intersection of each pair of a triply orthogonal
system are lines of curvature.\footnote{Three families of surfaces in a subset $V$ of a 3D space form a triply
orthogonal system if at each point $P$ of $V$ there is a single
surface of each family passing through $P$ such that each pair of
these surfaces intersect orthogonally at their curve of intersection.}

$\bullet$ Examples of lines of curvature are meridians and parallels
of surfaces of revolution of class $C^{2}$.

$\bullet$ At a non-umbilical point $P$ on a sufficiently smooth
surface $S$, the $u$ and $v$ coordinate curves are aligned with
the principal directions \textit{iff} $f=F=0$ at $P$. Hence, the
coordinate curves on $S$, excluding the umbilical points, are lines
of curvature \textit{iff} $f=F=0$ over the entire surface.

$\bullet$ When the $u$ and $v$ coordinate curves of a surface patch
are lines of curvature, the principal curvatures over the entire patch
will be given by:
\begin{equation}
\kappa_{1}=\frac{e}{E}\,\,\,\,\,\,\,\,\,\,\,\,\,\,\,\,\,\,\,\,\,\,\,\,\,\kappa_{2}=\frac{g}{G}
\end{equation}
where $E,G,e,g$ are the coefficients of the first and second fundamental
forms of the patch. These coefficients are functions of position in
general.

$\bullet$ On a surface of class $C^{3}$, there are two perpendicular
families of lines of curvature in the neighborhood of any non-umbilical
point.

$\bullet$ If the curve of intersection of two surfaces is a line
of curvature for one surface then it is a line of curvature for the
other surface when the two surfaces are intersecting each other at
a constant angle.

$\bullet$ The lines of curvature satisfy the following relation:
\begin{equation}
\left(a_{12}b_{11}-a_{11}b_{12}\right)du^{1}du^{1}+\left(a_{22}b_{11}-a_{11}b_{22}\right)du^{1}du^{2}+\left(a_{22}b_{12}-a_{12}b_{22}\right)du^{2}du^{2}=0
\end{equation}
where the indexed $a$ and $b$ are the coefficients of the surface
covariant metric and curvature tensors respectively.

$\bullet$ The condition that should be satisfied by a line of curvature
on a surface may be given in tensor notation by:
\begin{equation}
\underline{\epsilon}^{\gamma\delta}a_{\alpha\gamma}b_{\beta\delta}du^{\alpha}du^{\beta}=0
\end{equation}

$\bullet$ The generators of developable surfaces are lines of curvature.

$\bullet$ The lines of curvature form a real orthogonal grid over
the surface. On surfaces with constant Gaussian curvature, the lines
of curvature form an isometric conjugate grid.

$\bullet$ For a developable surface, the lines of curvature consist
of its generators and their orthogonal trajectories.

$\bullet$ On a smooth surface, excluding planes and spheres, if the
lines of curvature are selected as the net of coordinate curves then
$a_{12}=b_{12}=0$ over the entire surface.

$\bullet$ On a sufficiently smooth surface, any geodesic which is
a plane curve is a line of curvature.

$\bullet$ A curve is a line of curvature \textit{iff} the tangent
to the curve and the tangent to its spherical image (see $\S$ \ref{subGaussMapping})
at corresponding points are parallel.

$\bullet$ The lines of curvature on a surface, which is not a sphere
or minimal surface, are represented by an orthogonal net on its spherical
image.

$\bullet$ In the neighborhood of a non-umbilical point on a sufficiently
smooth surface there are two orthogonal families of lines of curvature.
Hence, at each point $P$ on such a surface a coordinate patch including
$P$ can be introduced in the neighborhood of $P$ where the coordinate
curves at $P$ are aligned with the principal directions.

$\bullet$ On a surface patch where the Gaussian curvature does not
vanish, the angles between the asymptotic lines (see $\S$ \ref{subAsymptoticDirections})
are bisected by the lines of curvature.

\subsection{Asymptotic Lines\label{subAsymptoticDirections}}

$\bullet$ An asymptotic direction of a surface at a point $P$ is
a direction for which the normal curvature vanishes, i.e. $\kappa_{n}=0$.\footnote{Asymptotic directions are defined only at points for which $K\le0$.}
Hence, at an asymptotic point we have:
\begin{equation}
\mathbf{K}=\mathbf{K}_{g}=\kappa_{g}\mathbf{u}
\end{equation}

$\bullet$ As a consequence of Eq. \ref{eqKapnBA}, $\kappa_{n}$
is zero for directions for which the second fundamental form is zero.
Hence the necessary and sufficient condition for the asymptotic directions
is that:
\begin{equation}
b_{\alpha\beta}du^{\alpha}du^{\beta}=b_{11}(du^{1})^{2}+2b_{12}\,du^{1}du^{2}+b_{22}(du^{2})^{2}=0\label{eqb0}
\end{equation}

$\bullet$ The number of asymptotic directions at elliptic, parabolic
and hyperbolic points is 0, 1 and 2 respectively, while at flat points
all directions are asymptotic. The two asymptotic directions of a
hyperbolic point separate the directions of positive normal curvature
from the directions of negative normal curvature. The sign of the
normal curvature at elliptic and parabolic points is the same in all
directions.\footnote{The directions here are those of the tangents to the normal sections
at the given point (excluding any asymptotic direction).}

$\bullet$ A $t$-parameterized surface curve $C(t):I\rightarrow S$,
where $I\subseteq\mathbb{R}$ is an open interval over which $K<0$
and $S$ represents the surface, is described as an asymptotic line
if at each point $t\in I$ the vector $\mathbf{T}$, which is the
tangent to the curve, is collinear with one of the asymptotic directions
at that point on the curve.

$\bullet$ From the previous points, it can be seen that the asymptotic
lines are characterized by the following:

(A) The normal component of the curvature vector is zero at each point
on the curve, that is $\mathbf{K}_{n}=\mathbf{0}$ identically.

(B) The tangent plane to the surface at each point of the curve coincides
with the osculating plane of the curve at that point.

$\bullet$ The differential equations representing asymptotic lines
can be obtained from the condition that the normal curvature vanishes
identically over the line, that is (see Eq. \ref{eqKapnBA}):
\begin{equation}
e\left(\frac{du^{1}}{ds}\right)^{2}+2f\frac{du^{1}}{ds}\frac{du^{2}}{ds}+g\left(\frac{du^{2}}{ds}\right)^{2}=0\label{eqAsymptotic}
\end{equation}

$\bullet$ From Eq. \ref{eqAsymptotic}, it can be seen that the necessary
and sufficient condition for the $u^{1}$ and $u^{2}$ coordinate
curves to become asymptotic lines is that $e=g=0$ on every point
on the curve.\footnote{This is because on a coordinate curve either $du^{1}=0$ and $du^{2}\ne0$
or $du^{1}\ne0$ and $du^{2}=0$ and hence the middle term is vanishing
anyway.}

$\bullet$ On a sufficiently smooth surface of class $C^{3}$, there
are two distinct families of asymptotic lines in the neighborhood
of a hyperbolic point.

$\bullet$ According to Eq. \ref{eqKapnKapg} $\kappa_{n}=\mathbf{n}\cdot\mathbf{K}$
where $\mathbf{n}$ and $\mathbf{K}$ are respectively the normal
vector to the surface and the curve curvature vector. Hence, a curve
on a sufficiently smooth surface is an asymptotic line \textit{iff}
$\mathbf{n}\cdot\mathbf{K}=0$ identically. This condition is realized
if at each point on the curve either $\mathbf{K}=\mathbf{0}$ or $\mathbf{K}$
and $\mathbf{n}$ are orthogonal vectors.\footnote{The vector $\mathbf{n}$ cannot vanish on a regular point.}
In the former case the point is an inflection point while in the latter
case the osculating plane is tangent to the surface at the point.
Therefore, all points on an asymptotic line should be one of these
types or the other. The reverse is also true, i.e. a curve whose all
points are one of these types or the other is an asymptotic line.

$\bullet$ As a results of the previous point, any straight line on
a surface is an asymptotic line since the curve curvature vector $\mathbf{K}$
vanishes identically on such a line.

$\bullet$ According to the theorem of Beltrami-Enneper, along an
asymptotic non-straight line on a sufficiently smooth surface the
square of the torsion $\tau$ is equal to the negative of the Gaussian
curvature $K$, that is:
\begin{equation}
\tau^{2}=-K
\end{equation}
where $\tau$ and $K$ are evaluated at each individual point along
the curve.\footnote{Since asymptotic directions are defined only at points for which $K\le0$,
the square of the torsion in the above equation is equal to the absolute
value of the Gaussian curvature at the point, that is: $\tau^{2}=\left|K\right|$.}

$\bullet$ The torsions of two asymptotic lines passing through a
given point on a surface are equal in magnitude and opposite in sign.

$\bullet$ As we will see (refer to $\S$ \ref{subConjugateDirections}),
asymptotic directions are self-conjugate.\footnote{Some authors take self-conjugation as the defining characteristic
for being asymptotic.}

$\bullet$ From the definition of the asymptotic direction plus the
Euler equation (Eq. \ref{eqEuler}), we see that the angle $\theta$
which an asymptotic direction makes with the principal direction of
$\kappa_{1}$ at a given point $P$ on a sufficiently smooth surface
$S$ is given by:
\begin{equation}
\tan^{2}\theta=-\frac{\kappa_{1}}{\kappa_{2}}
\end{equation}
where $\kappa_{1}$ and $\kappa_{2}$ are the principal curvatures
of $S$ at $P$.

$\bullet$ On a sufficiently smooth surface with orthogonal families
of asymptotic lines the mean curvature $H$ is zero.

$\bullet$ The principal directions at a given point on a smooth surface
bisect the asymptotic directions at the point.

$\bullet$ Eq. \ref{eqb0} is quadratic and hence it possesses two
solutions which are real and distinct, or real and coincident, or
conjugate imaginary depending on its discriminant $\Delta$ which
is opposite in sign to the determinant $b$ of the surface covariant
curvature tensor.\footnote{The discriminant is: $\Delta=4(b_{12})^{2}-4b_{11}b_{22}$ while the
determinant is: $b=b_{11}b_{22}-(b_{12})^{2}$, and hence $\Delta=-4b$.} Hence, the asymptotic directions at a given point on a surface can
be classified according to the determinant $b$ at the point as:

(A) Real and distinct for $\Delta>0$ and hence $b<0$.

(B) Real and coincident for $\Delta=0$ and hence $b=0$.

(C) Conjugate imaginary for $\Delta<0$ and hence $b>0$.

$\bullet$ From Eq. \ref{eqKRa}, the sign of the Gaussian curvature
$K$ is the same as the sign of $b$, since $a>0$. Hence, the classification
in the previous point can also be based on $K$ as stated for $b$
in the previous point.

$\bullet$ A straight line contained in a surface is an asymptotic
line. This is due to the fact that such a line is wholly contained
in a plane, which is the tangent space of each of its points, and
hence the normal curvature vanishes identically along the line.

$\bullet$ As a consequence of having identically vanishing normal
curvature, the osculating planes at each point of a curved asymptotic
line on a surface are tangent to the surface.

$\bullet$ On a sufficiently smooth surface, if a geodesic curve $C$
is a line of curvature then $C$ is a plane curve.

\subsection{Conjugate Directions\label{subConjugateDirections}}

$\bullet$ A direction $\frac{\delta u}{\delta v}$ at a point on
a sufficiently smooth surface is described as conjugate to the direction
$\frac{du}{dv}$ if the following relation holds:
\begin{equation}
d\mathbf{r}\cdot\delta\mathbf{n}=0
\end{equation}
where $d\mathbf{r}=\mathbf{E}_{1}du+\mathbf{E}_{2}dv$ and $\delta\mathbf{n}=\partial_{u}\mathbf{n}\delta u+\partial_{v}\mathbf{n}\delta v$.
Due to the symmetry, $\frac{du}{dv}$ is also conjugate to $\frac{\delta u}{\delta v}$,
and hence the two directions are described as conjugate directions.

$\bullet$ Two families of curves on a sufficiently smooth surface
are described as conjugate families if the directions of their tangents
at each point on the curves are conjugate directions.

$\bullet$ The $u$ and $v$ coordinate curves on a smooth surface
are conjugate families of curves \textit{iff} $f$, which is the coefficient
of the second fundamental form, vanishes identically.

$\bullet$ At a hyperbolic or elliptic point on a sufficiently smooth
surface, each direction has a unique conjugate direction.

$\bullet$ An asymptotic direction is self-conjugate direction.

\pagebreak{}

\section{Special Types of Surface\label{subTypesSurfaces}}

$\bullet$ There are many classifications of surfaces in 3D spaces
depending on their properties and relations; a few of these classifications
are given in the following subsections.

\subsection{Plane Surfaces}

$\bullet$ Planes are simple, ruled, connected, elementary surfaces.

$\bullet$ All the coefficients of the surface curvature tensor vanish
identically throughout plane surfaces.

$\bullet$ The Riemann-Christoffel curvature tensor vanishes identically
over plane surfaces.

$\bullet$ The Gaussian curvature $K$ and the mean curvature $H$
vanish identically over planes.

$\bullet$ Planes are minimal surfaces.

$\bullet$ All points on planes are flat umbilical.

$\bullet$ At any point on a plane surface, $\kappa_{1}=\kappa_{2}=0$
and hence all the directions are principal directions (or there is
no principal direction).

$\bullet$ At any point on a plane surface, all the directions are
asymptotic.

$\bullet$ A sufficient and necessary condition for a surface to be
isometric with the plane is having an identically vanishing Riemann-Christoffel
curvature tensor. The same applies for identically vanishing Gaussian
curvature.

\subsection{Quadratic Surfaces\label{subQuadraticSurfaces}}

$\bullet$ Quadratic surfaces are defined by the following quadratic
equation:\footnote{We assume a Euclidean 3D space with a rectangular Cartesian coordinate
system.}
\begin{equation}
A_{ij}x^{i}x^{j}+B_{i}x^{i}+C=0\,\,\,\,\,\,\,\,\,\,\,\,\,\,\,\,\,\,\,\,\,\text{(\ensuremath{i,j=1,2,3})}
\end{equation}
where the coefficients $A_{ij}$ and $B_{i}$ are real-valued tensors
of rank-2 and rank-1 respectively and $C$ is a real scalar.

$\bullet$ There are six main non-degenerate types of quadratic surfaces:
ellipsoid, hyperboloid of one sheet, hyperboloid of two sheets, elliptic
paraboloid, hyperbolic paraboloid, and quadric cone. By rigid motion
transformations, consisting of translation and rotation of coordinates
whose purpose is to put the center of symmetry of these surfaces at
the origin and orient their axes with the coordinate lines, these
types can be given in the following canonical forms:\footnote{In the following equations, $a,b,c$ are real parameters, and for
convenience we use $x,y,z$ for $x^{1},x^{2},x^{3}$ respectively.
Also for the first three of these surfaces the origin of coordinates
is not a valid point.}

(A) Ellipsoid:
\begin{equation}
\frac{x^{2}}{a^{2}}+\frac{y^{2}}{b^{2}}+\frac{z^{2}}{c^{2}}=1
\end{equation}

(B) Hyperboloid of one sheet:\label{HyperboloidOneSheet}
\begin{equation}
\frac{x^{2}}{a^{2}}+\frac{y^{2}}{b^{2}}-\frac{z^{2}}{c^{2}}=1
\end{equation}

(C) Hyperboloid of two sheets:
\begin{equation}
\frac{x^{2}}{a^{2}}-\frac{y^{2}}{b^{2}}-\frac{z^{2}}{c^{2}}=1
\end{equation}

(D) Elliptic paraboloid:
\begin{equation}
\frac{x^{2}}{a^{2}}+\frac{y^{2}}{b^{2}}-z=0
\end{equation}

(E) Hyperbolic paraboloid:
\begin{equation}
\frac{x^{2}}{a^{2}}-\frac{y^{2}}{b^{2}}-z=0
\end{equation}

(F) Quadric cone:
\begin{equation}
\frac{x^{2}}{a^{2}}+\frac{y^{2}}{b^{2}}-\frac{z^{2}}{c^{2}}=0
\end{equation}

\subsection{Ruled Surfaces\label{subRuledSurfaces}}

$\bullet$ A ``ruled surface'', or ``scroll'', is a surface generated
by a continuous translational-rotational motion of a straight line
in space. Hence, at each point of the surface there is a straight
line passing through the point and lying entirely in the surface.
Planes, cones, cylinders and Mobius strips are common examples of
ruled surface. The different perspectives of the generating line along
its movement are described as the rulings of the surface.

$\bullet$ A ruled surface that can be generated by two different
families of lines is called doubly-ruled surface. Examples of doubly-ruled
surface are hyperbolic paraboloids and hyperboloids of one sheet.

$\bullet$ At any point of a regular ruled surface, the Gaussian curvature
is non-positive ($K\le0$).

$\bullet$ The tangent surface (see $\S$ \ref{subSurfaces}) of a
smooth curve is a ruled surface generated by the tangent line of the
curve.

$\bullet$ The tangent plane is constant along a branch, represented
by the tangent line at a given point, of the tangent surface of a
curve. Examples are the tangent planes of cylinders and cones along
their generators.

$\bullet$ If $P$ is a point on a curve $C$ where $C$ has a tangent
surface $S$, then the tangent plane to $S$ along the ruling that
passes through $P$ coincides with the osculating plane of $C$ at
$P$. Hence, the tangent surface may be described as the envelope
of the osculating planes of the curve.

$\bullet$ The tangent plane to a cylinder or a cone is constant along
their generators.

\subsection{Developable Surfaces }

$\bullet$ As defined previously (see $\S$ \ref{subSurfaces}), a
surface that can be flattened into a plane without local distortion
is called developable surface.

$\bullet$ A developable surface can also be defined as a surface
that is isometric to the Euclidean plane.

$\bullet$ In 3D manifolds, all developable surfaces are ruled surfaces
but not all ruled surfaces are developable surfaces. A ruled surface
is developable if the tangent plane is constant along every ruling
of the surface as it is the case with cones and cylinders.

$\bullet$ The neighborhood of each point on a sufficiently smooth
surface with no flat points is developable \textit{iff} the Gaussian
curvature vanishes identically on the surface.

$\bullet$ The generators of a developable surface and their orthogonal
trajectories are their lines of curvature.

$\bullet$ A developable surface, excluding cylinder and cone, is
a tangent surface of a curve.

$\bullet$ If a developable surface $S$, excluding cylinder and cone,
is rolled out on a plane then all the points of an orthogonal trajectory
of the tangent planes of $S$ will map on a single point on the plane.

$\bullet$ The collection of normal lines to a surface $S$ along
a given curve $C$ on $S$ make a developable surface \textit{iff}
$C$ is a line of curvature (see $\S$ \ref{subLineofCurvature}).

$\bullet$ Intrinsically, any developable surface is equivalent\footnote{Equivalent means having the same metric characteristics.}
to a plane and hence any two developable surfaces are equivalent to
each other.

$\bullet$ The generators of developable surfaces are lines of curvature
(see $\S$ \ref{subLineofCurvature}).

\subsection{Isometric Surfaces\label{subIsometricSurfaces}}

$\bullet$ An isometry is an injective mapping from a surface $S$
to a surface $\bar{S}$ which preserves distances.

$\bullet$ As a consequence of preserving the lengths in isometric
mappings, the angles and areas are also preserved.

$\bullet$ Examples of isometric surfaces are cylinder and cone which
are both isometric to plane.

$\bullet$ Two isometric surfaces, such as a cylinder and a cone or
each one of these and a plane, appear identical to a 2D inhabitant.
Any difference between the two can only be perceived by an external
observer residing in a reference frame in the enveloping space.

$\bullet$ Isometry is an equivalence relation and hence it is reflective,
symmetric and transitive, that is for three surfaces $S_{1}$, $S_{2}$
and $S_{3}$ we have:\footnote{The symbol $\sim$ represents an isometric relation.}

(A) $S_{1}\sim S_{1}$.

(B) $S_{1}\sim S_{2}\,\,\,\,\,\,\,\,\,\Longleftrightarrow\,\,\,\,\,\,\,\,\,\,S_{2}\sim S_{1}$.

(C) $S_{1}\sim S_{2}$ and $S_{2}\sim S_{3}$ $\,\,\,\,\,\,\,\Longrightarrow\,\,\,\,\,\,S_{1}\sim S_{3}$.

$\bullet$ Two isometric surfaces possess identical first fundamental
forms and hence any difference between them, as viewed extrinsically
from the embedding space, is based on the difference between their
second fundamental forms.

$\bullet$ If two surfaces have constant equal Gaussian curvature
then they are isometric. The mapping relation between the two surfaces
then include three constants corresponding to the three independent
coefficients of the first fundamental form.

$\bullet$ A surface of revolution is isometric to itself in infinitely-many
ways, each of which corresponds to a rotation of the surface through
a given angle around its axis of symmetry. Hence, a surface $S_{1}$
which is isometric to a surface of revolution $S_{2}$ is equivalent
to $S_{2}$.

\subsection{Tangent Surfaces\label{subTangentSurfaces}}

$\bullet$ As stated previously, the tangent surface of a space curve
is a surface generated by the assembly of all the tangent lines to
the curve. The tangent lines of the curve are called the generators
of the tangent surface.

$\bullet$ Accordingly, the equation of a tangent surface $S_{T}$
to a curve $C$ is given by:
\begin{equation}
\mathbf{r}_{T}=\mathbf{r}_{i}+k\mathbf{T}_{i}\,\,\,\,\,\,\,\,\,\,\,\,\,\,\,\,\,(-\infty<k<\infty)
\end{equation}
where $\mathbf{r}_{T}$ is an arbitrary point on the tangent surface,
$\mathbf{r}_{i}$ is a given point on the curve $C$, $k$ is a real
variable, and $\mathbf{T}_{i}$ is the unit vector tangent to $C$
at $\mathbf{r}_{i}$. The tangent surface is generated by varying
$i$ along $C$.

$\bullet$ The tangent surface of a curve is made of two parts: one
part corresponding to $k>0$ and the other corresponding to $k<0$
where the curve is a border line between these two parts. The two
parts of the tangent surface are tangent to each other along the curve
which forms a sharp edge between the two.\footnote{The curve is called the edge of regression of the surface.}

$\bullet$ The tangent plane is constant along a branch\footnote{A branch in this context is the tangent line of the curve. It is also
called generator.} of the tangent surface of a curve. This tangent plane is the osculating
plane of the curve at the point of contact of the branch with the
curve.

$\bullet$ According to the definition of involute, all the involutes
of a curve $C_{e}$ are wholly embedded in the tangent surface of
$C_{e}$.

$\bullet$ The normal to the tangent surface of a space curve $C$
at a point of a given ruling $\mathfrak{R}$ is parallel to the binormal
line of $C$ at the point of contact of $C$ with $\mathfrak{R}$.

\subsection{Minimal Surfaces\label{subMinimalSurface}}

$\bullet$ A ``minimal surface'' is a surface whose area is minimum
compared to the area of any other surface sharing the same boundary.
Hence, the minimal surface is an extremal with regard to the integral
of area over its domain.

$\bullet$ A common physical example of a minimal surface is a soap
film formed between two coaxial rings where it takes the minimal surface
shape of a catenoid due to the surface tension. This problem, and
its alike of investigations related to the physical realization of
minimal surfaces, may be described as the Plateau problem.

$\bullet$ Since the mean curvature of a surface at a point $P$ is
a measure of the rate of change of area of the surface elements in
the neighborhood of $P$, a minimal surface is characterized by having
an identically vanishing mean curvature and hence the principal curvatures
at each point have the same magnitude and opposite signs.

$\bullet$ Examples of minimal surface shapes are planes, catenoids,
helicoids and ennepers.\footnote{The catenoid is a surface of revolution generated by revolving a catenary
around its directrix. The helicoid is a ruled surface (see $\S$ \ref{subRuledSurfaces})
with the property that for each point $P$ on the surface there is
a helix passing through $P$ and contained entirely in the surface.
The enneper is a self-intersecting surface which can be described
parametrically in various ways; one of which is:
\begin{eqnarray*}
x & = & -\frac{u^{3}}{3}+u+uv^{2}\\
y & = & -u^{2}v-v+\frac{v^{3}}{3}\\
z & = & u^{2}-v^{2}
\end{eqnarray*}
where $u,v\in(-\infty,\infty)$. It is noteworthy that the catenoid
and helicoid are locally isometric.}

$\bullet$ A minimal surface is characterized by having an orthogonal
net of asymptotic lines and a conjugate net of minimal lines.\footnote{Minimal lines are curves of minimal length. Having an orthogonal net
of asymptotic lines is a sufficient and necessary condition for having
zero mean curvature.}

\pagebreak{}

\section{Tensor Differentiation\label{secTensorDifferentiation}}

$\bullet$ Tensor differentiation, whether covariant or absolute,
over a 2D surface follows similar rules to those stated in \cite{SochiTC1, SochiTC2}
for general $n$D curved spaces. Some of these rules are:

(A) The sum and product rules of differentiation apply to covariant
and absolute differentiation as usual.

(B) The covariant and absolute derivatives of tensors are tensors.

(C) The covariant and absolute derivatives of scalars and invariant
tensors of higher ranks are the same as the ordinary derivatives.

(D) The covariant and absolute derivative operators commute with contraction
of indices.

(E) The covariant and absolute derivatives of the metric, Kronecker
and alternating tensors (and their associated tensors) vanish identically
in any coordinate system, that is:
\begin{equation}
\begin{aligned}a_{\alpha\beta|\gamma} & =0\,\,\,\,\,\,\,\,\,\,\,\,\,\, & a_{\,\,\,|\gamma}^{\alpha\beta} & =0\\
\delta_{\beta|\gamma}^{\alpha} & =0\,\,\,\,\,\,\,\,\,\,\,\,\,\, & \delta_{\beta\omega|\gamma}^{\alpha\delta} & =0\\
\underline{\epsilon}_{\alpha\beta|\gamma} & =0\,\,\,\,\,\,\,\,\,\,\,\,\,\, & \underline{\epsilon}_{\,\,\,|\gamma}^{\alpha\beta} & =0
\end{aligned}
\end{equation}
where the sign $|$ represents covariant or absolute differentiation
with respect to the surface coordinate $u^{\gamma}$. Hence, these
tensors should be treated like constants in tensor differentiation.

$\bullet$ An exception of these rules is the covariant derivative
of the space basis vectors in their covariant and contravariant forms
which is identically zero, as stated previously in \cite{SochiTC2},
that is:
\begin{equation}
\begin{aligned}\mathbf{E}_{i;j} & =\partial_{j}\mathbf{E}_{i}-\Gamma_{ij}^{k}\mathbf{E}_{k}=\Gamma_{ij}^{k}\mathbf{E}_{k}-\Gamma_{ij}^{k}\mathbf{E}_{k}=\mathbf{0}\\
\mathbf{E}_{;j}^{i} & =\partial_{j}\mathbf{E}^{i}+\Gamma_{kj}^{i}\mathbf{E}^{k}=-\Gamma_{kj}^{i}\mathbf{E}^{k}+\Gamma_{kj}^{i}\mathbf{E}^{k}=\mathbf{0}
\end{aligned}
\label{eqEspace}
\end{equation}
but this is not the case with the surface basis vectors in their covariant
and contravariant forms, $\mathbf{E}_{\alpha}$ and $\mathbf{E}^{\alpha}$,
whose covariant derivatives do not vanish identically. The reason
is that, due to curvature, the partial derivatives of the surface
basis vectors do not necessarily lie in the tangent plane and hence
the following relations:
\begin{equation}
\begin{aligned}\partial_{j}\mathbf{E}_{i} & =+\Gamma_{ij}^{k}\mathbf{E}_{k}\\
\partial_{j}\mathbf{E}^{i} & =-\Gamma_{kj}^{i}\mathbf{E}^{k}
\end{aligned}
\end{equation}
which are valid in the enveloping space and are used in Eq. \ref{eqEspace},
are not valid on the surface anymore.

$\bullet$ At a point $P$ on a sufficiently smooth surface with geodesic
surface coordinates and Cartesian rectangular space coordinates, the
covariant and absolute derivatives reduce respectively to the partial
and total derivatives at $P$.

$\bullet$ The covariant derivative of the surface basis vectors is
symmetric in its two indices, that is:
\begin{equation}
\begin{aligned}\mathbf{E}_{\alpha;\beta} & =\partial_{\beta}\mathbf{E}_{\alpha}-\Gamma_{\alpha\beta}^{\gamma}\mathbf{E}_{\gamma}\\
 & =\partial_{\alpha}\mathbf{E}_{\beta}-\Gamma_{\beta\alpha}^{\gamma}\mathbf{E}_{\gamma}\\
 & =\mathbf{E}_{\beta;\alpha}
\end{aligned}
\end{equation}

$\bullet$ The covariant derivative of the surface basis vectors,
$\mathbf{E}_{\alpha;\beta}$, represents space vectors which are normal
to the surface with no tangential component.

$\bullet$ The covariant derivative of a space tensor with respect
to a surface coordinate $u^{\alpha}$ is formed by the inner product
of the covariant derivative of the tensor with respect to the space
coordinates $x^{l}$ by the tensor $x_{\alpha}^{l}$. For example,
the covariant derivative of $A^{i}$ with respect to $u^{\alpha}$
is given by:
\begin{equation}
A_{\,\,;\alpha}^{i}=A_{\,\,;k}^{i}x_{\alpha}^{k}
\end{equation}

$\bullet$ The covariant derivative with respect to a surface coordinate
$u^{\beta}$ of a mixed tensor $A_{\alpha}^{i}$, which is contravariant
with respect to transformations in space coordinates $x^{i}$ and
covariant with respect to transformations in surface coordinates $u^{\alpha}$,
is given by:\footnote{An example of such a tensor is $x_{\alpha}^{i}$ which was discussed
in $\S$ \ref{subSurfaceMetric}.}
\begin{equation}
A_{\alpha;\beta}^{i}=\frac{\partial A_{\alpha}^{i}}{\partial u^{\beta}}+\Gamma_{jk}^{i}A_{\alpha}^{k}\frac{\partial x^{j}}{\partial u^{\beta}}-\Gamma_{\alpha\beta}^{\gamma}A_{\gamma}^{i}
\end{equation}
where the Christoffel symbols with Latin and Greek indices are derived
respectively from the space and surface metrics. This pattern can
be easily generalized to a mixed tensor of type ($m,n$) $A_{\alpha_{1}\ldots\alpha_{n}}^{i_{1}\ldots i_{m}}$
which is contravariant in transformations of space coordinates $x^{i}$
and covariant in transformations of surface coordinates $u^{\alpha}$.
For example, the covariant derivative of a tensor $T_{\alpha\beta}^{ij}$
with respect to $u^{\gamma}$ is given by:
\begin{equation}
A_{\alpha\beta;\gamma}^{ij}=\frac{\partial A_{\alpha\beta}^{ij}}{\partial u^{\gamma}}+\Gamma_{mk}^{i}A_{\alpha\beta}^{mj}\frac{\partial x^{k}}{\partial u^{\gamma}}+\Gamma_{mk}^{j}A_{\alpha\beta}^{im}\frac{\partial x^{k}}{\partial u^{\gamma}}-\Gamma_{\alpha\gamma}^{\delta}A_{\delta\beta}^{ij}-\Gamma_{\beta\gamma}^{\delta}A_{\alpha\delta}^{ij}
\end{equation}

$\bullet$ The above rules can be extended further to include tensors
with space and surface contravariant indices and space and surface
covariant indices. For Example, the covariant derivative of a tensor
$A_{j\beta}^{i\alpha}$ with respect to a surface coordinate $u^{\gamma}$,
where $i$ and $j$ are space indices and $\alpha$ and $\beta$ are
surface indices, is given by:
\begin{equation}
A_{j\beta;\gamma}^{i\alpha}=\frac{\partial A_{j\beta}^{i\alpha}}{\partial u^{\gamma}}+\Gamma_{mk}^{i}A_{j\beta}^{m\alpha}\frac{\partial x^{k}}{\partial u^{\gamma}}+\Gamma_{\delta\gamma}^{\alpha}A_{j\beta}^{i\delta}-\Gamma_{jk}^{m}A_{m\beta}^{i\alpha}\frac{\partial x^{k}}{\partial u^{\gamma}}-\Gamma_{\beta\gamma}^{\delta}A_{j\delta}^{i\alpha}
\end{equation}
This example can be easily extended to the most general form of a
tensor with any combination of covariant and contravariant space and
surface indices.

$\bullet$ The covariant derivative of the surface basis vector $\mathbf{E}_{\alpha}$,
which in tensor notation is denoted by $x_{\alpha}^{i}$, is given
by:
\begin{equation}
x_{\alpha;\beta}^{i}=\frac{\partial^{2}x^{i}}{\partial u^{\beta}\partial u^{\alpha}}+\Gamma_{jk}^{i}x_{\alpha}^{j}x_{\beta}^{k}-\Gamma_{\alpha\beta}^{\delta}x_{\delta}^{i}=x_{\beta;\alpha}^{i}
\end{equation}
where curvilinear space coordinates are in use.

$\bullet$ The mixed second order covariant derivative of the surface
basis vectors is given by:
\begin{equation}
x_{\alpha;\beta\gamma}^{i}=b_{\alpha\beta;\gamma}n^{i}+b_{\alpha\beta}n_{;\gamma}^{i}=b_{\alpha\beta;\gamma}n^{i}-b_{\alpha\beta}a^{\delta\omega}b_{\delta\gamma}x_{\omega}^{i}
\end{equation}
where the covariant derivative of the surface covariant curvature
tensor is given, as usual, by:
\begin{equation}
b_{\alpha\beta;\gamma}=\frac{\partial b_{\alpha\beta}}{\partial u^{\alpha}}-\Gamma_{\alpha\gamma}^{\delta}b_{\delta\beta}-\Gamma_{\beta\gamma}^{\delta}b_{\alpha\delta}
\end{equation}

$\bullet$ The covariant differentiation operators in mixed derivatives
are not commutative and hence for a surface covariant vector $A^{\gamma}$
we have:
\begin{equation}
A_{\,\,;\alpha\beta}^{\gamma}-A_{\,\,;\beta\alpha}^{\gamma}=R_{\,\,\delta\alpha\beta}^{\gamma}A^{\delta}
\end{equation}
where $R_{\,\,\delta\alpha\beta}^{\gamma}$ is the Riemann-Christoffel
curvature tensor of the second kind for the surface.

$\bullet$ The mixed second order covariant derivatives of the surface
basis vectors satisfy the following relation:\footnote{This is an instance of the relation: $A_{j;kl}-A_{j;lk}=R_{jkl}^{i}A_{i}$
which is given and explained in \cite{SochiTC2}.}
\begin{equation}
x_{\alpha;\beta\gamma}^{i}-x_{\alpha;\gamma\beta}^{i}=R_{\,\,\,\alpha\beta\gamma}^{\delta}x_{\delta}^{i}
\end{equation}

$\bullet$ As defined in \cite{SochiTC2}, the absolute or intrinsic
derivative of a tensor field along a $t$-parameterized curve in an
$n$D space with respect to the parameter $t$ is the inner product
of the covariant derivative of the tensor and the tangent vector to
the curve. This identically applies to the absolute derivative of
curves contained in 2D surfaces.

$\bullet$ The absolute derivative of a tensor field along a $t$-parameterized
curve on a surface with respect to the parameter $t$ follows similar
rules to those of a space curve in a general $n$D space, as outlined
in the previous notes \cite{SochiTC2}. Hence, the absolute derivative
of a differentiable surface vector field $\mathbf{A}$ in its covariant
and contravariant forms with respect to the parameter $t$ is given
by:
\begin{eqnarray*}
\frac{\delta A_{\alpha}}{\delta t} & = & \frac{dA_{\alpha}}{dt}-\Gamma_{\alpha\beta}^{\gamma}A_{\gamma}\frac{du^{\beta}}{dt}\\
\frac{\delta A^{\alpha}}{\delta t} & = & \frac{dA^{\alpha}}{dt}+\Gamma_{\beta\gamma}^{\alpha}A^{\gamma}\frac{du^{\beta}}{dt}
\end{eqnarray*}
where the Christoffel symbols are derived from the surface metric.
It should be remarked that if $\mathbf{A}$ is a space vector field
defined along the above surface curve then the above formulae will
take a similar form but with change from surface to space coordinates,
and hence the curve is treated as a space curve, that is:
\begin{eqnarray*}
\frac{\delta A_{i}}{\delta t} & = & \frac{dA_{i}}{dt}-\Gamma_{ik}^{j}A_{j}\frac{dx^{k}}{dt}\\
\frac{\delta A^{i}}{\delta t} & = & \frac{dA^{i}}{dt}+\Gamma_{jk}^{i}A^{j}\frac{dx^{k}}{dt}
\end{eqnarray*}
where the Christoffel symbols are derived from the space metric.

$\bullet$ The absolute derivative of the tensor $A_{\alpha}^{i}$,
defined in the previous points, along a $t$-parameterized surface
curve is given by:
\begin{equation}
\frac{\delta A_{\alpha}^{i}}{\delta t}=A_{\alpha;\beta}^{i}\frac{du^{\beta}}{dt}=\left(\frac{\partial A_{\alpha}^{i}}{\partial u^{\beta}}+\Gamma_{jk}^{i}A_{\alpha}^{k}\frac{\partial x^{j}}{\partial u^{\beta}}-\Gamma_{\alpha\beta}^{\gamma}A_{\gamma}^{i}\right)\frac{du^{\beta}}{dt}
\end{equation}

$\bullet$ To extend the idea of geodesic coordinates to deal with
mixed tensors of the type $A_{\alpha}^{i}$, a rectangular Cartesian
coordinate system over the space and a geodesic system on the surface
can be introduced and hence at the poles the absolute and covariant
derivatives become total and partial derivatives respectively.

$\bullet$ The covariant and absolute derivatives of space and surface
metric, permutation and Kronecker tensors in their covariant, contravariant
and mixed forms vanish identically and hence they behave as constants
with respect to tensor differentiation when involved in inner or outer
product operations with other tensors and commute with these operators.

$\bullet$ The \textit{surface} covariant and absolute derivatives
of \textit{space} metric tensor, space Kronecker tensor, space alternating
tensor and space basis vectors vanish identically, that is:
\begin{equation}
\begin{aligned}g_{ij|\gamma} & =0\,\,\,\,\,\,\,\,\,\,\,\,\,\, & g_{\,\,\,|\gamma}^{ij} & =0\\
\delta_{j|\gamma}^{i} & =0\,\,\,\,\,\,\,\,\,\,\,\,\,\, & \delta_{kl|\gamma}^{ij} & =0\\
\underline{\epsilon}_{ijk|\gamma} & =0\,\,\,\,\,\,\,\,\,\,\,\,\,\, & \underline{\epsilon}_{\,\,\,|\gamma}^{ijk} & =0\\
\mathbf{E}_{i|\gamma} & =\mathbf{0}\,\,\,\,\,\,\,\,\,\,\,\,\,\, & \mathbf{E}_{\,\,|\gamma}^{i} & =\mathbf{0}
\end{aligned}
\end{equation}
where the sign $|$ represents covariant or absolute differentiation
with respect to the surface coordinate $u^{\gamma}$. Hence, these
space tensors are in lieu of constants with respect to surface tensor
differentiation.

$\bullet$ The nabla $\nabla$ based differential operations, such
as gradient and divergence, apply to the surface as for any general
curved space and hence the formulae given in \cite{SochiTC2} can
be used with the substitution of the surface metric parameters. For
example, the divergence of a surface vector field $A^{\alpha}$ is
given by:
\begin{equation}
\nabla\cdot\mathbf{A}=\frac{1}{\sqrt{a}}\partial_{\alpha}\left(\sqrt{a}A^{\alpha}\right)
\end{equation}
and the Laplacian of a surface scalar field $f$ is given by:
\begin{equation}
\nabla^{2}f=\frac{1}{\sqrt{a}}\partial_{\alpha}\left(\sqrt{a}a^{\alpha\beta}\partial_{\beta}f\right)\label{eqLaplacianGeneral}
\end{equation}
where the indexed $a$ is the surface contravariant metric tensor,
$a$ is the determinant of the surface covariant metric tensor and
$\alpha,\beta=1,2$.

\pagebreak{}

\phantomsection
\addcontentsline{toc}{section}{References}
\bibliographystyle{unsrt}

\end{document}